\documentclass{jcmlatex}
\def\JCMvol{xx}
\def\JCMno{x}
\def\JCMyear{200x}
\usepackage{subfigure}
\usepackage{epstopdf}

\usepackage{algorithm}
\usepackage{algorithmic}
\usepackage{amssymb}
\usepackage{amsmath}
\usepackage{amsfonts}
\usepackage{mathrsfs}
\usepackage{mathtools}
\usepackage{bm}
\usepackage{enumerate}
\usepackage{textcomp}
\usepackage{cite}
\usepackage{cleveref}
\usepackage{color}

\newcommand{\F}{\bm F}
\newcommand{\U}{\bm U}
\newcommand{\B}{\bm B}
\newcommand{\K}{\bm K}
\newcommand{\R}{\bm R}

\newcommand{\E}{\bm E}
\newcommand{\G}{\bm G}

\newcommand{\Ep}{\mathcal{E}}
\newcommand{\hf}{\frac{1}{2}}

\newcommand{\dx}{\Delta x}
\newcommand{\dy}{\Delta y}
\newcommand{\dt}{\Delta t}
\newcommand{\jph}{{j+\frac{1}{2}}}
\newcommand{\jmh}{{j-\frac{1}{2}}}
\newcommand{\kph}{{k+\frac{1}{2}}}
\newcommand{\kmh}{{k-\frac{1}{2}}}

\newcommand\eref[1]{(\ref{#1})}

\def\JCMvol{xx}
\def\JCMno{x}
\def\JCMyear{200x}
\def\JCMreceived{xxx}
\def\JCMrevised{xxx}
\def\JCMaccepted{xxx}

\begin{document}

\markboth{Y.-P. QIU, Z. GAO, A. KURGANOV, B.-S. WANG, AND X. WEN}
{Fifth-Order Well-Balanced Path-Conservative A-WENO Scheme for the Ripa Model}

\title{FIFTH-ORDER WELL-BALANCED PATH-CONSERVATIVE A-WENO SCHEME FOR THE RIPA MODEL}

\author{Yan-Ping Qiu
\thanks{School of Mathematical Sciences and Laboratory of Marine Mathematics, Ocean University of China, Qingdao, China\\ Email:
qyp@stu.ouc.edu.cn}
\and
Zhen Gao
\thanks{School of Mathematical Sciences and Laboratory of Marine Mathematics, Ocean University of China, Qingdao, China\\ Email:
zhengao@ouc.edu.cn}
\and
Alexander Kurganov
\thanks{Department of Mathematics and Shenzhen International Center for Mathematics, Southern University of Science and Technology,
Shenzhen, China\\ Email: alexander@sustech.edu.cn}
\and
Bao-Shan Wang
\thanks{School of Mathematical Sciences and Laboratory of Marine Mathematics, Ocean University of China, Qingdao, China\\ Email:
wbs@ouc.edu.cn}
\and
Xiao Wen
\thanks{College of Mathematics and Systems Science, Shandong University of Science and Technology, Qingdao, China\\ Email:
xiaowen@sdust.edu.cn}}

\maketitle

\begin{abstract}
In this work, we introduce a fifth-order well-balanced (WB) path-conservative A-WENO scheme with the central-upwind numerical fluxes
(PCCU-5) for the Ripa model. The proposed scheme is capable of exactly preserving a variety of steady states, including still-water,
moving-water, isobaric, and constant water height ones. This goal is achieved with the help of a flux globalization technique: The source
terms are incorporated into the fluxes, resulting in a quasi-conservative system, for which central-upwind numerical fluxes are computed
using the path-conservative integration. The proposed A-WENO scheme utilizes a WENO interpolation of the equilibrium variables rather than
the conservative ones to ensure the WB property. In addition, we perform the WENO interpolation of the local characteristic equilibrium
variables to mitigate numerical oscillations near discontinuities. We perform a series of numerical experiments, which demonstrate that the
proposed fifth-order WB PCCU-5 scheme achieves high resolution and clearly outperforms its second-order counterpart. Our numerical results
also demonstrate the importance of the local characteristic projection for significantly reducing (eliminating) numerical oscillations near
discontinuities.
\end{abstract}



\section{Introduction}\label{sec1}
We consider the shallow water equations, in which the water temperature fluctuations are taken into account. The studied system was
introduced in \cite{Ripa93,Rip95} for modeling ocean currents and is often referred to as the Ripa model, which in the one-dimensional (1-D)
case reads as
\begin{equation}
\left\{\begin{aligned}
&h_t+q_x=0,\\
&q_t+\Big(hu^2+\hf h^2\theta\Big)_x=-h\theta Z_x,\\
&(h\theta)_t+(hu\theta)_x=0,
\end{aligned}\right.
\label{1.1}
\end{equation}
where $x$ is a spatial coordinate, $t$ is time, $h(x,t)$ is the water depth, $u(x,t)$ is the velocity, $q(x,t)=h(x,t)u(x,t)$ is the
discharge, $Z(x)$ is the bottom topography, and $\theta(x,t)$ is the potential temperature defined to be the reduced gravity
$g\Delta\Theta/{\Theta_{\rm ref}}$, where $\Delta\Theta$ is the difference in potential temperature from a certain reference value
$\Theta_{\rm ref}$ and $g$ is the acceleration due to gravity. Notice that when the potential temperature field $\theta\equiv g$, the system
\eref{1.1} reduces to the classical Saint-Venant system of shallow water equations.

The Ripa model \eref{1.1} admits several steady states, including both {\em still-} and {\em moving-water} equilibria. There are three types
of still-water equilibria:

\noindent
$\bullet$ {\em ``Lake-at-rest''} steady states satisfying
\begin{equation}
u\equiv0,\quad h+Z\equiv{\rm Const},\quad\theta\equiv{\rm Const};
\label{1.2}
\end{equation}
$\bullet$ {\em Isobaric} steady states satisfying
\begin{equation}
u\equiv0,\quad P:=\hf h^2\theta\equiv{\rm Const},\quad Z\equiv{\rm Const},
\label{1.3}
\end{equation}
where $P$ is the pressure;

\noindent
$\bullet$ {\em Constant water depth} steady states satisfying
\begin{equation}
u\equiv0,\quad h\equiv{\rm Const},\quad Z+\hf h\ln\theta\equiv{\rm Const}.
\label{1.4}
\end{equation}
Moving-water steady states satisfy
\begin{equation}
q=hu\equiv{\rm Const},\quad\frac{u^2}{2}+\theta(h+Z)\equiv{\rm Const},\quad\theta\equiv{\rm Const}.
\label{1.5}
\end{equation}
Note that the ``lake-at-rest'' states \eref{1.2} can be viewed as a particular case of \eref{1.5} with $u=0$.

It is well-known that many physically relevant solutions are, in fact, small perturbations of the steady states \eref{1.2}--\eref{1.5}.
Good numerical methods for the Ripa model \eref{1.1} are supposed to accurately capture such solutions on practically feasible and thus
relatively coarse meshes. The way to construct such methods is to enforce them to be able to exactly preserve (some of the) steady states.
Our goal is to develop a non-oscillatory high-order scheme, which is fully well-balanced (WB) in the sense that it is capable of preserving
all of the steady states \eref{1.2}--\eref{1.5}.

A variety of WB methods for the Ripa model have been proposed within the past decade. Several second-order WB finite-volume (FV) methods
were introduced in \cite{Chertock14,Sanchez16,Touma15}. For higher-order WB weighted essentially non-oscillatory (WENO) schemes and
discontinuous Galerkin (DG) methods, we refer the reader to \cite{Han16,Rehman21} and \cite{Huang23,Qian18}, respectively. The
aforementioned methods can, however, preserve still-water equilibria only. Several numerical methods, which are capable of preserving both
still- and moving-water steady states, have been recently introduced in \cite{Britton20,Dong22,Li20,Thanh21}.

We would like to stress that the second-order FV method presented in \cite{Li20} relies on the flux globalization approach, which we use in
this paper to ensure the WB property. We would also like to emphasize that the second-order FV method presented in \cite{Sanchez16} utilizes
path-conservative technique, which we use in this paper to accurately treat the nonconservative term arising on the right-hand side of
\eref{1.1} in the case of discontinuous bottom topography $Z$.

We aim at developing a high-order flux globalization based fully WB path-conservative central-upwind (PCCU) scheme. Central-upwind (CU)
schemes were originally introduced in \cite{KNP, Kurganov00} for general multidimensional hyperbolic systems of conservation laws, and later
extended to a variety of hyperbolic systems of balance laws, including the classical Saint-Venant system
\cite{Cheng16,Kurganov02,Kurganov07} and several other shallow water models \cite{Chertock14,CKW,KM14,KP09}; see also the review paper
\cite{Kur_Acta}. A general flux globalization approach of development fully WB CU schemes was proposed in \cite{CHO18} (we refer the reader
to \cite{GC01}, where the first flux globalization based FV method was introduced) and then applied to the Saint-Venant system
\cite{Cheng19} and thermal-rotating shallow water (TRSW) equations \cite{Kurganov20}, which can be viewed as a generalization of the Ripa
model. A PCCU scheme was originally introduced in \cite{Castro19} as a general numerical tool for designing robust and quite accurate
numerical methods for general nonconservative hyperbolic systems. In \cite{Kurganov23}, the idea of path-conservative integration was
incorporated into the flux globalization framework, and second-order flux globalization based WB PCCU were introduced. In these schemes, one
first identifies equilibrium variables, uses them to perform piecewise linear reconstruction, recovers the other point values out of the
equilibrium ones, and then evaluates global fluxes in a very accurate way, which takes into account a proper treatment of nonconservative
products. Flux globalization based WB PCCU schemes, which are capable of exactly preserving a wide variety of steady states including
discontinuous ones, were successfully applied to several nonconservative systems, including the nozzle flow system \cite{Kurganov23},
Saint-Venant system with the Manning friction term \cite{Cao22}, shallow water equations in channels \cite{CKN23,KN25a}, shallow water
system with Coriolis forces \cite{Cao22}, two-layer shallow water equations \cite{Kurganov23}, TRSW equations \cite{CKL23}, two-layer TRSW
equations \cite{Cao23b,CKLZ2}, and rotating shallow water magnetohydrodynamics equations \cite{Chertock24}.

A fifth-order extension of the flux globalization based WB PCCU schemes has been recently proposed in \cite{CKX_Beijing} in the framework of
finite-difference (FD) alternative weighted essentially non-oscillatory (A-WENO) schemes (in \cite{CKX_Beijing}, the fifth-order WB PCCU
schemes were applied to the nozzle flow and two-layer shallow water equations). High-order A-WENO schemes, introduced in
\cite{Jiang13,Wang18}, are based on a direct application of ``standard'' FV numerical fluxes together with high-order correction terms.
Compared with the ``clasical'' FD WENO schemes (see, e.g., \cite{ShuActa} and references therein), the A-WENO schemes generate
higher-resolution solutions without introducing spurious oscillations. In the past few years, A-WENO schemes have been successfully applied
to various hyperbolic systems; see, e.g., \cite{Balsara24,Chu22,CKMZ,CKX_Beijing,Fu24,Li20,Liu22,WDGK,Wang22,Wang23,Yang24}.

In this paper, we present a fifth-order flux globalization based fully WB PCCU scheme (which will be referred to as a PCCU-5 scheme) for the
Ripa model. While designing this scheme, we face several challenges. The PCCU-5 scheme should be capable of exactly preserving both
``lake-at-rest'' (satisfying either \eref{1.2}, \eref{1.3}, or \eref{1.4}) and moving-water \eref{1.5} steady states. To achieve this goal,
we perform WENO interpolation for the equilibrium variables rather than the conservative ones. Such interpolation, however, may be
oscillatory, and thus one needs to perform a local characteristic decomposition (LCD) and then to interpolate in the local characteristic
variables. This is done in a WB way using an LCD of equilibrium variables recently introduced in \cite{CKNX}. We test the proposed PCCU-5
scheme on a number of 1-D and two-dimensional (2-D) numerical examples, which demonstrate high accuracy, WB, and essentially non-oscillatory
properties of the PCCU-5 scheme, as well as the enhanced resolution compared with those achieved by the second-order flux globalization WB
PCCU scheme and the PCCU-5 scheme constructed local characteristic decomposition.

The rest of the paper is organized as follows. In \S\ref{sec2}, we describe the 1-D Ripa model and introduce the equilibrium variables. In
\S\ref{sec3}, we design the 1-D PCCU-5 scheme. In \S\ref{sec4}, we describe the 2-D Ripa model, several of its steady states, and the
corresponding sets of equilibrium variables. In \S\ref{sec5}, we extend the 1-D PCCU-5 scheme to the 2-D case. The 1-D and 2-D PCCU-5
schemes are tested on a variety of numerical experiments in \S\ref{sec6}. Finally, brief conclusions are provided in \S\ref{sec7}.

\section{One-Dimensional Ripa Model}\label{sec2}
The 1-D Ripa model \eref{1.1} can be written in the following vector form:
\begin{equation}
\U_t+\F(\U)_x=B(\U)\U_x
\label{2.1}
\end{equation}
with the conservative variables $\U$, flux function $\F(\U)$, and matrix $B(\U)$ given by
\begin{equation}
\U=\begin{pmatrix}h\\q\\h\theta\\Z\end{pmatrix},\quad\F(\U)=\begin{pmatrix}q\\qu+P\\q\theta\\0\end{pmatrix},\quad
B(\U)=\begin{pmatrix}0&0&0&0\\0&0&0&-h\theta\\0&0&0&0\\0&0&0&0\end{pmatrix},
\label{2.2}
\end{equation}
where the equation $Z_t=0$ is added to the Ripa model \eref{1.1}. For \eref{2.1}--\eref{2.2}, the eigenvalues of
${\cal A}{(\U)}:=\frac{\partial\F}{\partial\U}-B(\U)$ are $\lambda_1=u-c<\lambda_2=u<\lambda_3=u+c$ and $\lambda_4=0$, where
$c:=\sqrt{h\theta}$ is the speed of sound.

The Ripa model \eref{2.1}--\eref{2.2} can be recast in the quasi-conservative form
\begin{equation}
\U_t+\K(\U)_x=\bm0,
\label{2.3}
\end{equation}
with the global flux
\begin{equation}
\begin{aligned}
&\K(\U)=\F(\U)-\R(\U),\\
&\R(\U)=\int\limits_{\widehat x}^xB(\U(\xi,t))\U_\xi(\xi,t)\,{\rm d}\xi=
\Big(0,-\int\limits_{\widehat x}^xh(\xi,t)\theta(\xi,t)Z_\xi(\xi)\,{\rm d}\xi,0,0\Big)^\top,
\end{aligned}
\label{2.4}
\end{equation}
where $\widehat x$ is an arbitrary number.

\paragraph{Equilibrium Variables.} We first note that it is easy to verify that the quantity ${\cal E}$, which is defined by
\begin{equation*}
{\cal E}:=\frac{u^2}{2}+\theta(h+Z)+{\cal Q},\quad{\cal Q}:=
-\int\limits_{\widehat x}^x\left[\sqrt{2P(\xi,t)}\big(\sqrt{\theta(\xi,t)}\big)_\xi+Z(\xi)\theta_\xi(\xi,t)\right]{\rm d}\xi,
\end{equation*}
remains constant at all of the steady states \eref{1.2}--\eref{1.5}. We then follow \cite{Kurganov23} and rewrite the global flux derivative
as
\begin{equation*}
\K(\U)_x=M(\U)\E(\U)_x,
\end{equation*}
where $M(\U)$ is a $4\times4$ matrix and $\E(\U)$ is a vector of equilibrium variables given by
\begin{equation*}
M(\U)=\begin{pmatrix}1&0&0&0\\u&h&0&0\\\theta&0&q&0\\0&0&0&0\end{pmatrix},\quad\E(\U)=\begin{pmatrix}q\\\Ep\\\theta\\Z\end{pmatrix}.
\end{equation*}
Notice that we call $\E$ the equilibrium variables since $q$ and $\Ep$ are constant at all of the steady states \eref{1.2}--\eref{1.5},
$\theta$ is constant at steady states \eref{1.2} and \eref{1.5}, and the bottom topography $Z$ is included in $\E$ for a purely technical
reason.

We also remark that the system \eref{2.3} can be rewritten in the equivalent (for smooth solutions) form using the equilibrium variables
$\bm E$ as follows:
\begin{equation*}
\E(\U)_t+C(\U)\E(\U)_x=\bm I(\bm U),
\end{equation*}
where
\begin{equation}
C(\U)=\begin{pmatrix}u&h&0&0\\\theta&u&\dfrac{q}{2}&0\\0&0&u&0\\0&0&0&0\end{pmatrix}
\label{2.8a}
\end{equation}
and $\bm I$ is the following global (integral) term, which vanishes at the steady states \eref{1.2}--\eref{1.5}:
$$
\bm I(\bm U)=\bigg(0,\int\limits_{\widehat x}^x\Big(\hf h(\xi,t)u_\xi(\xi,t)-u(\xi,t)Z_\xi(\xi)\Big)\theta_\xi(\xi,t)\,{\rm d}\xi,0,0
\bigg)^\top.
$$
We remark that the matrix $C$ has a complete eigensystem and it can be diagonalized using
\begin{equation}
Q(\U)=\begin{pmatrix}0&-\dfrac{q}{2\theta}&\dfrac{c}{\theta}&-\dfrac{c}{\theta}\\[1.5ex]0&0&1&1\\0&1&0&0\\1&0&0&0\end{pmatrix},\quad
Q^{-1}(\U)=\begin{pmatrix}0&0&0&1\\0&0&1&0\\[0.5ex]\dfrac{\theta}{2c}&\dfrac{1}{2}&\dfrac{q}{4c}&0\\[2.0ex]
-\dfrac{\theta}{2c}&\dfrac{1}{2}&-\dfrac{q}{4c}&0\end{pmatrix},
\label{2.9}
\end{equation}
so that $Q^{-1}(\U)C(\U)Q(\U)$ is a diagonal matrix.

\section{One-Dimensional Scheme}\label{sec3}
In this section, we present the 1-D scheme. We first introduce a uniform mesh consisting of the cells $I_j:=[x_\jmh,x_\jph]$ of size
$x_\jph-x_\jmh\equiv\dx$ centered at $x_j=\hf(x_\jmh+x_\jph),~j=1,\dots,N$. We then introduce the point values $\U_j(t):\approx\U(x_j,t)$
and evolve them in time using the modified version of the semi-discrete WB A-WENO scheme from \cite{CKX_Beijing}:
\begin{equation}
\frac{{\rm d}\U_j}{{\rm d}t}=-\frac{\bm{{\cal K}}_\jph-\bm{{\cal K}}_\jmh}{\dx},
\label{3.1}
\end{equation}
where the global numerical fluxes $\bm{{\cal K}}_\jph$ are given by
\begin{equation}
\bm{{\cal K}}_\jph=\bm{{\cal K}}_\jph^{\rm FV}-\frac{(\dx)^2}{24}(\K_{xx})_\jph+\frac{7(\dx)^4}{5760}(\K_{xxxx})_\jph,
\label{3.2}
\end{equation}
and $\bm{{\cal K}}_\jph^{\rm FV}$ are the CU numerical fluxes from \cite{Kurganov23}:
\begin{equation}
\bm{{\cal K}}_\jph^{\rm FV}=\frac{a_\jph^+\K_\jph^--a_\jph^-\K_\jph^+}{a_\jph^+-a_\jph^-}+
\frac{a_\jph^+a_\jph^-}{a_\jph^+-a_\jph^-}\left(\widehat{\U}_\jph^+-\widehat{\U}_\jph^-\right),
\label{3.3}
\end{equation}
and $(\K_{xx})_\jph$ and $(\K_{xxxx})_\jph$ are the high-order correction terms.

In fact, Eqs. \eref{3.2}--\eref{3.3} should be used for the first three components of $\bm{{\cal K}}_\jph$ only as the fourth equation in
\eref{2.3}--\eref{2.4} is trivial ($Z_t=0$), and hence we take $\bm{{\cal K}}_\jph^{(4)}\equiv0$ for all $j$.

Notice that all of the indexed quantities in \eref{3.1}--\eref{3.4}, as well as most of the indexed quantities below, are time-dependent,
but we omit this dependence for the sake of brevity.

In \eref{3.3}, the one-sided point values of the global flux function are computed as
\begin{equation}
\K_\jph^\pm=\F_\jph^\pm-\R_\jph^\pm,\quad\F_\jph^\pm:=\F\big(\U_\jph^\pm\big),
\label{3.4}
\end{equation}
where the one-sided values of the global part of the flux, $\R_\jph^\pm$ are determined precisely as in \cite{CKX_Beijing} and
$\U_\jph^\pm$ are obtained with the help of the fifth-order WENO interpolation applied using the LCD of the equilibrium variables as in
\cite{CKNX}.

In order to apply the fifth-order Ai-WENO-Z interpolation to the equilibrium variables $\bm E$, we first compute
\begin{equation*}
{\cal E}_j=\frac{u_j^2}{2}+\theta_j(h_j+Z_j)+{\cal Q}_j,
\end{equation*}
where the values ${\cal Q}_j$ are computed using the path-conservative integration as in \cite{CKL23}.

We first set $\widehat x=x_\hf$, hence ${\cal Q}_\hf^-=0$, and then ${\cal Q}_\hf^+={\cal B}_{\bm\Psi,\hf}$ and the rest of the terms
${\cal Q}_\jph^\pm$ are evaluated recursively by
\begin{equation*}
{\cal Q}_\jph^-={\cal Q}_\jmh^++{\cal B}_j,\quad{\cal Q}_\jph^+={\cal Q}_\jph^-+{\cal B}_{\bm\Psi,\jph},\quad j=1,\ldots,N.
\end{equation*}
Here,
\begin{equation}
{\cal B}_j\approx-\int\limits_{I_j}\left[\sqrt{2P}(\sqrt{\theta})_x+Z\theta_x\right]{\rm d}x
\label{3.5a}
\end{equation}
are computed using the Newton-Cotes quadrature described in \cite{FGKW}. The terms
$$
{\cal B}_{\bm\Psi,\jph}\approx-\int\limits_0^1\bigg[\sqrt{2P\big(\bm\Psi_\jph(s)\big)}\left(\sqrt{\theta\big(\bm\Psi_\jph(s)\big)}\right)'+
Z\big(\bm\Psi_\jph(s)\big)\left(\theta\big(\bm\Psi_\jph(s)\big)\right)'\bigg]{\rm d}s
$$
are obtained by approximating the path integrals using the linear segment paths for $\theta$, $\sqrt{2P}$, and $Z$, which results in
\begin{equation}
\begin{aligned}
{\cal B}_{\bm\Psi,\jph}&=-\hf\left(\sqrt{2P_\jph^-}+\sqrt{2P_\jph^+}\,\right)\left(\sqrt{\theta_\jph^+}-\sqrt{\theta_\jph^-}\,\right)\\
&-\hf\left(Z_\jph^-+Z_\jph^+\right)\left(\theta_\jph^+-\theta_\jph^-\right);
\end{aligned}
\label{3.6a}
\end{equation}
see \cite{CKL23} for details. Next, the terms ${\cal Q}_j$ are computed by
\begin{equation}
{\cal Q}_j\approx{\cal Q}_\jmh^+-\int\limits_{x_\jmh}^{x_j}\left[\sqrt{2P}(\sqrt{\theta})_x+Z\theta_x\right]{\rm d}x,
\label{3.7a}
\end{equation}
where the integral is also approximated by the Newton-Cotes quadrature described in \cite{FGKW}.

It should be pointed out that the point values of $\theta$, $P$, and $Z$ used in \eref{3.6a} as well as in the Newton-Cotes quadratures
applied to evaluate the integrals in \eref{3.5a} and \eref{3.7a}, are obtained by applying the Ai-WENO-Z interpolation directly to the
fields $\theta$, $P$, and $Z$. These values are, however, not WB, and we now proceed with the non-oscillatory WB interpolation performed in
the local characteristic equilibrium variables following the approach recently introduced in \cite{CKNX}.

To this end, we evaluate the matrix $C_\jph:=C\big(\overline{\U}_\jph\big)$, where the matrix $C$ is given in \eref{2.8a} and
$\overline{\U}_\jph:=\hf(\U_j+\U_{j+1})$. We then introduce $Q_\jph:=Q\big(\overline{\U}_\jph\big)$ and
$\big(Q_\jph\big)^{-1}:=Q^{-1}\big(\overline{\U}_\jph\big)$, where $Q$ and $Q^{-1}$ are given in \eref{2.9}, switch to the local
characteristic equilibrium variables in the neighborhood of $x=x_\jph$:
\begin{equation*}
\bm\Gamma_\ell=\big(Q_\jph\big)^{-1}\E_{j+\ell},\quad\ell=-2,\dots,3,
\end{equation*}
apply the Ai-WENO-Z interpolation from \cite{FGKW,ZGW} to each component of $\bm\Gamma_\ell$ to reconstruct the one-sided point values
$\bm\Gamma_\hf^\pm$, and then transform them back to equilibrium variables to obtain
\begin{equation*}
\E_\jph^\pm=Q_\jph\bm\Gamma_\hf^\pm.
\end{equation*}

Equipped with $\E_\jph^\pm$, we now evaluate $\U_\jph^\pm$ by solving the nonlinear equations
\begin{equation*}
\E\big(\U_\jph^\pm\big)=\E_\jph^\pm,
\end{equation*}
which, in fact, require solving the cubic equations
\begin{equation}
\hf\left(\frac{q_\jph^\pm}{h_\jph^\pm}\right)^2+\theta_\jph^\pm\big(h_\jph^\pm+Z_\jph^\pm\big)+{\cal Q}_\jmh^\pm-\Ep_\jph^\pm=0
\label{3.7}
\end{equation}
for $h_\jph^\pm$. These equations can be solved exactly; see Appendix \ref{appA} for details.

The obtained point values $\U_\jph^\pm$ are to be used to evaluate the fluxes in \eref{3.4}. However, as it was pointed out in
\cite{Kurganov23}, modified point values $\widehat{\U}_\jph^\pm$ have to be used in the numerical diffusion term in \eref{3.3} to guarantee
that the difference $\widehat{\U}_\jph^+-\widehat{\U}_\jph^-$ vanishes at steady states. These values are obtained by solving the following
modified versions of the nonlinear equations \eref{3.7}:
\begin{equation}
\hf\left(\frac{q_\jph^\pm}{\widehat h_\jph^\pm}\right)^2+\theta_\jph^\pm\big(\widehat h_\jph^\pm+\widehat Z_\jph\big)-
\Ep_\jph^\pm=0,\quad\widehat Z_\jph:=\hf\big(Z_\jph^-+Z_\jph^+\big),
\label{3.8}
\end{equation}
which we solve for $\widehat h_\jph^\pm$. After that, we obtain $\widehat{(h\theta)}_\jph^\pm=\widehat h_\jph^\pm\theta_\jph^\pm$ and also
set $\widehat q_\jph^{\,\pm}=q_\jph^\pm$. Notice that at steady states, both $\theta$ and $q$ are constant and thus their reconstructed
values satisfy $\theta_\jph^-=\theta_\jph^+$ and $q_\jph^-=q_\jph^+$.

Next, $a_\jph^\pm$ in \eref{3.3} are the one-sided local speed of propagation, which can be estimated using the largest and smallest
eigenvalues of ${\cal A}(\U)$ as follows:
\begin{equation*}
a_\jph^-=\min\big\{u_\jph^--c_\jph^-,u_\jph^+-c_\jph^+,0\big\},\quad a_\jph^+=\max\big\{u_\jph^-+c_\jph^-,u_\jph^++c_\jph^+,0\big\}.
\end{equation*}

Finally, the high-order correction terms in \eref{3.2} are approximated as in \cite{CKX2025}:
\begin{equation*}
\begin{aligned}
&(\K_{xx})_\jph=\frac{1}{12(\dx)^{2}}\left[-\bm{{\cal K}}_{j-\frac{3}{2}}^{\rm FV}+16\bm{{\cal K}}_\jmh^{\rm FV}-
30\bm{{\cal K}}_\jph^{\rm FV}+16\bm{{\cal K}}_{j+\frac{3}{2}}^{\rm FV}-\bm{{\cal K}}_{j+\frac{5}{2}}^{\rm FV}\right],\\
&(\K_{xxxx})_\jph=\frac{1}{(\dx)^{4}}\left[\bm{{\cal K}}_{j-\frac{3}{2}}^{\rm FV}-4\bm{{\cal K}}_\jmh^{\rm FV}+6\bm{{\cal K}}_\jph^{\rm FV}-
4\bm{{\cal K}}_{j+\frac{3}{2}}^{\rm FV}+\bm{{\cal K}}_{j+\frac{5}{2}}^{\rm FV}\right].
\end{aligned}
\end{equation*}
\begin{remark}
\rm One can show that the resulting PCCU-5 scheme is WB in the sense that it is capable of exactly preserving moving-water steady states.
The proof is similar to that in \cite{Kurganov23}, and thus we omit it for the sake of brevity.
\end{remark}

Note that at steady states \eref{1.2} and \eref{1.5}, both $\theta$ and $q$ are constant, which implies
$\widehat{(h\theta)}_\jph^-=\widehat{(h\theta)}_\jph^+$. However, at steady states \eref{1.3} and \eref{1.4}, the values
$\widehat{(h\theta)}_\jph^-$ and $\widehat{(h\theta)}_\jph^+$ may be different. We therefore replace the third component of the numerical
flux \eref{3.3} with
$$
{\cal K}_{\jph}^{{\rm FV},(3)}=\frac{a_\jph^+q_\jph^-\theta_\jph^--a_\jph^-q_\jph^+\theta_\jph^+}{a_\jph^+-a_\jph^-}+
\frac{a_\jph^+a_\jph^-}{a_\jph^+-a_\jph^-}\,H_\jph\left(\widehat{(h\theta)}_\jph^+-\widehat{(h\theta)}_\jph^-\right),
$$
where $H_\jph=H\big(\varphi_\jph\big)$ is the same numerical diffusion switch used in \cite{CKL23} with
\begin{equation}\label{3.10}
H(\varphi):=\frac{400\varphi^8}{1+400\varphi^8},
\end{equation}
and
\begin{equation*}
\varphi_\jph:=\frac{\big|K_{j+1}^{(2)}-K_j^{(2)}\big|}{\dx}\cdot\frac{x_{N+\hf}-x_{\hf}}{\max\big\{K_{j+1}^{(2)},K_j^{(2)}\big\}},\quad
K_j^{(2)}:=\frac{\big(K^{(2)}\big)_\jph^-+\big(K^{(2)}\big)_\jmh^+}{2}.
\end{equation*}
We emphasize that the factor $H_\jph$ vanishes at steady states and is very small when the solution is (locally) close to a steady state.

\section{Two-Dimensional Ripa Model}\label{sec4}
The 2-D Ripa model reads as
\begin{equation}
\U_t+\F(\U)_x+\G(\U)_y= \B(\U)\U_x+\bm D(\U)\U_y ,
\label{4.1}
\end{equation}
where the conservative variables $\U$, flux functions $\F(\U)$ and $\G(\U)$, and matrices $B(\U)$ and $D(\U)$ are given by
\begin{equation}
\begin{aligned}
&\U=\begin{pmatrix}h\\q^x\\q^y\\h\theta\\Z\end{pmatrix},\quad\F(\U)=\begin{pmatrix}q^x\\q^xu+P\\[0.5ex]\dfrac{q^xq^y}{h}\\[1.0ex]
q^x\theta\\0\end{pmatrix},\quad\G(\U)=\begin{pmatrix}q^y\\[0.5ex]\dfrac{q^xq^y}{h}\\[0.5ex]
q^yv+P\\q^y\theta\\0\end{pmatrix},\\
&B(\U)=\begin{pmatrix}0&0&0&0&0\\0&0&0&0&-h\theta\\0&0&0&0&0\\0&0&0&0&0\\0&0&0&0&0\end{pmatrix},\quad
D(\U)=\begin{pmatrix}0&0&0&0&0\\0&0&0&0&0\\0&0&0&0&-h\theta\\0&0&0&0&0\\0&0&0&0&0\end{pmatrix}.
\label{4.2}
\end{aligned}
\end{equation}
Here, $y$ is the second spatial coordinate, $v(x,y,t)$ is the $y$-velocity, $q^x(x,t)=h(x,y,t)u(x,y,t)$ and $q^y(x,t)=h(x,y,t)v(x,y,t)$ are
the $x$- and $y$-discharges. The eigenvalues of ${\cal A}^x{(\U)}:=\frac{\partial\F}{\partial\U}-B(\U)$ and
${\cal A}^y{(\U)}:=\frac{\partial\G}{\partial\U}-D(\U)$ are
$$
\begin{aligned}
&\lambda_1({\cal A}^x)=u-c<\lambda_2({\cal A}^x)=\lambda_3({\cal A}^x)=u<\lambda_4({\cal A}^x)=u+c,&&\lambda_5({\cal A}^x)=0,\\
&\lambda_1({\cal A}^y)=v-c<\lambda_2({\cal A}^y)=\lambda_3({\cal A}^y)=v<\lambda_4({\cal A}^y)=v+c,&&\lambda_5({\cal A}^y)=0.
\end{aligned}
$$

As in the 1-D case, the 2-D Ripa system \eref{4.1}--\eref{4.2} can be rewritten in the quasi-conservative form
\begin{equation}
\U_t+\K(\U)_x+\bm L(\U)_y=\bm0,
\label{4.3}
\end{equation}
with the global fluxes
\begin{equation}
\begin{aligned}
&\K(\U)=\F(\U)-\R^x(\U),\quad\R^x(\U)=\int\limits_{\widehat x}^xB(\U(\xi,y,t))\U_\xi(\xi,y,t)\,{\rm d}\xi\\
&\hspace{5.4cm}=\Big(0,-\int\limits_{\widehat x}^xh(\xi,y,t)\theta(\xi,y,t)Z_\xi(\xi,y)\,{\rm d}\xi,0,0,0\Big)^\top,\\
&\bm L(\U)=\G(\U)-\R^y(\U)\quad\R^y(\U)=\int\limits_{\widehat y}^yD(\U(x,\eta,t))\U_\eta(x,\eta,t)\,{\rm d}\eta\\
&\hspace{5.15cm}=\Big(0,-\int\limits_{\widehat y}^yh(x,\eta,t)\theta(x,\eta,t)Z_\eta(x,\eta)\,{\rm d}\eta,0,0,0,0\Big)^\top,
\end{aligned}
\label{4.4a}
\end{equation}
where $\widehat x$ and $\widehat y$ are arbitrary numbers.

\paragraph{Steady States and Equilibrium Variables.} The 2-D Ripa system \eref{4.1}--\eref{4.2} admits several steady states including both
still- and moving-water ones. There are the same three types of still-water equilibria as in the 1-D case:

\noindent
$\bullet$ {\em ``Lake-at-rest''} steady states satisfying
\begin{equation}
u=v\equiv0,\quad h+Z\equiv{\rm Const},\quad\theta\equiv{\rm Const};
\label{4.4}
\end{equation}
$\bullet$ {\em Isobaric} steady states satisfying
\begin{equation}
u=v\equiv0,\quad P=\hf h^2\theta\equiv{\rm Const},\quad Z\equiv{\rm Const},
\label{4.5}
\end{equation}
where $P$ is the pressure;

\noindent
$\bullet$ {\em Constant water} depth steady states satisfying
\begin{equation}
u=v\equiv0,\quad h\equiv{\rm Const},\quad Z+\hf h\ln\theta\equiv{\rm Const}.
\label{4.6}
\end{equation}
On the other hand, 2-D moving-water steady states are substantially more complicated than their 1-D counterparts. We will restrict our
consideration to the quasi 1-D moving-water equilibria: the $x$-directional,
\begin{equation}
\begin{aligned}
&q^x_x=q^y=\theta_x=\Ep^x_x\equiv0,\quad\Ep^x:=\frac{u^2}{2}+\theta(h+Z)+{\cal Q}^x,\\
&{\cal Q}^x:=-\int\limits_{\widehat x}^x\left[\sqrt{2P(\xi,y,t)}\big(\sqrt{\theta(\xi,y,t)}\big)_\xi+Z(\xi,y)\theta_\xi(\xi,y,t)\right]
{\rm d}\xi,
\end{aligned}
\label{4.7}
\end{equation}
and the $y$-directional,
\begin{equation}
\begin{aligned}
&q^x=q^y_y=\theta_y=\Ep^y_y\equiv0,\quad\Ep^y:=\frac{v^2}{2}+\theta(h+Z)+{\cal Q}^y,\\
&{\cal Q}^y:=-\int\limits_{\widehat y}^y\left[\sqrt{2P(x,\eta,t)}\big(\sqrt{\theta(x,\eta,t)}\big)_\eta+Z(x,\eta)\theta_\eta(x,\eta,t)
\right]{\rm d}\eta,
\end{aligned}
\label{4.8}
\end{equation}
ones.

It is easy to verify that:

\noindent
(i) at all of the still-water steady states \eref{4.4}--\eref{4.6}, $\Ep^x_x=\Ep^y_y\equiv0$;

\noindent
(ii) at the steady states \eref{4.4}--\eref{4.7},
\begin{equation}
\K(\U)_x=M^x(\U)\E^x(\U)_x\equiv\bm0,
\label{4.9}
\end{equation}
where
\begin{equation}
M^x(\U)=\begin{pmatrix}1&0&0&0&0\\u&h&0&0&0\\v&0&q^x&0&0\\\theta&0&0&q^x&0\\0&0&0&0&0\end{pmatrix},\quad
\E^x(\U)=\begin{pmatrix}q^x\\\Ep^x\\v\\\theta\\Z\end{pmatrix};
\label{4.10}
\end{equation}

\noindent
(iii) at the steady states \eref{4.4}--\eref{4.6} and \eref{4.8},
\begin{equation}
\bm L(\U)_y=M^y(\U)\E^y(\U)_y\equiv\bm0,
\label{4.11}
\end{equation}
where
\begin{equation}
M^y(\U)=\begin{pmatrix}1&0&0&0&0\\u&0&q^y&0&0\\v&h&0&0&0\\\theta&0&0&q^y&0\\0&0&0&0&0\end{pmatrix},\quad
\E^y(\U)=\begin{pmatrix}q^y\\\Ep^y\\u\\\theta\\Z\end{pmatrix}.
\label{4.12}
\end{equation}
According to \eref{4.9}--\eref{4.10} and \eref{4.11}--\eref{4.12}, we identify $\E^x$ and $\E^y$ as the equilibrium variables since at
steady states the first four components of $\E^x$ are constant along the $x$-direction and the first four components of $\E^y$ are constant
along the $y$-direction.

We also remark that the system \eref{4.3}--\eref{4.4a} can be rewritten in the equivalent (for smooth solutions) form using the equilibrium
variables $\E^x$
$$
\E^x(\U)_t+C^x(\U)\E^x(\U)_x+D^x(\U)\U_y=\bm I^x(\U),
$$
where
\begin{equation}
C^x(\U)=\begin{pmatrix}u&h&0&0&0\\[0.5ex]\theta&u&0&\dfrac{q^x}{2}&0\\[1.5ex]0&0&u&0&0\\0&0&0&u&0\\0&0&0&0&0\end{pmatrix},\quad
D^x(\U)=\frac{1}{h}\begin{pmatrix}-huv&v&u&0&0\\[0.5ex]-\dfrac{\theta q^y}{2}&0&h\theta&\dfrac{q^y}{2}&0\\[2.0ex]
\dfrac{h\theta}{2}-v^2&0&v&\dfrac{h}{2}&\theta\\[1.0ex]-v\theta&0&0&v&0\\0&0&0&0&0
\end{pmatrix},
\label{4.14}
\end{equation}
and $\bm I^x(\bm U)=\big(0,I^x(\bm U),0,0,0\big)^\top$ with
$$
\begin{aligned}
&I^x(\bm U)=\int\limits_{\widehat x}^x\bigg[\Big(\hf h(\xi,y,t)u_\xi(\xi,y,t)-u(\xi,y,t)Z_\xi(\xi,y)+\hf q^y_y(\xi,y,t)\Big)
\theta_\xi(\xi,y,t)\\
&\hspace*{1.6cm}-\Big(\hf h_\xi(\xi,y,t)+Z_\xi(\xi,y)\Big)v(\xi,y,t)\theta_y(\xi,y,t)\bigg]{\rm d}\xi.
\end{aligned}
$$
One can also rewrite the system \eref{4.3}--\eref{4.4a} in another equivalent form using the equilibrium variables $\E^y$:
\begin{align*}
&\E^y(\U)_t+D^y(\U)\U_x+C^y(\U)\E^y(\U)_y=\bm I^y(\U),
\end{align*}
where
\begin{equation}
C^y(\U)=\begin{pmatrix}v&h&0&0&0\\\theta&v&0&\dfrac{q^y}{2}&0\\0&0&v&0&0\\0&0&0&v&0\\0&0&0&0&0\end{pmatrix},\quad
D^y(\U)=\frac{1}{h}\begin{pmatrix}-huv&v&u&0&0\\[0.5ex]-\dfrac{\theta q^x}{2}&h\theta&0&\dfrac{q^x}{2}&0\\[2.0ex]
\dfrac{h\theta}{2}-u^2&u&0&\dfrac{h}{2}&\theta\\[1.0ex]-u\theta&0&0&u&0\\0&0&0&0&0
\end{pmatrix},
\label{4.15a}
\end{equation}
and $\bm I^y(\bm U)=\big(0,I^y(\bm U),0,0,0\big)^\top$ with
$$
\begin{aligned}
&I^y(\bm U)=\int\limits_{\widehat y}^y\bigg[\Big(\hf h(x,\eta,t)v_\eta(x,\eta,t)-v(x,\eta,t)Z_\eta(x,\eta)+\hf q^x_x(x,\eta,t)\Big)
\theta_\eta(x,\eta,t)\\
&\hspace*{1.6cm}-\Big(\hf h_\eta(x,\eta,t)+Z_\eta(x,\eta)\Big)u(x,\eta,t)\theta_x(x,\eta,t)\bigg]{\rm d}\eta.
\end{aligned}
$$

We stress that the matrices $C^x$ and $C^y$ have complete eigensystems and they can be diagonalized using
$$
Q^x(\U)=\begin{pmatrix}0&0&-\dfrac{q^x}{2\theta}&\dfrac{c}{\theta}&-\dfrac{c}{\theta}\\[1.5ex]0&0&0&1&1\\0&1&0&0&0\\0&0&1&0&0\\1&0&0&0&0
\end{pmatrix},\quad[Q^x(\U)]^{-1}=\begin{pmatrix}0&0&0&0&1\\0&0&1&0&0\\0&0&0&1&0\\[0.5ex]
\dfrac{\theta}{2c}&\dfrac{1}{2}&0&-\dfrac{q^x}{4c}&0\\[2.0ex]-\dfrac{\theta}{2c}&\dfrac{1}{2}&0&-\dfrac{q^x}{4c}&0\end{pmatrix},
$$
and
$$
Q^y(\U)=\begin{pmatrix}0&0&-\dfrac{q^y}{2\theta}&\dfrac{c}{\theta}&-\dfrac{c}{\theta}\\[1.5ex]0&0&0&1&1\\0&1&0&0&0\\0&0&1&0&0\\1&0&0&0&0
\end{pmatrix},\quad[Q^y(\U)]^{-1}=\begin{pmatrix}0&0&0&0&1\\0&0&1&0&0\\0&0&0&1&0\\[0.5ex]
\dfrac{\theta}{2c}&\dfrac{1}{2}&0&-\dfrac{q^y}{4c}&0\\[2.0ex]-\dfrac{\theta}{2c}&\dfrac{1}{2}&0&-\dfrac{q^y}{4c}&0\end{pmatrix},
$$
respectively, namely, $[Q^x(\U)]^{-1}C^x(\U)Q^x(\U)$ and $[Q^y(\U)]^{-1}C^y(\U)Q^y(\U)$ are diagonal matrices.

\section{Two-Dimensional Scheme}\label{sec5}
In this section, we present the 2-D scheme. We first introduce a uniform mesh consisting of the cells
$I_{j,k}=[x_\jmh,x_\jph]\times[y_\kmh,y_\kph]$ centered at $(x_j,y_k)$ with $x_j=\hf(x_\jmh+x_\jph)$, $y_k=\hf(y_\kmh+y_\kph)$,
$x_\jph-x_\jmh\equiv\dx$, $y_\kph-y_\kmh\equiv\dy$, $j=1,\dots,N_x$, $k=1,\dots,N_y$. We then introduce the point values
${\U_{j,k}}:\approx\U(x_j,y_k,t)$ and evolve them in time using the following semi-discretization:
\begin{equation}
\frac{{\rm d}\U_{j,k}}{{\rm d}t}=-\frac{\bm{{\cal K}}_{\jph,k}-\bm{{\cal K}}_{\jmh,k}}{\dx}-
\frac{\bm{{\cal L}}_{j,\kph}-\bm{{\cal L}}_{j,\kmh}}{\dy},
\label{5.1a}
\end{equation}
where the global numerical fluxes $\bm{{\cal K}}_{\jph,k}$ and $\bm{{\cal L}}_{j,\kph}$ are defined within the A-WENO framework as
$$
\begin{aligned}
&\bm{{\cal K}}_{\jph,k}=\bm{{\cal K}}_{\jph,k}^{\rm FV}-\frac{(\dx)^2}{24}(\K_{xx})_{\jph,k}+\frac{7(\dx)^4}{5760}(\K_{xxxx})_{\jph,k},\\
&\bm{{\cal L}}_{j,\kph}=\bm{{\cal L}}_{j,\kph}^{\rm FV}-\frac{(\dy)^2}{24}(\bm L_{yy})_{j,\kph}+
\frac{7(\dy)^4}{5760}(\bm L_{yyyy})_{j,\kph}.
\end{aligned}
$$
Here, $\bm{{\cal K}}_{\jph,k}^{\rm FV}$ and $\bm{{\cal L}}_{j,\kph}^{\rm FV}$ are the CU numerical fluxes
\begin{equation}
\begin{aligned}
&\bm{{\cal K}}_{\jph,k}^{\rm FV}=\frac{a_{\jph,k}^+\K_{\jph,k}^--a_{\jph,k}^-\K_{\jph,k}^+}{a_{\jph,k}^+-a_{\jph,k}^-}+
\frac{a_{\jph,k}^+a_{\jph,k}^-}{a_{\jph,k}^+-a_{\jph,k}^-}\left(\widehat\U_{\jph,k}^+-\widehat\U_{\jph,k}^-\right),\\
&\bm{{\cal L}}_{j,\kph}^{\rm FV}=\frac{a_{j,\kph}^+\bm L_{j,\kph}^--a_{j,\kph}^-\bm L_{j,\kph}^+}{a_{j,\kph}^+-a_{j,\kph}^-}+
\frac{a_{j,\kph}^+a_{j,\kph}^-}{a_{j,\kph}^+-a_{j,\kph}^-}\left(\widehat\U_{j,\kph}^+-\widehat\U_{j,\kph}^-\right),
\label{5.1}
\end{aligned}
\end{equation}
where the one-sided point values of the global flux function are computed as
\begin{equation*}
\begin{aligned}
&\K_{\jph,k}^\pm=\F_{\jph,k}^\pm-(\R^x)_{\jph,k}^\pm,&&\F_{\jph,k}^\pm:=\F\big(\U_{\jph,k}^\pm\big),\\
&\bm L_{\jph,k}^\pm=\G_{j,\kph}^\pm-(\R^y)_{j,\kph}^\pm,&&\G_{j,\kph}^\pm:=\G\big(\U_{j,\kph}^\pm\big),
\end{aligned}
\end{equation*}
where the one-sided values of the global parts of the fluxes, $(\R^x)_{\jph,k}^\pm$ and $(\R^y)_{j,\kph}^\pm$, are determined as in the 1-D
case (details are omitted for the sake of brevity) and the one-sided point values $\U_{\jph,k}^\pm$ and $\U_{j,\kph}^\pm$ are obtained with
the help of the fifth-order Ai-WENO-Z interpolation applied using the LCD of the equilibrium variables in a dimension-by-dimension manner.

In order to apply the fifth-order Ai-WENO-Z interpolation to the equilibrium variables, we first compute
\begin{equation*}
\Ep^x_{j,k}=\frac{u_{j,k}^2}{2}+\theta_{j,k}(h_{j,k}+Z_{j,k})+{\cal Q}^x_{j,k},\quad
\Ep^y_{j,k}=\frac{v_{j,k}^2}{2}+\theta_{j,k}(h_{j,k}+Z_{j,k})+{\cal Q}^y_{j,k},
\end{equation*}
where the values ${\cal Q}^x_{j,k}$ and ${\cal Q}^y_{j,k}$ are computed using the path-conservative integration as in \cite{CKL23}. We
first set $\widehat x= x_\hf$ and $\widehat y=y_\hf$, hence $({\cal Q}^x)_{\hf,k}^-=0$ $\forall k$, $({\cal Q}^y)_{j,\hf}^-=0$ $\forall j$,
and then $({\cal Q}^x)_{\hf,k}^+={\cal B}^x_{\bm\Psi,\hf,k}$ $\forall k$, $({\cal Q}^y)_{j,\hf}^+={\cal B}^y_{\bm\Psi,j,\hf}$ $\forall j$,
and the rest of the terms $({\cal Q}^x)_{\jph,k}^\pm$ and $({\cal Q}^y)_{j,\kph}^\pm$ are evaluated recursively by
\begin{equation*}
\begin{aligned}
&({\cal Q}^x)_{\jph,k}^-=({\cal Q}^x)_{\jmh,k}^++{\cal B}^x_{j,k},&&({\cal Q}^x)_{\jph,k}^+=({\cal Q}^x)_{\jph,k}^-+
{\cal B}^x_{\bm\Psi,\jph,k},\\
&({\cal Q}^y)_{j,\kph}^-=({\cal Q}^y)_{j,\kmh}^++{\cal B}^y_{j,k},&&({\cal Q}^y)_{j,\kph}^+=({\cal Q}^y)_{j,\kph}^-+
{\cal B}^y_{\bm\Psi,j,\kph},
\end{aligned}
\end{equation*}
for $j=1,\dots,N_x$ and $k=1,\dots,N_y$. Here,
\begin{equation*}
{\cal B}^x_{j,k}\approx-\int\limits_{I_{j,k}}\left[\sqrt{2P}(\sqrt{\theta})_x+Z\theta_x\right]{\rm d}x,\quad
{\cal B}^y_{j,k}\approx-\int\limits_{I_{j,k}}\left[\sqrt{2P}(\sqrt{\theta})_y+Z\theta_y\right]{\rm d}y
\end{equation*}
are computed using the Newton-Cotes quadrature described in \cite{FGKW}, and the terms
$$
\begin{aligned}
&{\cal B}^x_{\bm\Psi,\jph,k}\approx-\int\limits_0^1\left[\sqrt{2P\big(\bm\Psi_{\jph,k}(s)\big)}
\left(\sqrt{\theta\big(\bm\Psi_{\jph,k}(s)\big)}\right)'+Z\big(\bm\Psi_{\jph,k}(s)\big)\left(\theta\big(\bm\Psi_{\jph,k}(s)\big)\right)'
\right]{\rm d}s,\\
&{\cal B}^y_{\bm\Psi,j,\kph}\approx-\int\limits_0^1\left[\sqrt{2P\big(\bm\Psi_{j,\kph}(s)\big)}
\left(\sqrt{\theta\big(\bm\Psi_{j,\kph}(s)\big)}\right)'+Z\big(\bm\Psi_{j,\kph}(s)\big)\left(\theta\big(\bm\Psi_{j,\kph}(s)\big)\right)'
\right]{\rm d}s
\end{aligned}
$$
are obtained by approximating the path integrals using the linear segment paths for $\theta$, $\sqrt{2P}$, and $Z$, which results in
\begin{equation}
\begin{aligned}
{\cal B}^x_{\bm\Psi,\jph,k}&=-\hf\left(\sqrt{2P_{\jph,k}^-}+\sqrt{2P_{\jph,k}^+}\right)
\left(\sqrt{\theta_{\jph,k}^+}-\sqrt{\theta_{\jph,k}^-}\right)\\
&-\hf\left(Z_{\jph,k}^-+Z_{\jph,k}^+\right)\left(\theta_{\jph,k}^+-\theta_{\jph,k}^-\right),\\
{\cal B}^y_{\bm\Psi,j,\kph}&=-\hf\left(\sqrt{2P_{j,\kph}^-}\sqrt{2P_{j,\kph}^+}\right)
\left(\sqrt{\theta_{j,\kph}^+}-\sqrt{\theta_{j,\kph}^-}\right)\\
&-\hf\left(Z_{j,\kph}^-+Z_{j,\kph}^+\right)\left(\theta_{j,\kph}^+-\theta_{j,\kph}^-\right).
\end{aligned}
\label{5.2}
\end{equation}
Next, the terms ${\cal Q}^x_{j,k}$ and ${\cal Q}^y_{j,k}$ are computed by
$$
\begin{aligned}
&({\cal Q}^x)_{j,k}\approx({\cal Q}^x)_{\jmh,k}^+-\int\limits_{x_\jmh}^{x_j}\left[\sqrt{2P}(\sqrt{\theta})_x+Z\theta_x\right]{\rm d}x,\\
&({\cal Q}^y)_{j,k}\approx({\cal Q}^x)_{j,\kmh}^+-\int\limits_{y_\kmh}^{y_k}\left[\sqrt{2P}(\sqrt{\theta})_y+Z\theta_y\right]{\rm d}y,
\end{aligned}
$$
where the integrals are also approximated by the Newton-Cotes quadrature described in \cite{FGKW}.

It should be pointed out that the point values of $\theta$, $P$, and $Z$ used in \eref{5.2} as well as in the Newton-Cotes quadratures
applied to evaluate the integrals, are obtained by applying the Ai-WENO-Z interpolation directly to the fields $\theta$, $P$, and $Z$. Since
these values are not WB, we now proceed with the non-oscillatory WB interpolation performed in the local characteristic equilibrium
variables.

To this end, we evaluate the matrices $C_{\jph,k}:=C^x(\overline{\U}_{\jph,k})$ and $C_{j,\kph}:=C^y(\overline{\U}_{j,\kph})$, where $C^x$
and $C^y$ are given in \eref{4.14} and \eref{4.15a} and
$$
\overline{\U}_{\jph,k}:=\hf(\U_{j,k}+\U_{j+1,k}),\quad\overline{\U}_{j,\kph}:=\hf(\U_{j,k}+\U_{j,k+1}).
$$
We then introduce the matrices
$$
\begin{aligned}
&Q^x_{\jph,k}:=Q^x\big(\overline{\U}_{\jph,k}\big),&&\big(Q^x_{\jph,k}\big)^{-1}:=\big(Q^x\big(\overline{\U}_{\jph,k}\big)\big)^{-1},\\
&Q^y_{j,\kph}:=Q^y\big(\overline{\U}_{j,\kph}\big),&&\big(Q^y_{j,\kph}\big)^{-1}:=\big(Q^y\big(\overline{\U}_{j,\kph}\big)\big)^{-1},
\end{aligned}
$$
switch to the local characteristic equilibrium variables in the neighborhoods of $(x_\jph,y_k)$ and $(x_j,y_\kph)$:
\begin{equation*}
\bm\Gamma^x_\ell=\big(Q^x_{\jph,k}\big)^{-1}\E^x_{j+\ell,k},\quad\bm\Gamma^y_\ell=\big(Q^y_{j,\kph}\big)^{-1}\E^y_{j,k+\ell},\quad
\ell=-2,\ldots,3,
\end{equation*}
apply the Ai-WENO-Z interpolation from \cite{FGKW,ZGW} to each component of $\bm\Gamma^x_\ell$ and $\bm\Gamma^y_\ell$ to reconstruct the
one-sided point values $(\bm\Gamma^x)_\hf^\pm$ and $(\bm\Gamma^y)_\hf^\pm$, and then transform them back to equilibrium variables to obtain
\begin{equation*}
(\E^x)_{\jph,k}^\pm=Q^x_{\jph,k}(\bm\Gamma^x)_\hf^\pm,\quad(\E^y)_{j,\kph}^\pm=Q^y_{j,\kph}(\bm\Gamma^y)_\hf^\pm.
\end{equation*}

Equipped with $(\E^x)_{\jph,k}^\pm$ and $(\E^y)_{j,\kph}^\pm$, we now evaluate $\U_{\jph,k}^\pm$ and $\U_{j,\kph}^\pm$ by solving the
nonlinear equations
\begin{equation*}
\E^x\big(\U_{\jph,k}^\pm\big)=(\E^x)_{\jph,k}^\pm,\quad\E^y\big(\U_{j,\kph}^\pm\big)=(\E^y)_{j,\kph}^\pm,
\end{equation*}
which, in fact, require solving the cubic equations
\begin{equation*}
\hf\left(\frac{(q^x)_{\jph,k}^\pm}{h_{\jph,k}^\pm}\right)^2+\theta_{\jph,k}^\pm\big(h_{\jph,k}^\pm+Z_{\jph,k}^\pm\big)+
({\cal Q}^x)_{\jph,k}^\pm-(\Ep^x)_{\jph,k}^\pm=0
\end{equation*}
and
\begin{equation*}
\hf\left(\frac{(q^y)_{j,\kph}^\pm}{h_{j,\kph}^\pm}\right)^2+\theta_{j,\kph}^\pm\big(h_{j,\kph}^\pm+Z_{j,\kph}^\pm\big)+
({\cal Q}^y)_{j,\kph}^\pm-(\Ep^y)_{j,\kph}^\pm=0
\end{equation*}
for $h_{\jph,k}^\pm$ and $h_{j,\kph}^\pm$, respectively. These equations can be solved exactly; see Appendix \ref{appA} for details.

Next, to guarantee the WB property of the resulting scheme, we need to ensure that the differences
$\widehat\U_{\jph,k}^+-\widehat\U_{\jph,k}^-$ and $\widehat\U_{j,\kph}^+-\widehat\U_{j,\kph}^-$ vanish at steady states. To this end, we
first solve
\begin{equation}
\hf\left(\frac{(q^x)_{\jph,k}^\pm}{\widehat h_{\jph,k}^\pm}\right)^2+\theta_{\jph,k}^\pm\big(\widehat h_{\jph,k}^\pm+\widehat Z_{\jph,k}^\pm
\big)+({\cal Q}^x)_{\jph,k}^\pm-(\Ep^x)_{\jph,k}^\pm=0,
\label{5.3}
\end{equation}
and
\begin{equation}
\hf\left(\frac{(q^y)_{j,\kph}^\pm}{\widehat h_{j,\kph}^\pm}\right)^2+\theta_{j,\kph}^\pm\big(\widehat h_{j,\kph}^\pm+\widehat Z_{j,\kph}^\pm
\big)+({\cal Q}^y)_{j,\kph}^\pm-(\Ep^y)_{j,\kph}^\pm=0,
\label{5.4}
\end{equation}
for $\widehat h_{\jph,k}^\pm$ and $\widehat h_{j,\kph}^\pm$, respectively. In \eref{5.3} and \eref{5.4},
\begin{equation*}
\widehat Z_{\jph,k}^\pm:=\hf\big(Z_{\jph,k}^-+Z_{\jph,k}^+\big),\quad\widehat Z_{j,\kph}^\pm:=\hf\big(Z_{j,\kph}^-+Z_{j,\kph}^+\big).
\end{equation*}
After that, we obtain $(\widehat{h\theta})_{\jph,k}^\pm=\widehat h_{\jph,k}^\pm\theta_{\jph,k}^\pm$ and $(\widehat{h\theta})_{j,\kph}^\pm=
\widehat h_{j,\kph}^\pm\theta_{j,\kph}^\pm$, and also set $(\widehat{q^x})_{\jph,k}^\pm=(q^x)_{\jph,k}^\pm$ and
$(\widehat{q^y})_{j,\kph}^\pm=(q^y)_{j,\kph}^\pm$.

Next, $a_{\jph,k}^\pm$ and $a_{j,\kph}^\pm$ in \eref{5.1} are one-sided local speed of propagation, which can be estimated using the largest
and smallest eigenvalues as follows:
\begin{equation*}
\begin{aligned}
&a_{\jph,k}^-=\min\big\{u_{\jph,k}^--c_{\jph,k}^-,u_{\jph,k}^+-c_{jph,k}^+,0\big\},\\
&a_{\jph,k}^+=\max\big\{u_{\jph,k}^-+c_{\jph,k}^-,u_{\jph,k}^++c_{\jph,k}^+,0\big\},\\
&a_{j,\kph}^-=\min\big\{v_{j,\kph}^--c_{j,\kph}^-,v_{j,\kph}^+-c_{j,\kph}^+, 0\big\},\\
&a_{j,\kph}^+=\max\big\{v_{j,\kph}^-+c_{j,\kph}^-,v_{j,\kph}^++c_{j,\kph}^+,0\big\}.
\end{aligned}
\end{equation*}

Finally, $(\K_{xx})_{\jph,k}$, $(\K_{xxxx})_{\jph,k}$, $(\bm L_{yy})_{j,\kph}$, $(\bm L_{yyyy})_{j,\kph}$ are the high-order correction
terms, which we compute as in \cite{CKX2025}:
$$
\begin{aligned}
&(\K_{xx})_{\jph,k}=\frac{1}{12(\dx)^{2}}\left[-\bm{{\cal K}}_{j-\frac{3}{2},k}^{\rm FV}+16\bm{{\cal K}}_{\jmh,k}^{\rm FV}-
30\bm{{\cal K}}_{\jph,k}^{\rm FV}+16\bm{{\cal K}}_{j+\frac{3}{2},k}^{\rm FV}-\bm{{\cal K}}_{j+\frac{5}{2},k}^{\rm FV}\right],\\
&(\K_{xxxx})_{\jph,k}=\frac{1}{(\dx)^{4}}\left[\bm{{\cal K}}_{j-\frac{3}{2},k}^{\rm FV}-4\bm{{\cal K}}_{\jmh,k}^{\rm FV}+
6\bm{{\cal K}}_{\jph,k}^{\rm FV}-4\bm{{\cal K}}_{j+\frac{3}{2},k}^{\rm FV}+\bm{{\cal K}}_{j+\frac{5}{2},k}^{\rm FV}\right],\\
&(\bm L_{yy})_{j,\kph}=\frac{1}{12(\dy)^{2}}\left[-\bm{{\cal L}}_{j,k-\frac{3}{2}}^{{FV}}+16\bm{{\cal L}}_{j,\kmh}^{\rm FV}-
30\bm{{\cal L}}_{j,\kph}^{\rm FV}+16\bm{{\cal L}}_{j,k+\frac{3}{2}}^{\rm FV}-\bm{{\cal L}}_{j+\frac{5}{2}}^{\rm FV}\right],\\
&(\bm L_{yyyy})_{j,\kph}=\frac{1}{(\dy)^{4}}\left[\bm{{\cal L}}_{j,k-\frac{3}{2}}^{\rm FV}-4\bm{{\cal L}}_{j,\kmh}^{\rm FV}+
6\bm{{\cal L}}_{j,\kph}^{\rm FV}-4\bm{{\cal L}}_{j,k+\frac{3}{2}}^{\rm FV}+\bm{{\cal L}}_{j,k+\frac{5}{2}}^{\rm FV}\right].
\end{aligned}
$$

Similarly to the 1-D scheme, we switch off the numerical diffusion when the solution is at or near steady states. To this end, we replace
the fourth component of the numerical fluxes in \eref{5.1} with
$$
\begin{aligned}
{\cal K}_{\jph,k}^{{\rm FV},(4)}&=\frac{a_{\jph,k}^+(q^x)_{\jph,k}^-\theta_{\jph,k}^--a_{\jph,k}^-(q^x)_{\jph,k}^+\theta_{\jph,k}^+}
{a_{\jph,k}^+-a_{\jph,k}^-}\\
&+\frac{a_{\jph,k}^+a_{\jph,k}^-}{a_{\jph,k}^+-a_{\jph,k}^-}\,H_{\jph,k}
\left(\widehat{(h\theta)}_{\jph,k}^+-\widehat{(h\theta)}_{\jph,k}^-\right),\\[1.0ex]
{\cal L}_{j,\kph}^{{\rm FV},(4)}&=\frac{a_{j,\kph}^+(q^y)_{j,\kph}^-\theta_{j,\kph}^--a_{j,\kph}^-(q^y)_{j,\kph}^+\theta_{j,\kph}^+}
{a_{j,\kph}^+-a_{j,\kph}^-}\\
&+\frac{a_{j,\kph}^+a_{j,\kph}^-}{a_{j,\kph}^+-a_{j,\kph}^-}\,H_{j,\kph}
\left(\widehat{(h\theta)}_{j,\kph}^+-\widehat{(h\theta)}_{j,\kph}^-\right),
\end{aligned}
$$
where $H_{\jph,k}=H\big(\varphi_{\jph,k}\big)$, $H_{j,\kph}=H\big(\varphi_{j,\kph}\big)$, the function $H$ is given by \eref{3.10}, and
\begin{equation*}
\begin{aligned}
&\varphi_{\jph,k}:=\frac{\big|K_{j+1,k}^{(2)}-K_{j,k}^{(2)}\big|}{\dx}\cdot
\frac{x_{N_x+\hf}-x_{\hf}}{\max\big\{K_{j+1,k}^{(2)},K_{j,k}^{(2)}\big\}},&&
K_{j,k}^{(2)}:=\frac{\big(K^{(2)}\big)_{\jph,k}^-+\big(K^{(2)}\big)_{\jmh.k}^+}{2},\\
&\varphi_{j,\kph}:=\frac{\big|L_{j,k+1}^{(3)}-L_{j,k}^{(3)}\big|}{\dy}\cdot
\frac{y_{N_y+\hf}-y_{\hf}}{\max\big\{L_{j,k+1}^{(3)},L_{j,k}^{(3)}\big\}},&&
L_{j,k}^{(3)}:=\frac{\big(L^{(3)}\big)_{j,\kph}^-+\big(L^{(3)}\big)_{j,\kmh}^+}{2}.
\end{aligned}
\end{equation*}

\section{Numerical Examples}\label{sec6}
In this section, we test the proposed PCCU-5 scheme on a number of 1-D and 2-D numerical examples. We also compare the PCCU-5 scheme with
its second-order counterpart (PCCU-2 scheme developed in \cite{CKL23}) and the PCCU-5 scheme, in which we reconstruct the equilibrium
variables without the characteristic decomposition (this scheme will be referred to as the PCCU-5-NCD scheme).

In all of the numerical examples, we have numerically solved the ODE systems \eref{3.1} and \eref{5.1a} using the three-stages third-order
strong stability preserving Runge-Kutta method; see, e.g., \cite{GKS,GST}. In all of the examples except for the accuracy tests (Examples
\ref{ex61} and \ref{ex66}), the time steps are selected adaptively using
$$
\dt=0.45\cdot\frac{\dx}{\max\limits_j\left\{|u_j|+\sqrt{h_j\theta_j}\right\}}
$$
and
$$
\dt=0.45\cdot\min\left[\frac{\dx}{\max\limits_{j,k}\left\{|u_{j,k}|+\sqrt{h_{j,k}\theta_{j,k}}\right\}},
\frac{\dy}{\max\limits_{j,k}\left\{|v_{j,k}|+\sqrt{h_{j,k}\theta_{j,k}}\right\}}\right]
$$
in the 1-D and 2-D examples, respectively.

\subsection{One-Dimensional Examples}
\begin{example}[1-D numerical accuracy test]\label{ex61}
\rm In the first example, we test the experimental order of accuracy of the proposed PCCU-5 scheme. The following smooth initial conditions,
\begin{equation*}
(h,q,\theta)\Big|_{(x,0)}=\big(1-Z(x),0.1,9.812[1-0.01\cos(2\pi x)]\big),
\end{equation*}
are prescribed in the computational domain $[0,1]$ subject to the periodic boundary conditions. The bottom topography is also smooth:
$Z(x)=0.1\sin(4\pi x)-1$, and the time step is selected as $\dt=0.45(\dx)^{5/3}$ to ensure that the spatial and temporal contributions of
the error are of the same order.

We compute the solution by the PCCU-5 scheme until the final time $t=1$ on a sequence of uniform meshes with $N=25$, $50$, $100$, $200$, and
$400$ cells. We also compute the reference solution on a much finer mesh with $N=6400$ cells. The $L^\infty$-errors and experimental rates
of convergence for $h$, $q$, and $h\theta$ are presented in Table \ref{tab61}. As one can see, the proposed PCCU-5 scheme achieves the
expected fifth order of accuracy.
\begin{table}[!htb]
\centering\small
\caption{Example \ref{ex61}: $L^\infty$-errors and experimental rates of convergence for $h$, $q$, and $h\theta$.}\label{tab61}
\begin{tabular}{cccccccccc}\hline
&&\multicolumn{2}{c}{$h$}&&\multicolumn{2}{c}{$q$}&&\multicolumn{2}{c}{$h\theta$}\\\cline{3-4}\cline{6-7}\cline{9-10}
$N$  &&Error   &Rate&&Error   &Rate&&Error   &Rate\\\hline
$25$ &&4.39E-05& ---&&1.55E-04& ---&&4.41E-04& ---\\
$50$ &&1.99E-06&4.46&&8.12E-06&4.26&&1.80E-05&4.53\\
$100$&&6.41E-08&4.96&&2.84E-07&4.84&&5.86E-07&4.94\\
$200$&&1.96E-09&5.03&&9.14E-09&4.96&&1.80E-08&5.02\\
$400$&&6.31E-11&4.96&&2.61E-10&5.13&&5.77E-10&4.96\\\hline
\end{tabular}
\end{table}
\end{example}

\begin{example}[1-D moving-water equilibria]\label{ex62}
\rm In this example taken from \cite{ZXX}, we demonstrate the ability of the proposed PCCU-5 scheme to exactly preserve moving-water
equilibria.

The computational domain in this example is $[0,25]$ and we consider the following continuous,
\begin{equation}
Z(x)=\left\{\begin{aligned}
&0.2-0.05(x-10)^2,&&8\le x\le12,\\
&0,&&\mbox{otherwise},
\end{aligned}\right.
\label{6.1}
\end{equation}
and discontinuous,
\begin{equation}
Z(x)=\left\{\begin{aligned}
&0.2,&& 8\le x\le12,\\
&0,&&\mbox{otherwise}.
\end{aligned}\right.
\label{5.2}
\end{equation}
bottom topographies. We consider the following three sets of initial and boundary conditions, which correspond to the subcritical,
supercritical, and transcritical equilibria:

\medskip
\noindent
$\bullet$ {\bf Subcritical equilibrium}:
\begin{equation}
\begin{aligned}
&(\Ep,q,\theta)\Big|_{(x,0)}=(\Ep_{\rm eq}(x),q_{\rm eq}(x),\theta_{\rm eq}(x))\equiv(110.33025,4.42\sqrt{5},49.06),\\
&h(25,t)=2,~~q(0,t)=4.42\sqrt{5};
\end{aligned}
\label{6.3}
\end{equation}

\medskip
\noindent
$\bullet$ {\bf Supercritical equilibrium}:
\begin{equation}
\begin{aligned}
&(\Ep,q,\theta)\Big|_{(x,0)}=(\Ep_{\rm eq}(x),q_{\rm eq}(x),\theta_{\rm eq}(x))\equiv(458.12,24\sqrt{5},49.06),\\
&h(0,t)=2,~~q(0,t)=24\sqrt{5};
\end{aligned}
\label{6.4}
\end{equation}

\medskip
\noindent
$\bullet$ {\bf Transcritical equilibrium}:
\begin{equation}
\begin{aligned}
&(\Ep,q,\theta)\Big|_{(x,0)}=(\Ep_{\rm eq}(x),q_{\rm eq}(x),\theta_{\rm eq}(x))\equiv(55.453570198891,1.53\sqrt{5},49.06),\\
&h(25,t)=0.405748088283403,~~q(0,t)=1.53\sqrt{5},
\end{aligned}
\label{6.5}
\end{equation}
where the downstream boundary condition ($h(25,t)=0.405748088283403$) is imposed only when the flow is subcritical there.

In all of the three cases, the rest of the boundary conditions are free (zero-order extrapolation).

Note that the initial data in \eref{6.3}--\eref{6.5} are prescribed in terms of the equilibrium variables. We therefore begin the
simulations by solving the cubic equations
\begin{equation}
\hf\left(\frac{(q_{\rm eq})_j}{(h_{\rm eq})_j}\right)^2+{(\theta_{\rm eq})_j}\big[(h_{\rm eq})_j+Z_j\big]-(\Ep_{\rm eq})_j=0
\label{5.6}
\end{equation}
for $(h_{\rm eq})_j$, and then set the initial data $h_j(0)=(h_{\rm eq})_j$. We solve \eref{5.6} exactly and we single out physically
relevant roots based on the initial Froude number, which is smaller than $1$ in the subcritical case, larger than $1$ in the supercritical
case, and in the transcritical case, the Froude number is smaller than $1$ on the left, larger than $1$ on the right, and equal to $1$ over
the flat part of the discontinuous bottom topography.

We compute the solution by the PCCU-5 scheme until the final time $t=1$ on a uniform mesh with $N=200$ cells and present the differences
$\|\Ep(\cdot,1)-\Ep_{\rm eq}\|_\infty$, $\|q(\cdot,1)-q_{\rm eq}\|_\infty$, and $\|\theta(\cdot,1)-\theta_{\rm eq}\|_\infty$ in Tables
\ref{tab62} and \ref{tab63} for the continuous and discontinuous bottom topographies, respectively. As one can see, all of the entries in
both the tables are close to the machine errors, and hence the proposed PCCU-5 scheme can preserve moving-water equilibria in all of the
six studied cases.
\begin{table}[!htb]
\centering\small
\caption{Example \ref{ex62}: $\|\Ep(\cdot,1)-\Ep_{\rm eq}\|_\infty$, $\|q(\cdot,1)-q_{\rm eq}\|_\infty$, and
$\|\theta(\cdot,1)-\theta_{\rm eq}\|_\infty$ for the smooth bottom topography \eref{6.1}.}\label{tab62}
\begin{tabular}{lcccccc}
\hline
Case&&$\|\Ep(\cdot,1)-\Ep_{\rm eq}\|_\infty$&&$\|q(\cdot,1)-q_{\rm eq}\|_\infty$&&$\|\theta(\cdot,1)-\theta_{\rm eq}\|_\infty$\\\hline
Subcritical  &&3.13E-13&&3.55E-14&&1.71E-13\\
Supercritical&&1.42E-12&&3.13E-13&&3.27E-13\\
Transcritical&&1.56E-13&&1.02E-14&&1.42E-13\\
\hline
\end{tabular}
\end{table}
\begin{table}[!htb]
\centering\small
\caption{Example \ref{ex62}: $\|\Ep(\cdot,1)-\Ep_{\rm eq}\|_\infty$, $\|q(\cdot,1)-q_{\rm eq}\|_\infty$, and
$\|\theta(\cdot,1)-\theta_{\rm eq}\|_\infty$ for the discontinuous bottom topography \eref{5.2}.}\label{tab63}
\begin{tabular}{lcccccc}
\hline
Case&&$\|\Ep(\cdot,1)-\Ep_{\rm eq}\|_\infty$&&$\|q(\cdot,1)-q_{\rm eq}\|_\infty$&&$\|\theta(\cdot,1)-\theta_{\rm eq}\|_\infty$\\\hline
Subcritical  &&1.85E-13&&3.20E-14&&1.71E-13\\
Supercritical&&2.56E-12&&3.34E-13&&3.34E-13\\
Transcritical&&1.14E-13&&1.15E-14&&1.35E-13\\
\hline
\end{tabular}
\end{table}
\end{example}

\begin{example}[Small perturbations of the 1-D moving-water equilibria]\label{ex63}
\rm In this example also taken from \cite{ZXX}, we test the ability of the PCCU-5 scheme to capture the propagation of small perturbations
of the moving-water equilibria from Example \ref{ex62}.

We use the same setting as in Example \ref{ex62}, but impose the free boundary conditions and add a small perturbation to the initial water
depth and take the following initial data:
$$
h(x,0)=h_{\rm eq}(x)+\left\{\begin{aligned}&0.0001,&&x\in[5.75,6.25],\\&0,&&\mbox{otherwise},\end{aligned}\right.\quad
q(x,0)=q_{\rm eq}(x),\quad\theta(x,0)=\theta_{\rm eq}(x).
$$

We compute the solutions by the PCCU-5 and PCCU-2 schemes in the subcritical and transcritical cases until the final time $t=0.75$, and in
the supercritical case, the final time is $t=0.45$. We use a uniform mesh with $N=200$ cells and compute the reference solutions by the
PCCU-5 scheme on a finer mesh with $N=10000$ cells. Figures \ref{fig61}, \ref{fig62}, and \ref{fig63} display the deviation of the numerical
solutions from the corresponding equilibria in the sub-, super-, and transcritical flows, respectively. These results show that the PCCU-5
scheme captures the propagation of small perturbations in a non-oscillatory manner and also clearly outperforms its second-order
PCCU-2 counterpart.
\begin{figure}[!htb]
\centering
\includegraphics[width=0.32\textwidth]{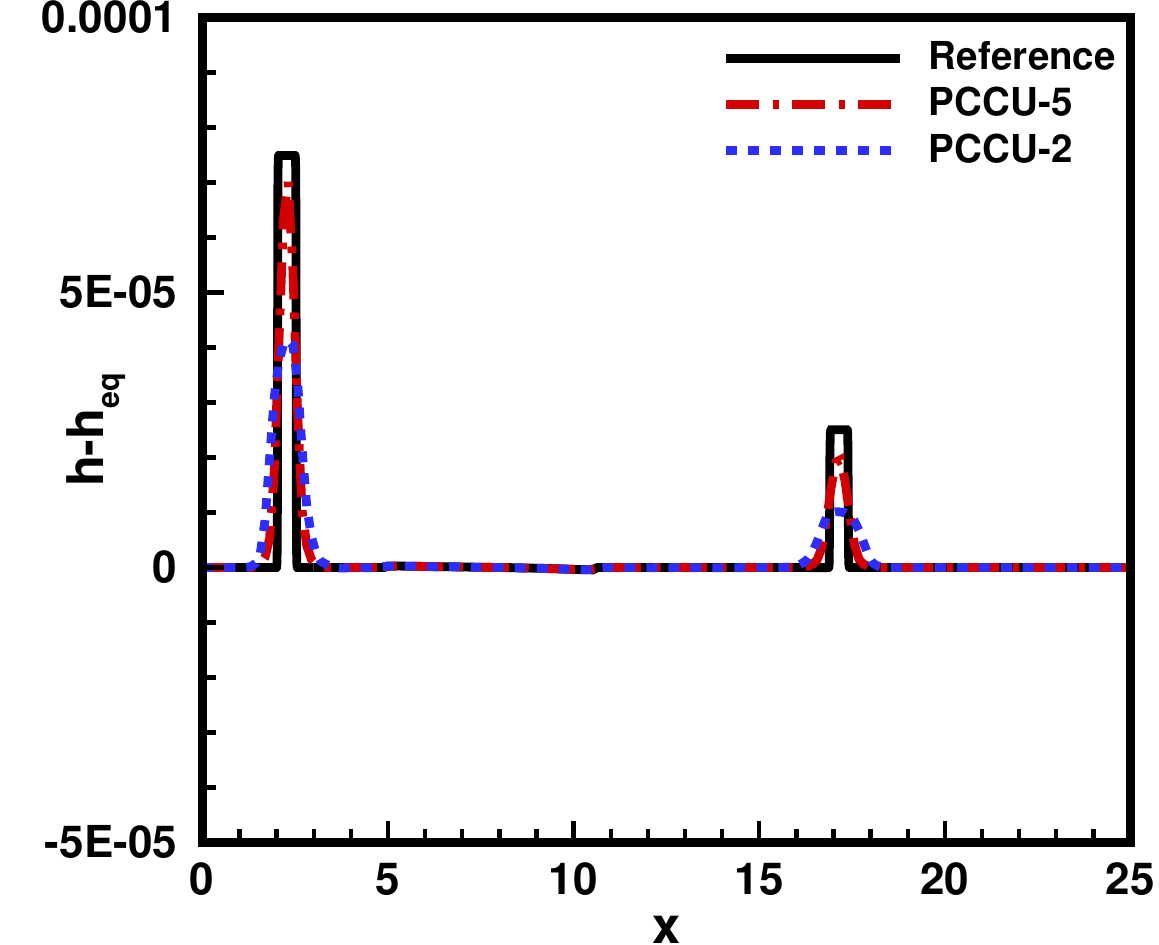}
\includegraphics[width=0.32\textwidth]{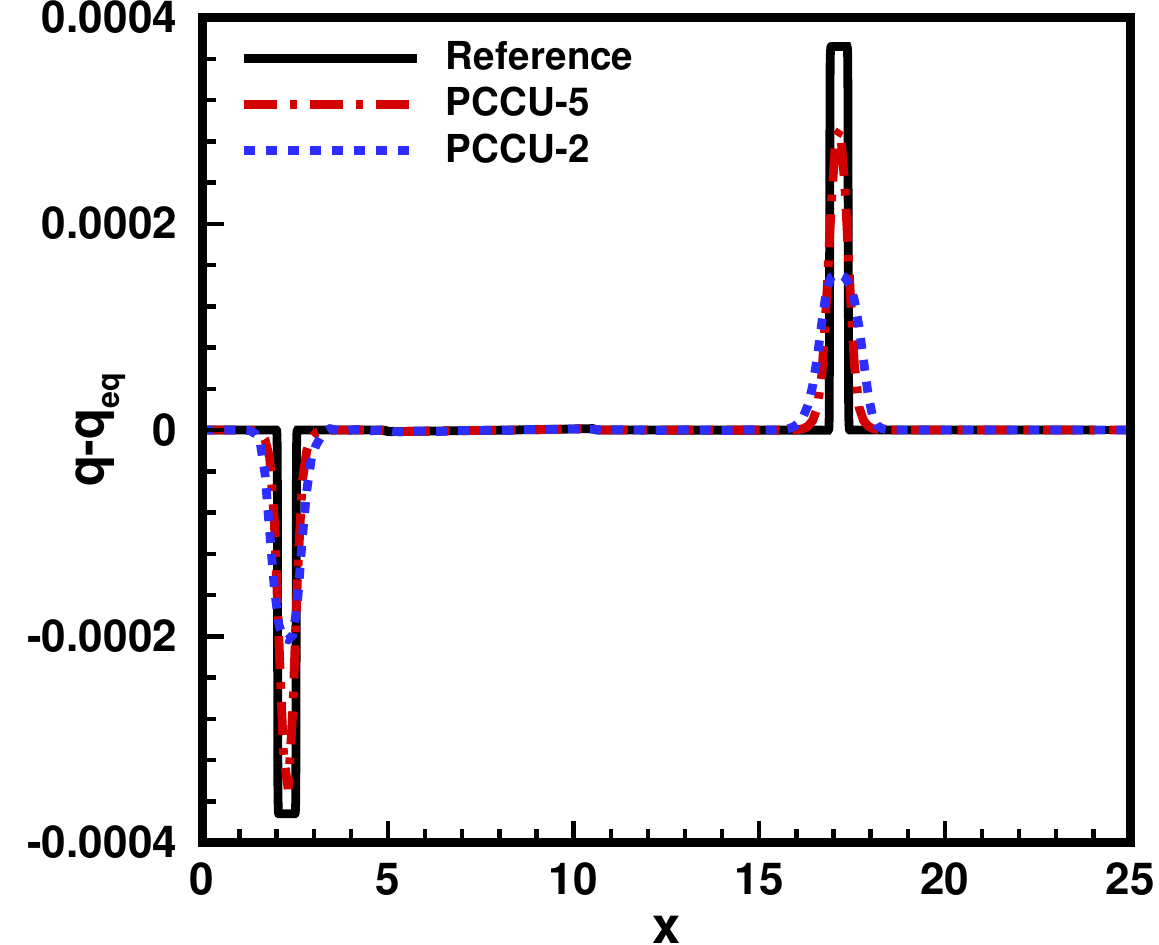}
\includegraphics[width=0.32\textwidth]{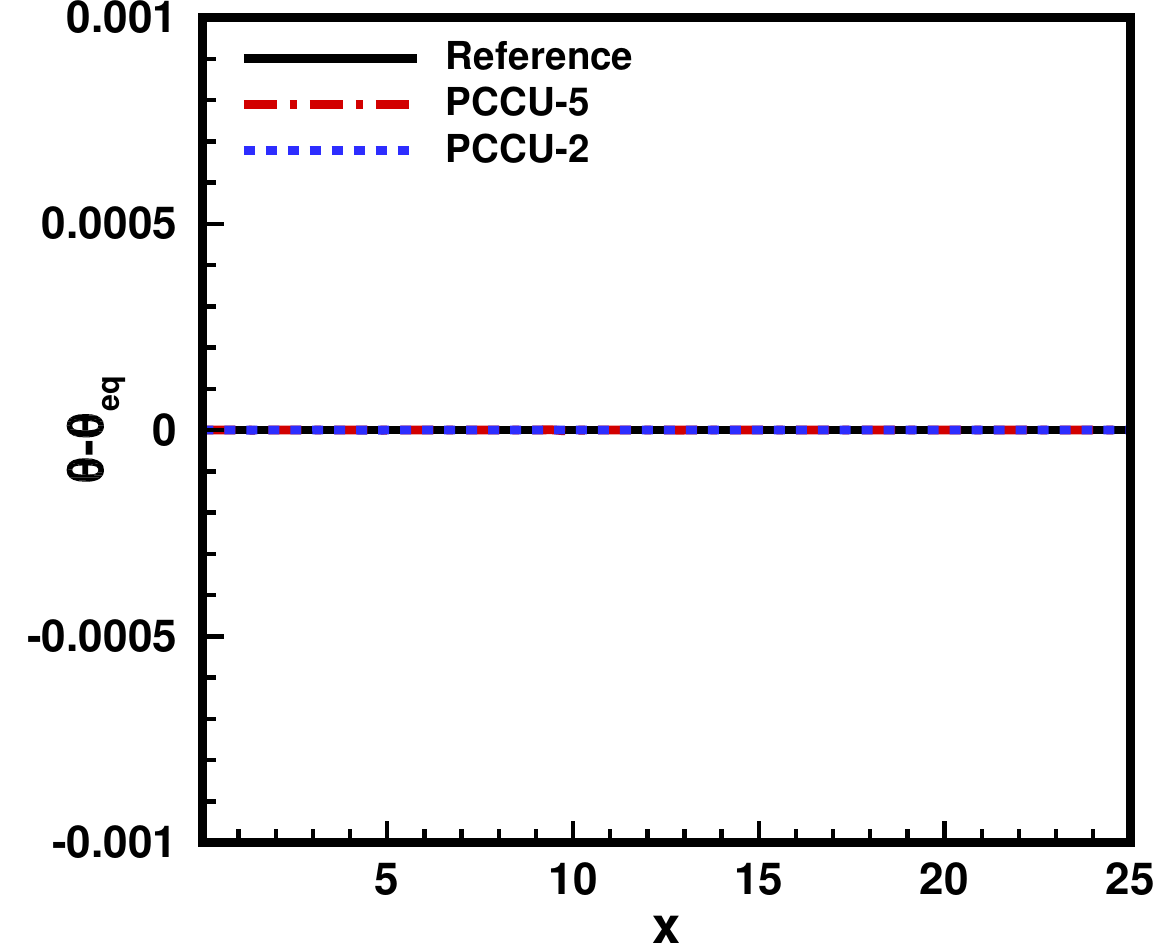}
\vskip-1mm
\caption{Example \ref{ex63} (subcritical case): $h-h_{\rm eq}$, $q-q_{\rm eq}$, and $\theta-\theta_{\rm eq}$ computed by the PCCU-5 and
PCCU-2 schemes.}
\label{fig61}
\end{figure}
\begin{figure}[!htb]
\centering
\includegraphics[width=0.32\textwidth]{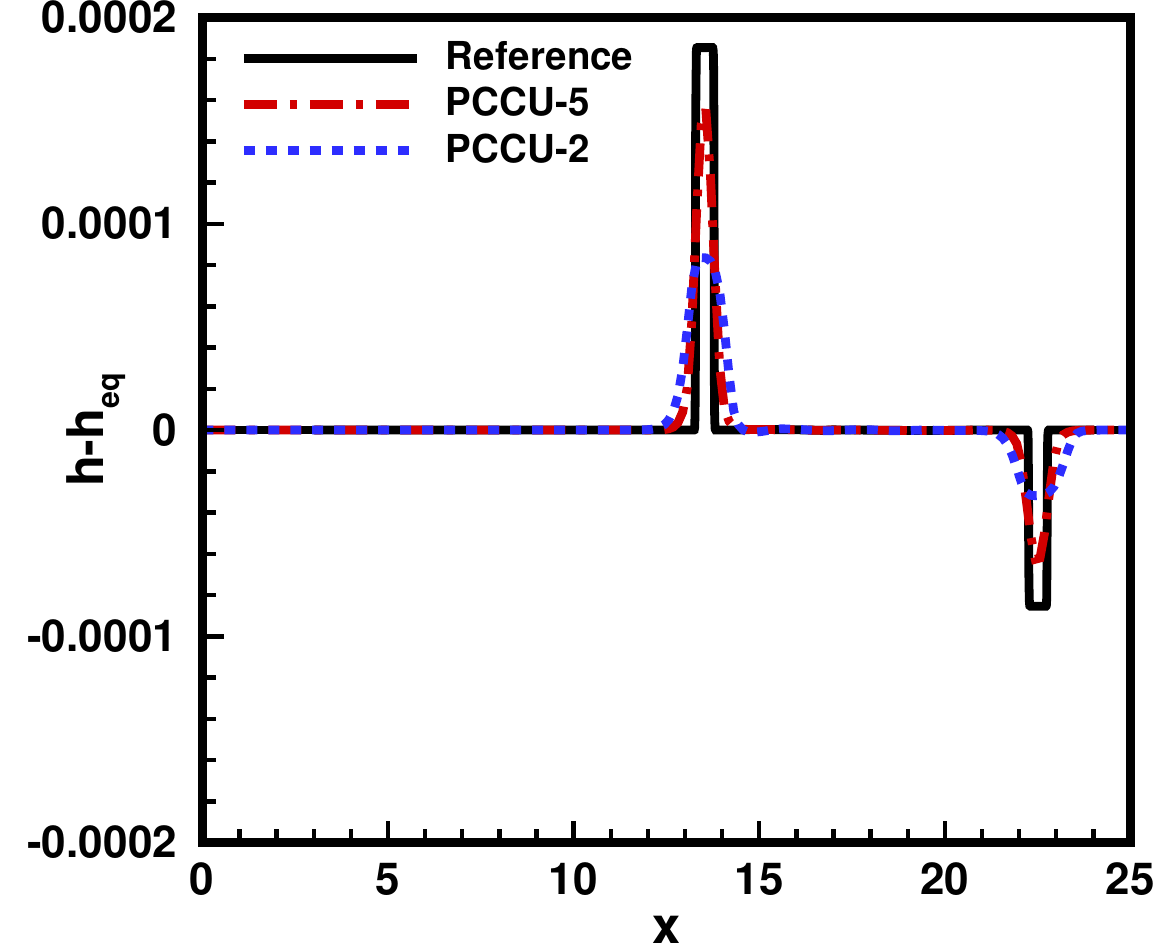}
\includegraphics[width=0.32\textwidth]{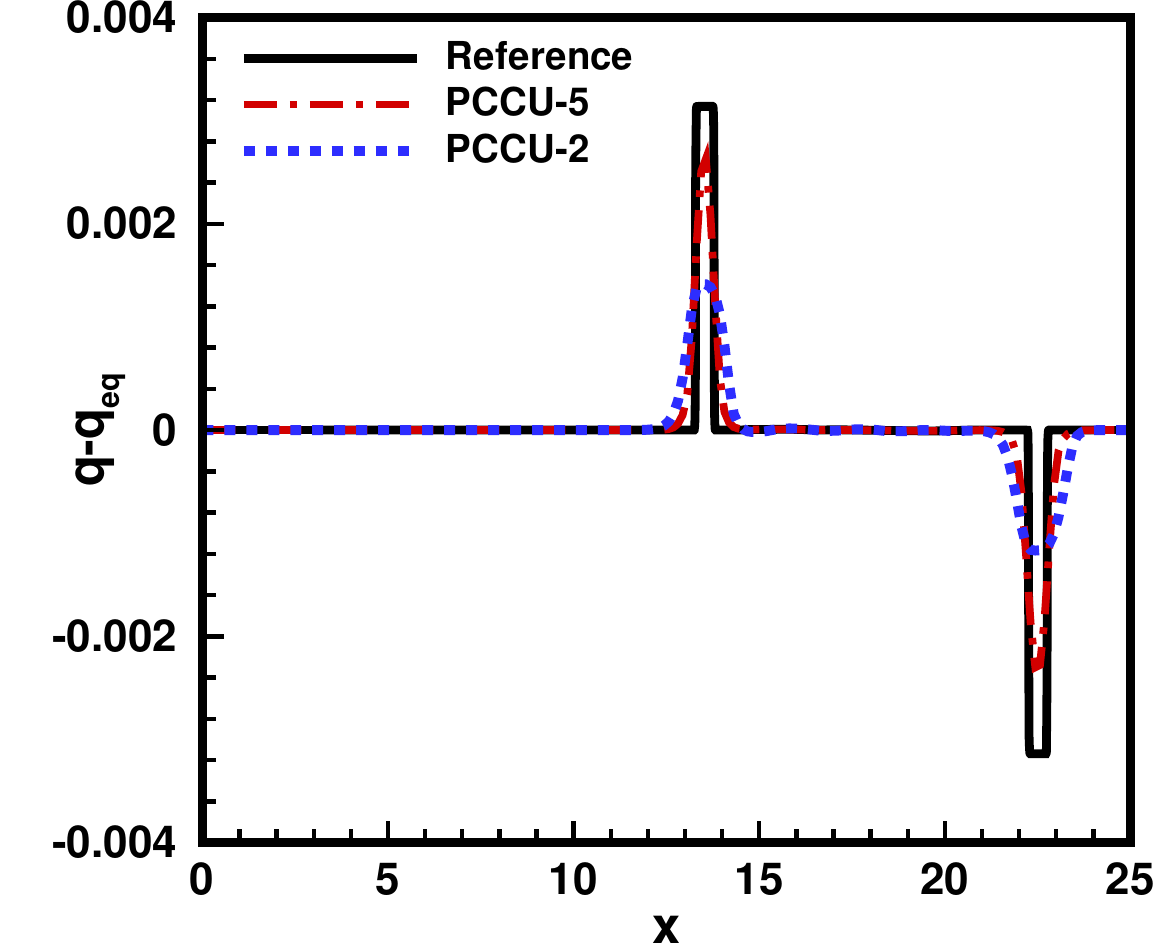}
\includegraphics[width=0.32\textwidth]{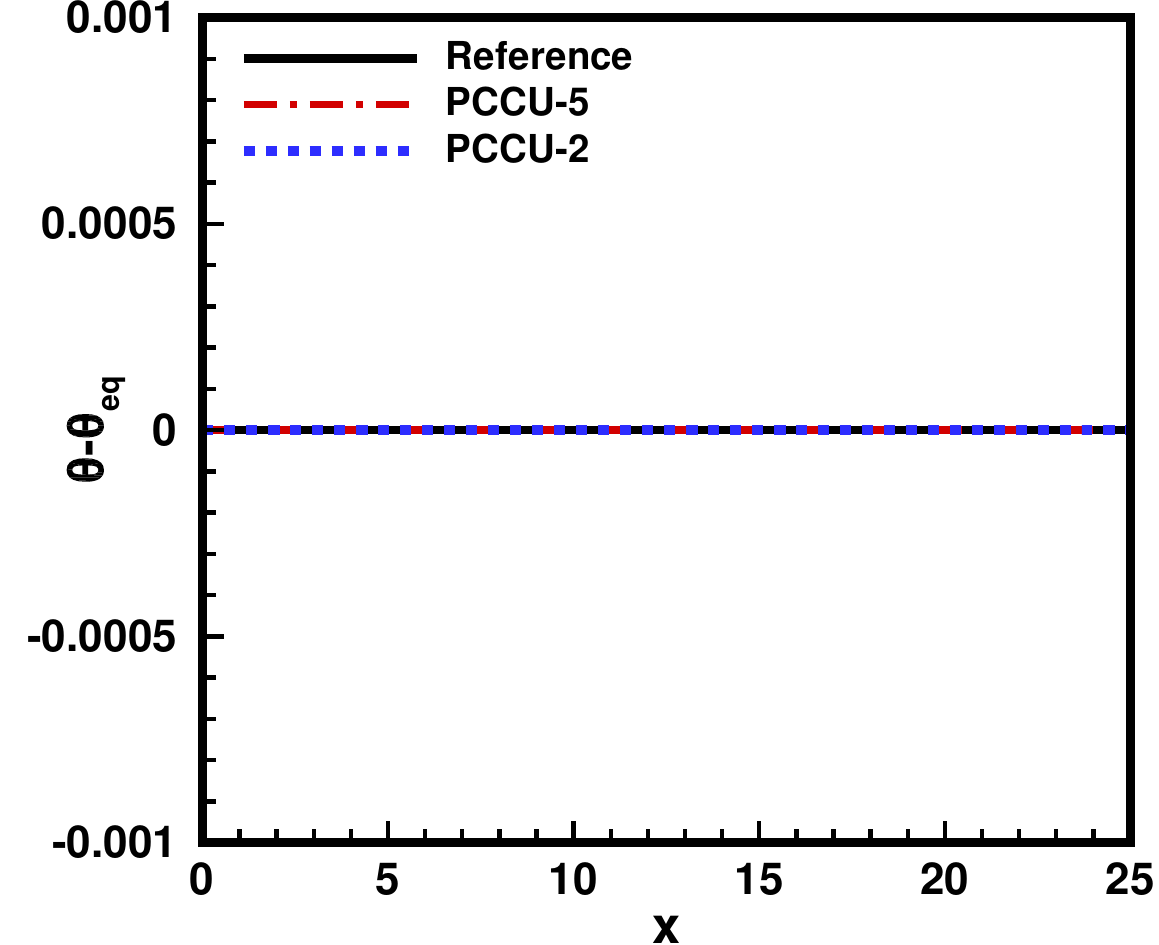}
\vskip-1mm
\caption{Same as in Figure \ref{fig61}, but in the supercritical case.}
\label{fig62}
\end{figure}
\begin{figure}[!htb]
\centering
\includegraphics[width=0.32\textwidth]{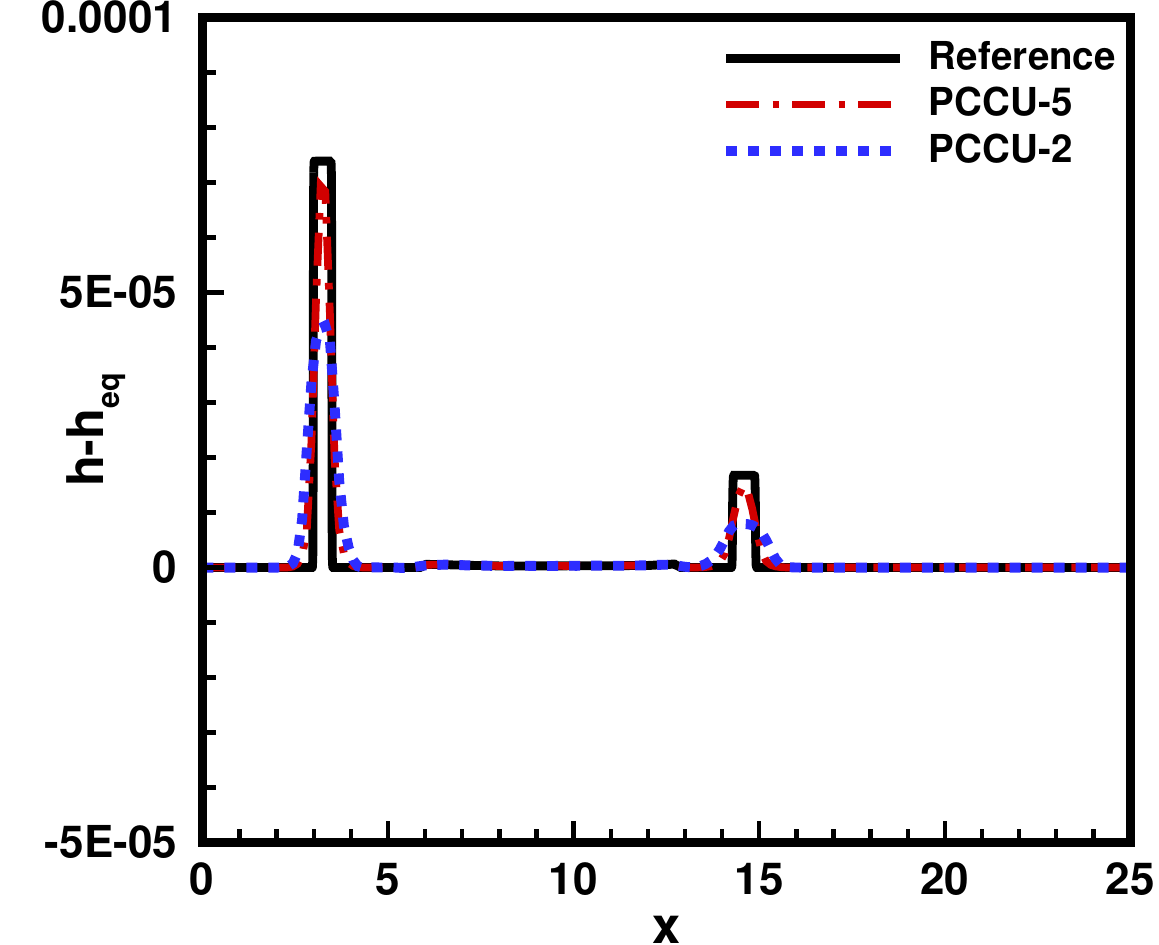}
\includegraphics[width=0.32\textwidth]{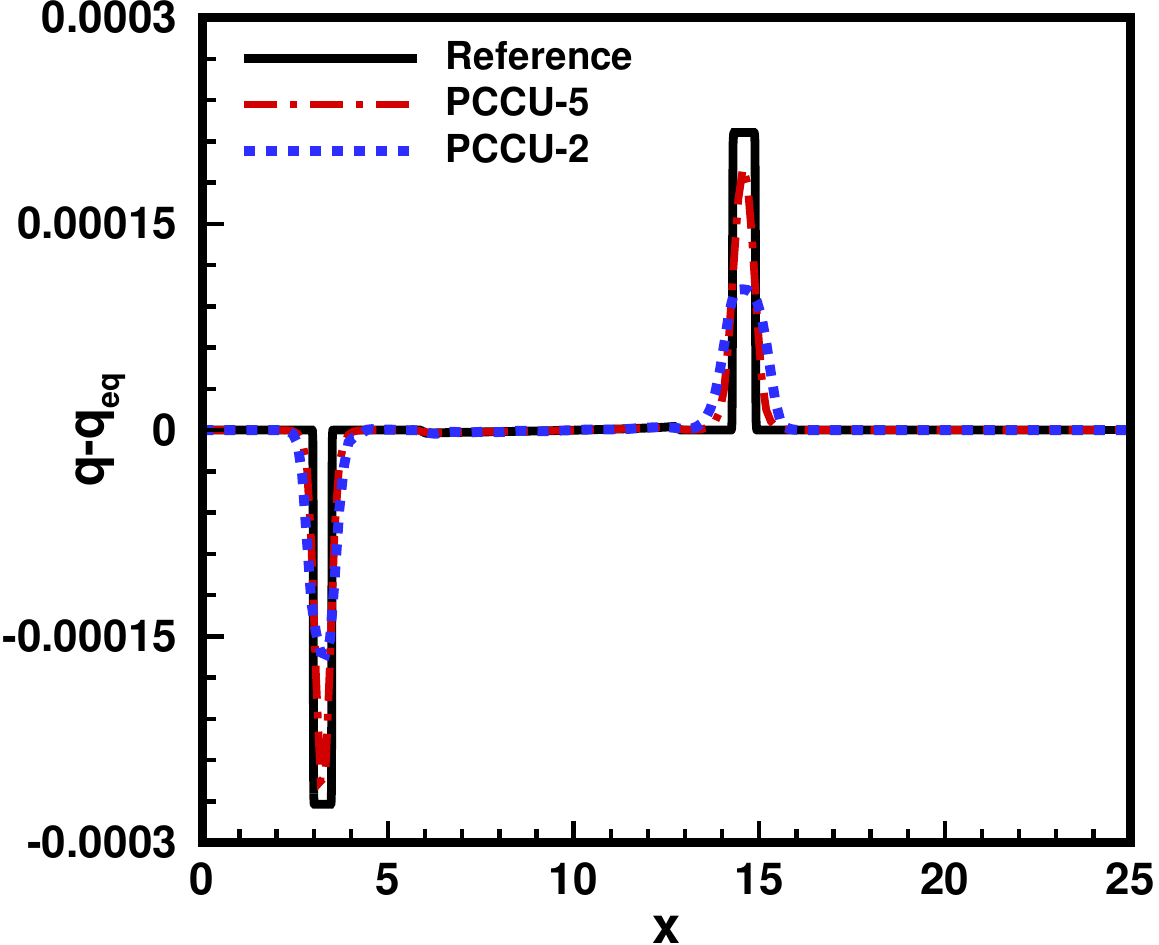}
\includegraphics[width=0.32\textwidth]{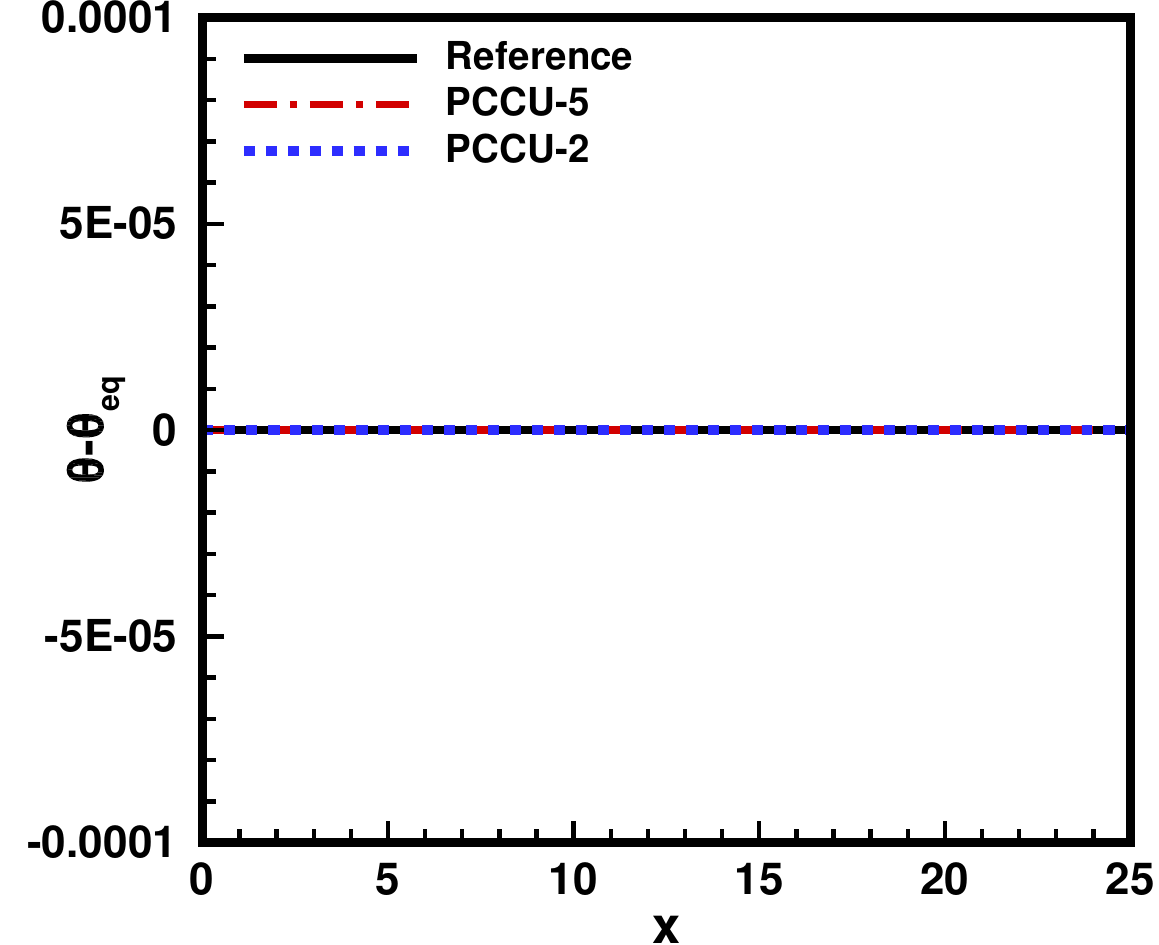}
\vskip-1mm
\caption{Same as in Figure \ref{fig61}, but in the transcritical case.}
\label{fig63}
\end{figure}
\end{example}

\begin{example}[Dam-break problems]\label{ex64}
\rm In this example, the following Riemann initial data:
\begin{equation*}
(h,u,\theta)\Big|_{(x,0)}=\left\{\begin{aligned}
&(5,0.5,9.812),&&-0.5\le x\le0.5,\\
&(3,2.75,15.2086),&&\quad~~~\mbox{otherwise},
\end{aligned}\right.
\end{equation*}
are prescribed in the computational domain $[-1,1]$ subject to the free boundary conditions. We consider three different bottom
topographies: a flat, $Z(x)\equiv0$, smooth,
$$
Z(x)=\left\{\begin{aligned}
&0.50(1-\cos(10\pi x)),&&-0.4\le x\le-0.2, \\
&0.75(1-\cos(10\pi x)),&&\quad~0.2\le x\le0.4,\\
&0,&&\quad~~~\mbox{otherwise},
\end{aligned}\right.
$$
and discontinuous
$$
Z(x)=\left\{\begin{aligned}
&0.3,&&-0.3\le x\le0.3,\\
&0,&&\quad~~~\mbox{otherwise},
\end{aligned}\right.
$$
ones.

We compute the solutions by the PCCU-5 and PCCU-5-NCD schemes until the final time $t=0.075$ on a uniform mesh with $N=200$ cells, and the
reference solution by the PCCU-5 scheme on a finer mesh with $N=3000$ cells. We plot the obtained solutions ($h+Z$, $u$, and $h\theta$) in
Figures \ref{fig64}, \ref{fig65}, and \ref{fig66}, where one can see a good agreement between the PCCU-5 and reference solutions. At the
same time, zooming in on several solution parts shows that while the PCCU-5 solution is oscillation-free, the PCCU-5-NCD scheme develops
spurious oscillations near discontinuities. This clearly indicates the necessity of the interpolation being performed in the local
characteristic equilibrium variables.
\begin{figure}[!htb]
\centering
\includegraphics[width=0.32\textwidth]{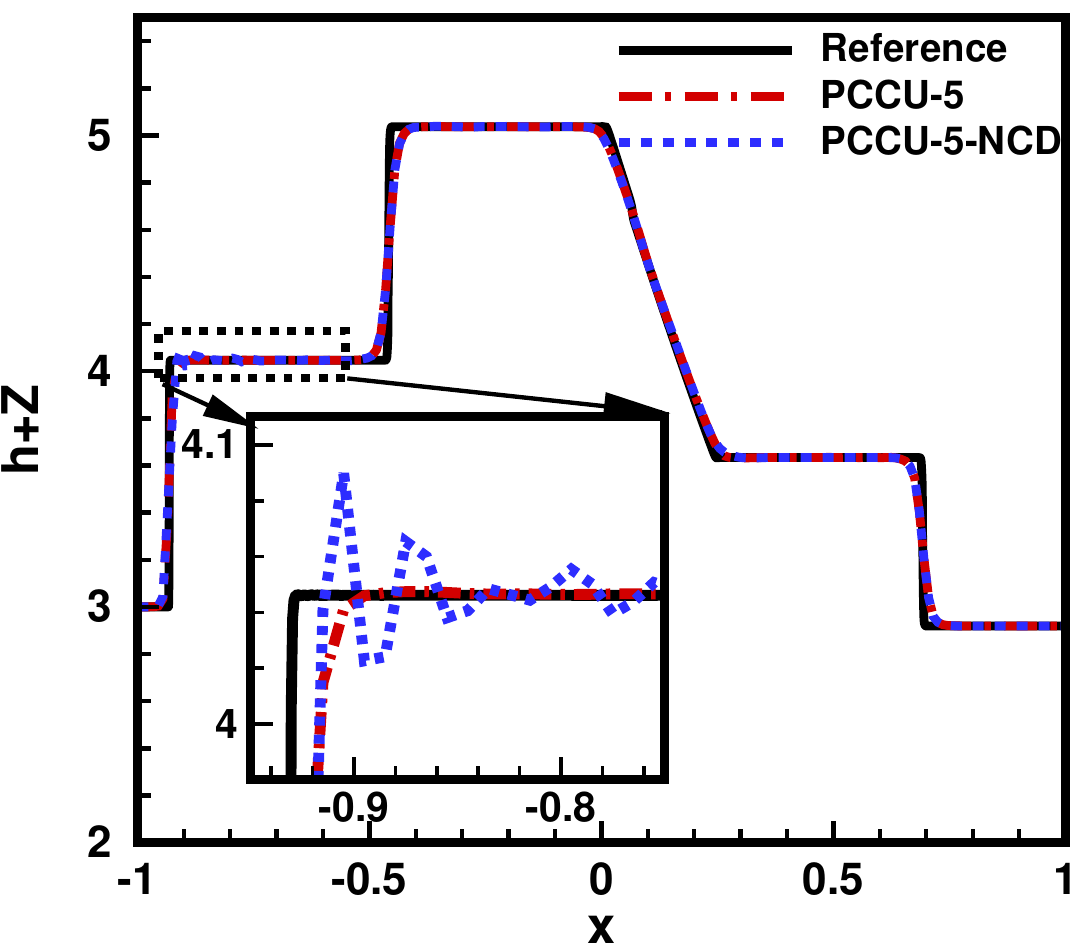}
\includegraphics[width=0.32\textwidth]{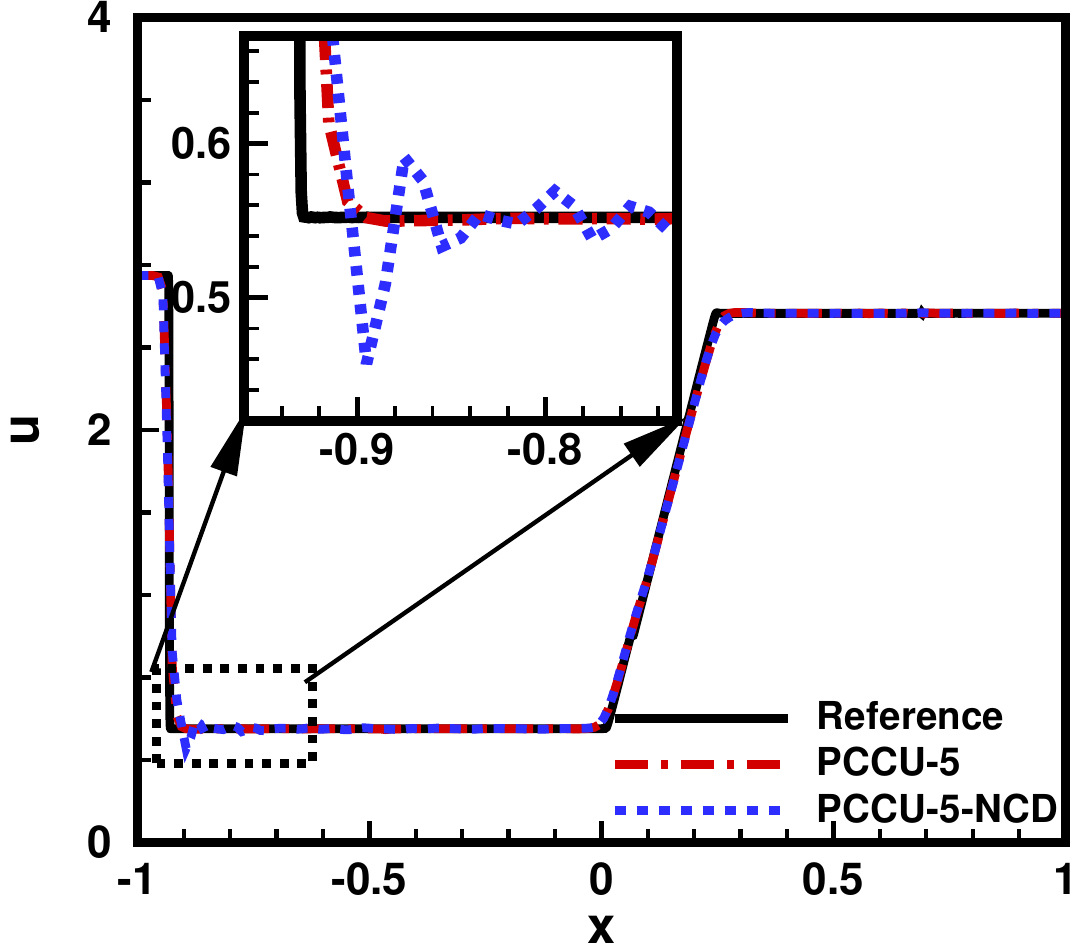}
\includegraphics[width=0.32\textwidth]{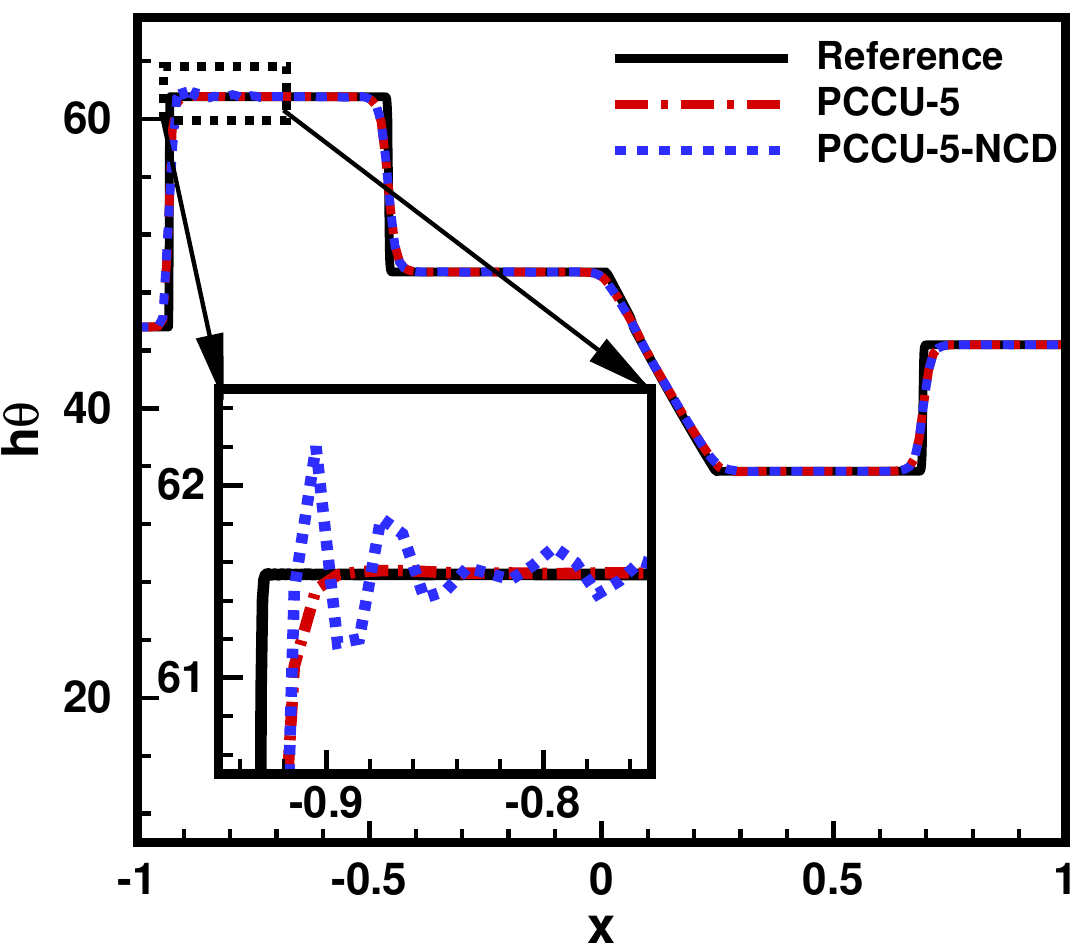}
\vskip-1mm
\caption{Example \ref{ex64} ($Z(x)\equiv0$): $h+Z$, $u$, and $h\theta$ computed by the PCCU-5 and PCCU-5-NCD schemes.}
\label{fig64}
\end{figure}
\begin{figure}[!htb]
\centering
\includegraphics[width=0.32\textwidth]{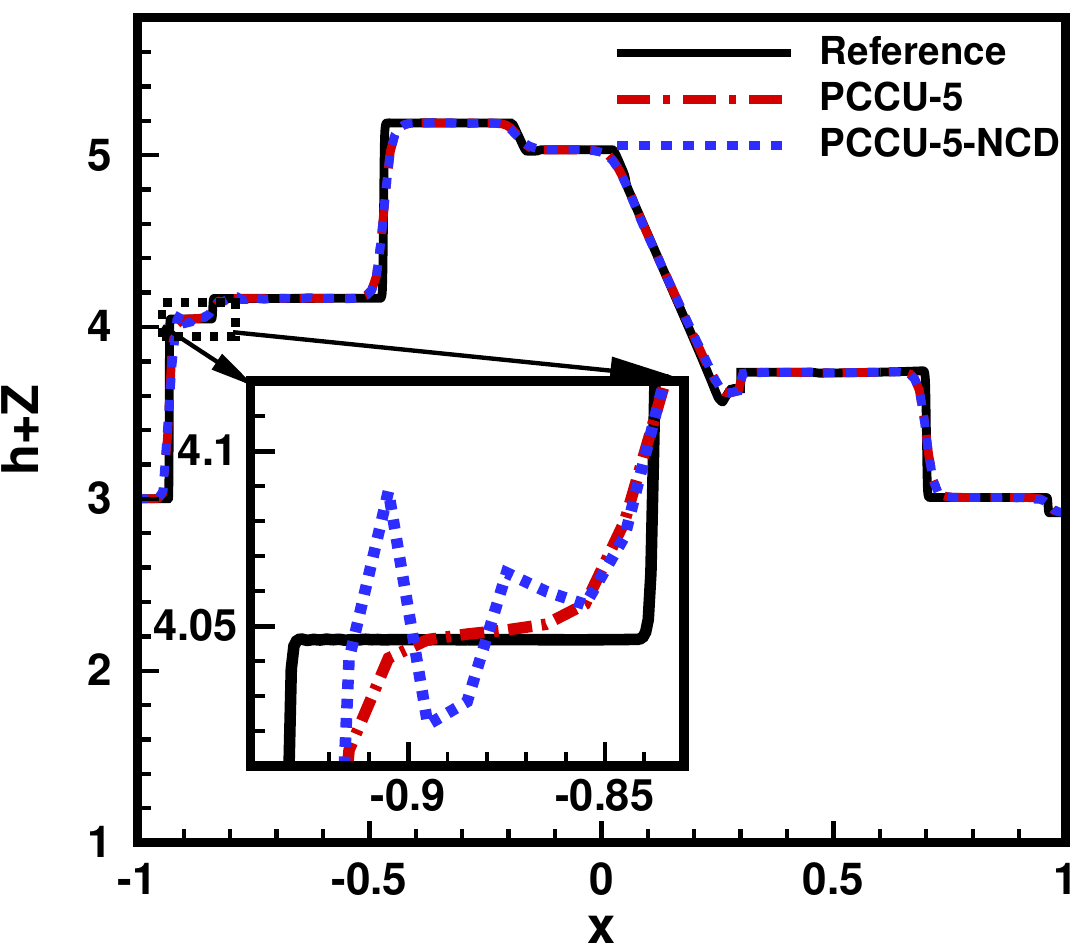}
\includegraphics[width=0.32\textwidth]{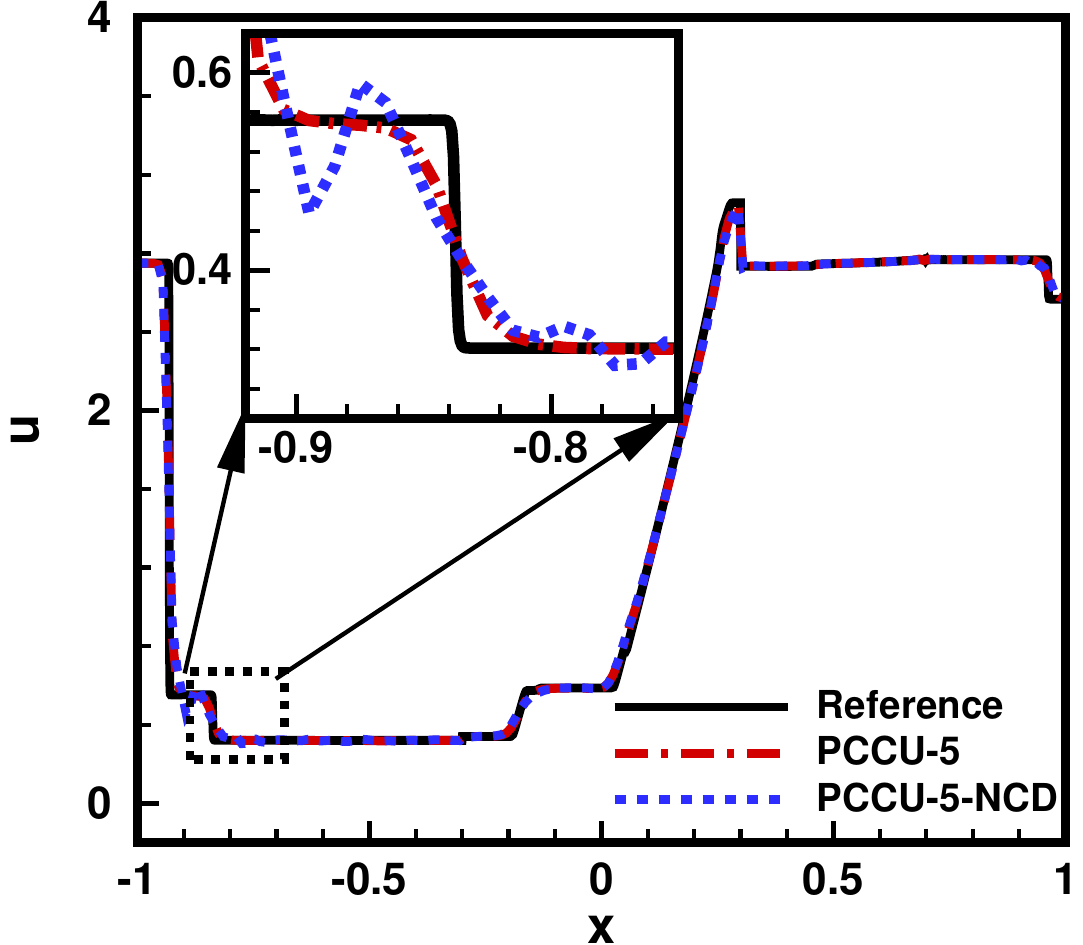}
\includegraphics[width=0.32\textwidth]{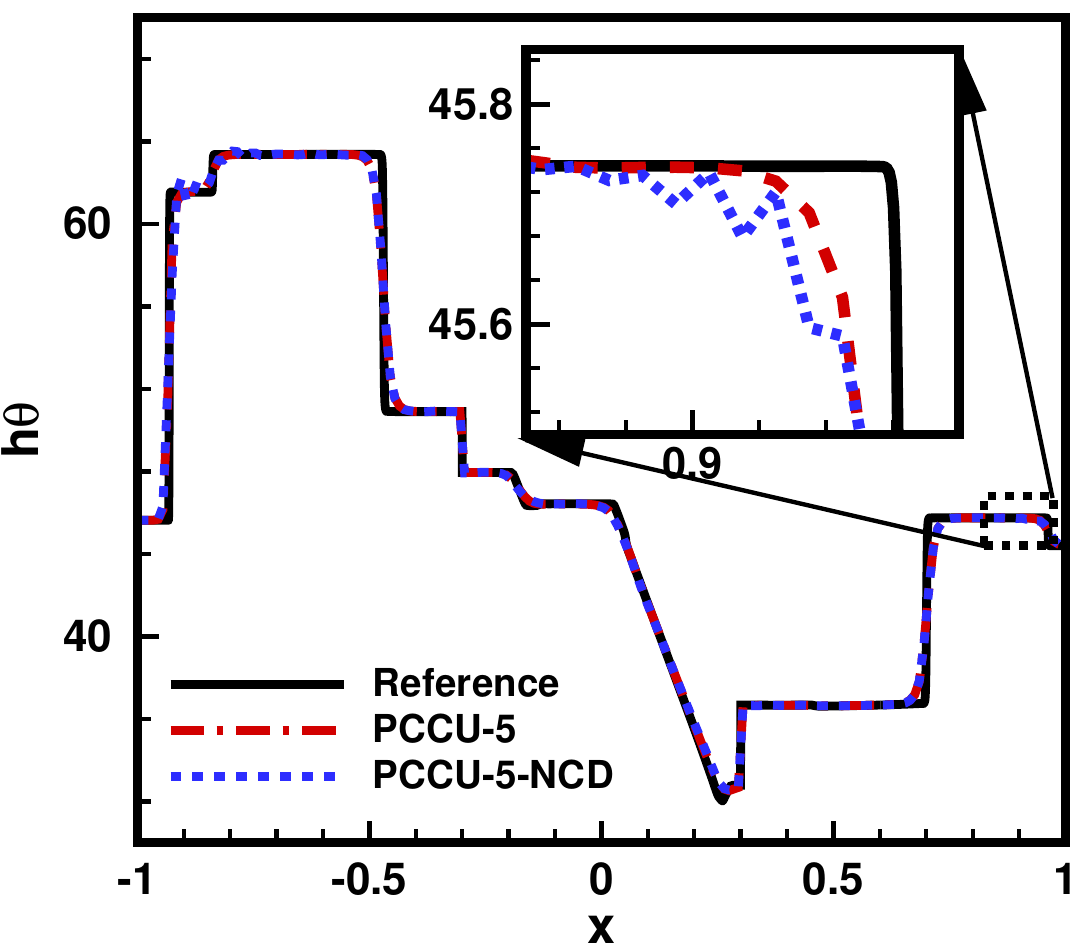}
\vskip-1mm
\caption{Same as in Figure \ref{fig64}, but for the smooth $Z$.}
\label{fig65}
\end{figure}
\begin{figure}[!htb]
\centering
\includegraphics[width=0.32\textwidth]{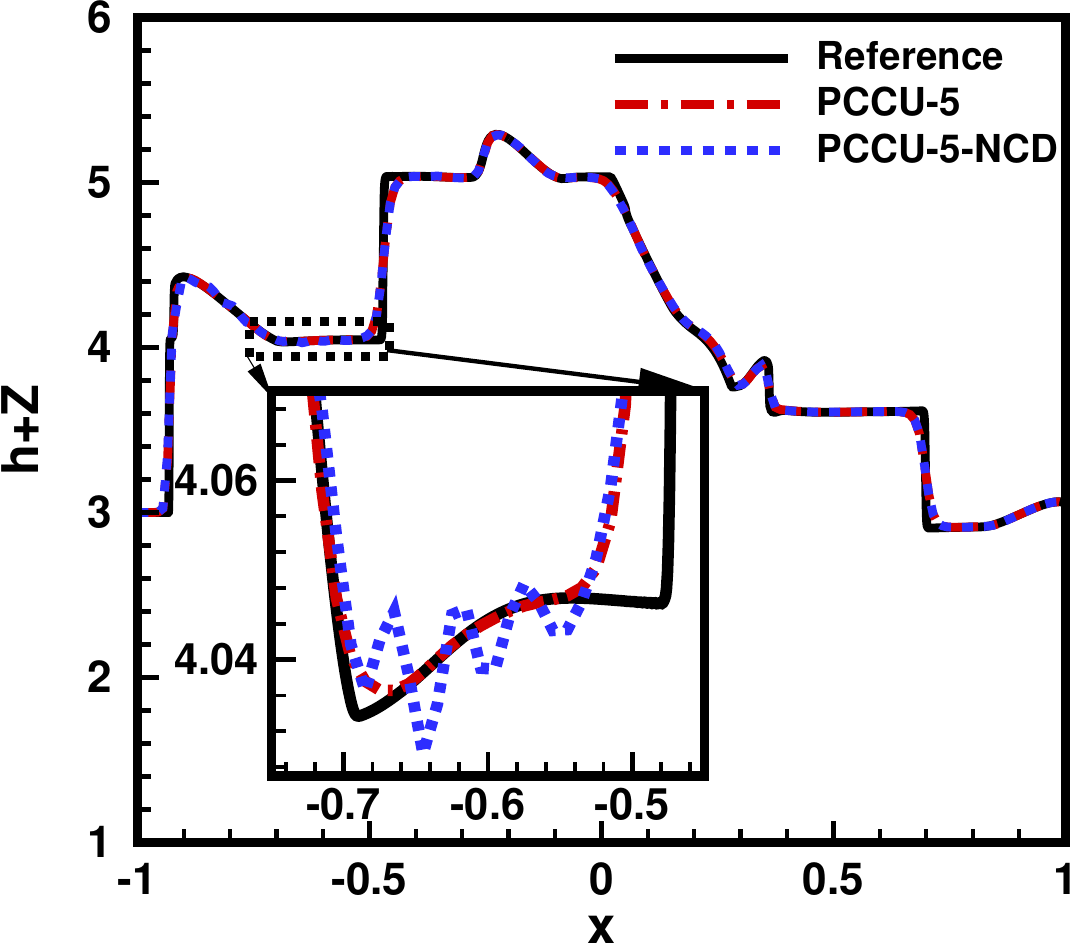}
\includegraphics[width=0.32\textwidth]{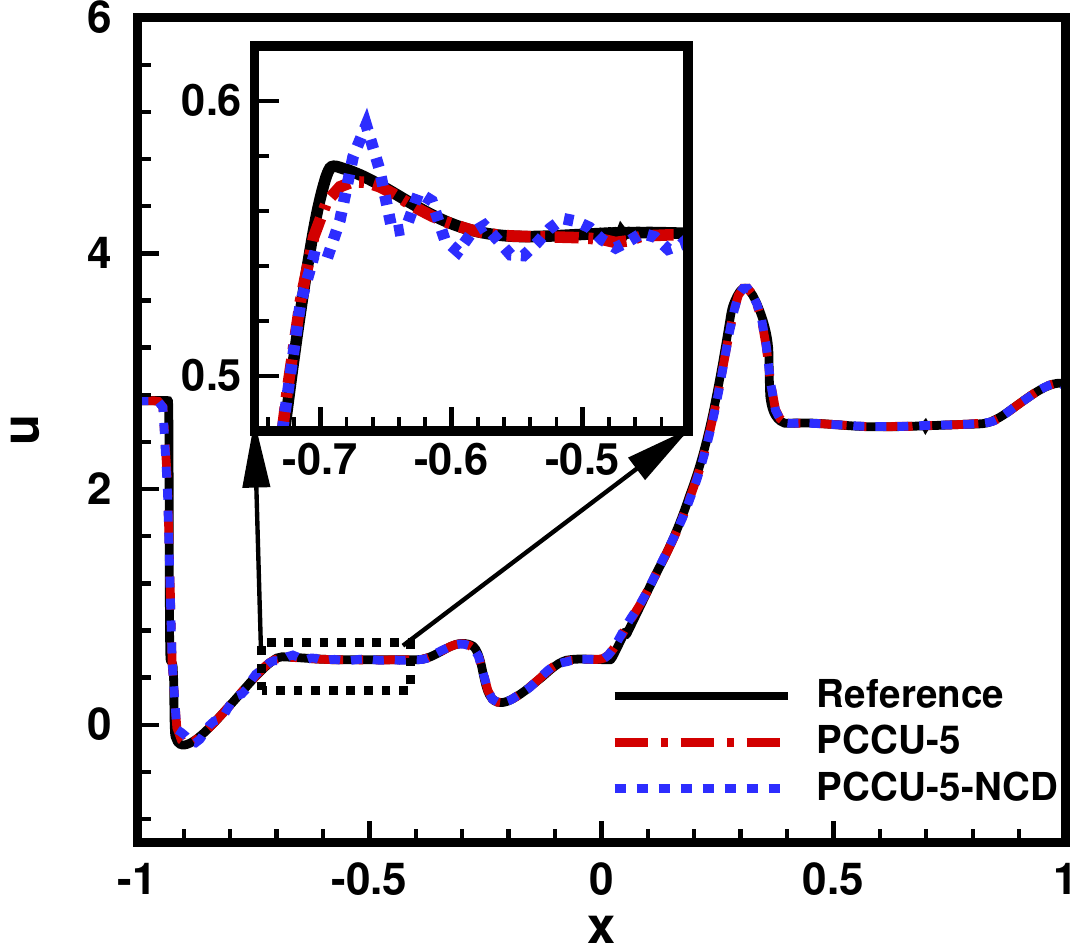}
\includegraphics[width=0.32\textwidth]{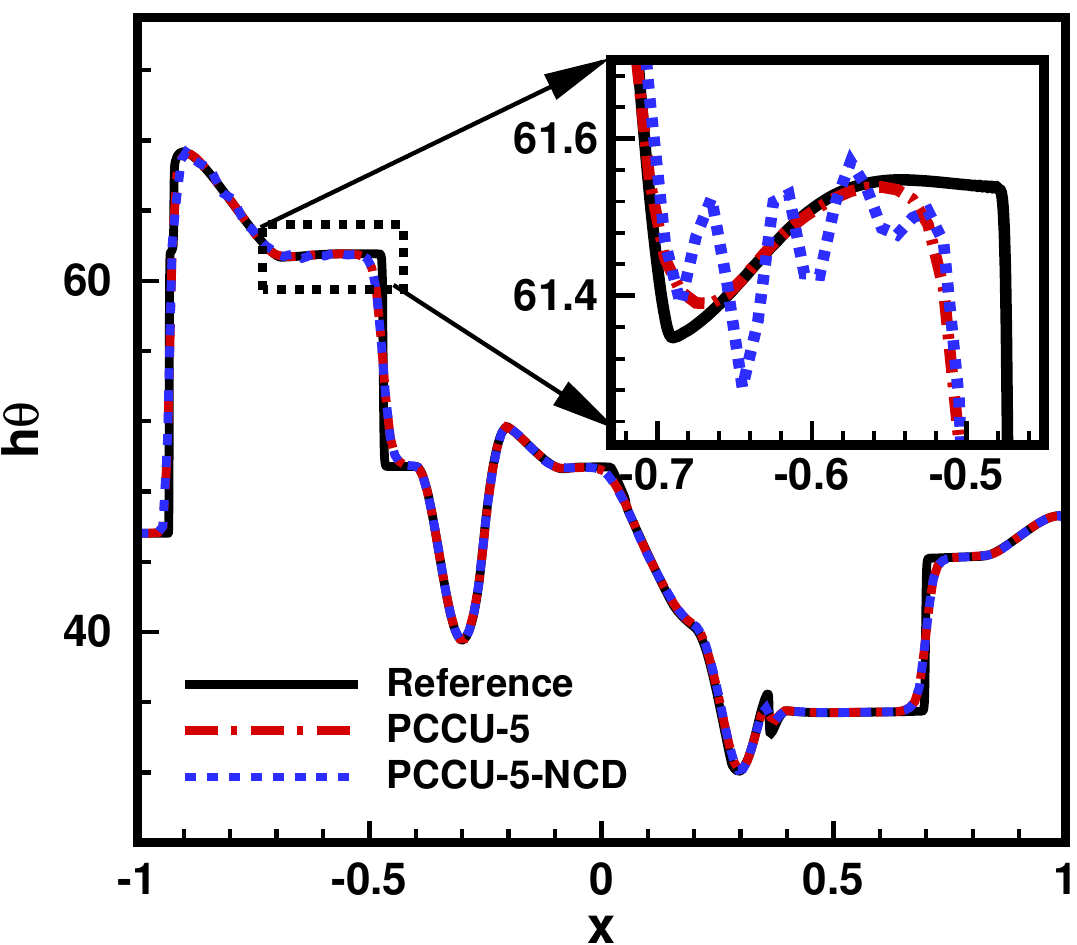}
\vskip-1mm
\caption{Same as in Figure \ref{fig64}, but for the discontinuous $Z$.}
\label{fig66}
\end{figure}
\end{example}

\begin{example}[1-D isobaric equilibrium and its small perturbation]\label{ex65}
\rm The goal of this example is to demonstrate the ability of the PCCU-5 scheme to preserve an isobaric equilibrium and accurately capture
propagation of its small perturbation. The bottom topography is flat ($Z(x)\equiv0$) and the initial conditions, given by
\begin{equation*}
(h,q,P)\Big|_{(x,0)}=\big(h_{\rm eq}(x),q_{\rm eq}(x),P_{\rm eq}(x)\big)
= \left(1+0.0001 e^{-100(x+1.8)^2},0,2\right),
\end{equation*}
are prescribed in the computational domain $[-5,5]$ subject to the free boundary conditions.

We compute the solution by the PCCU-5 scheme until the final time $t=10$ on a uniform grid with $N=200$ and $400$ cells and verify that the
initial steady state is preserved within the machine error; see Table \ref{tab63_2}.
\begin{table}[!htb]
\centering\small
\caption{Example \ref{ex65}: $\|\Ep(\cdot,10)-\Ep_{\rm eq}\|_\infty$, $\|q(\cdot,10)-q_{\rm eq}\|_\infty$,
$\|\theta(\cdot,10)-\theta_{\rm eq}\|_\infty$, and $\|P(\cdot,10)-P_{\rm eq}\|_\infty$.}\label{tab63_2}
\begin{tabular}{lccccccc}
\hline
$N$&&$\|\Ep(\cdot,10)-\Ep_{\rm eq}\|_\infty$&&$\|q(\cdot,10)-q_{\rm eq}\|_\infty$&&$\|\theta(\cdot,10)-\theta_{\rm eq}\|_\infty$&$\|P(\cdot,10)-P_{\rm eq}\|_\infty$ \\ \hline
200  &&1.46E-14&&7.97E-15&&3.91E-14&1.07E-14\\
400  &&1.95E-14&&1.61E-14&&8.66E-14&3.29E-14\\
\hline
\end{tabular}
\end{table}

We then modify the initial conditions by adding a small perturbation to the water depth so that the initial conditions become
\begin{equation*}
h(x,0)=h_{\rm eq}(x)+\left\{\begin{aligned}
&0.0001,&&-0.2< x<0.2,\\
&0,&&\qquad~\mbox{otherwise}.
\end{aligned}\right.
\end{equation*}
We compute the solutions using the PCCU-5 and PCCU-2 schemes until the final time $t=1.6$ on a uniform mesh with $N=200$ cells. In Figure
\ref{fig67}, we displays the difference of $h(x,t)-h_{\rm eq}(x)$ and $P(x,t)-P_{\rm eq}(x)$ at times $t=0.4$, $0.8$, $1.2$, and $1.6$ along
with the reference solution obtained by the PCCU-2 scheme on a much finer mesh consisting of $N=20000$ uniform cells. As one can clearly
see, the PCCU-5 scheme accurately captures the small perturbation of the isobaric equilibrium without producing spurious oscillations, and
it also clearly outperforms its second-order PCCU-2 counterpart, especially when the small perturbation passes the peak of the water depth
$h$ near $x\approx-1.8$.
\begin{figure}[!htb]
\centering
\includegraphics[width=0.24\textwidth]{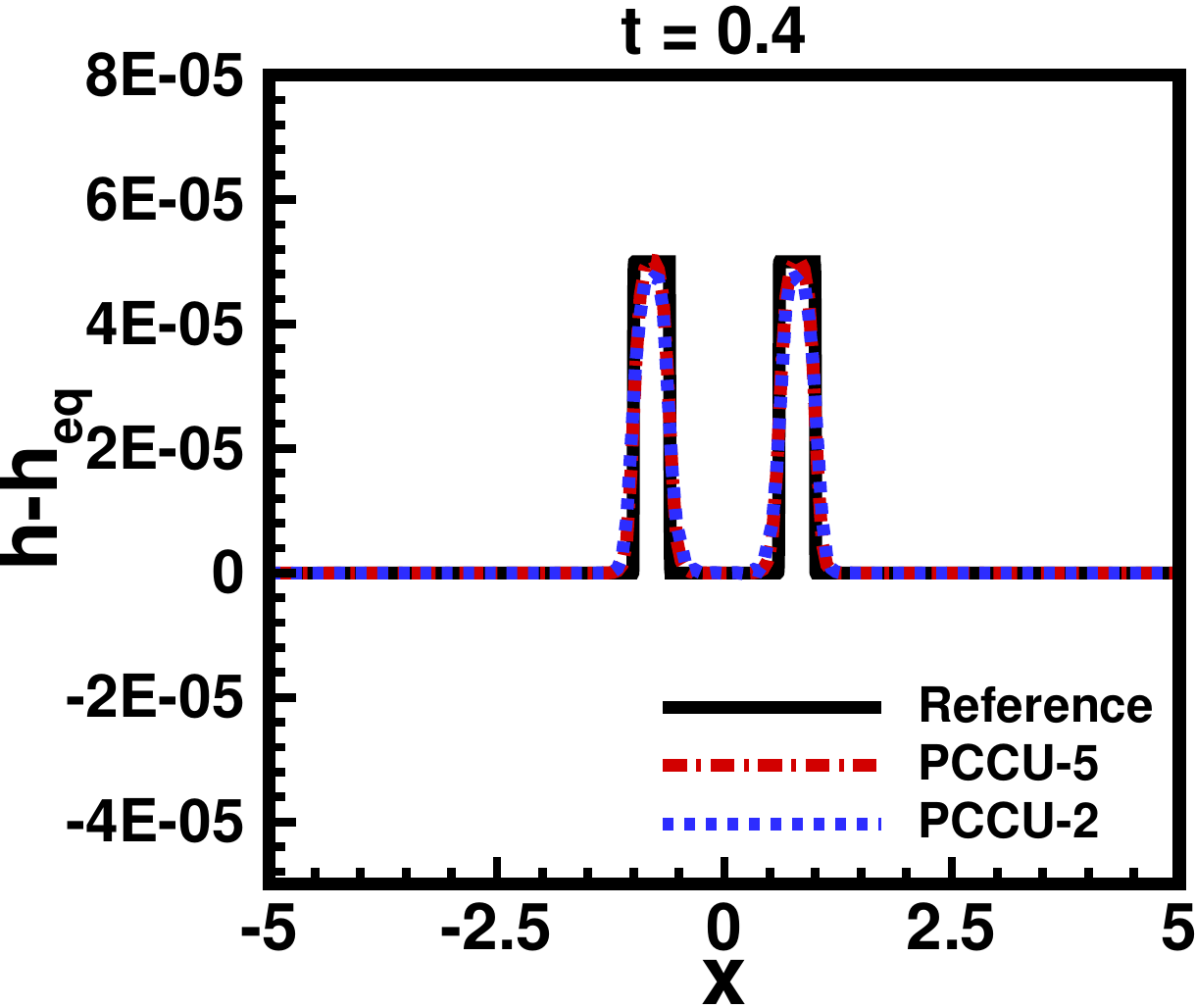}
\includegraphics[width=0.24\textwidth]{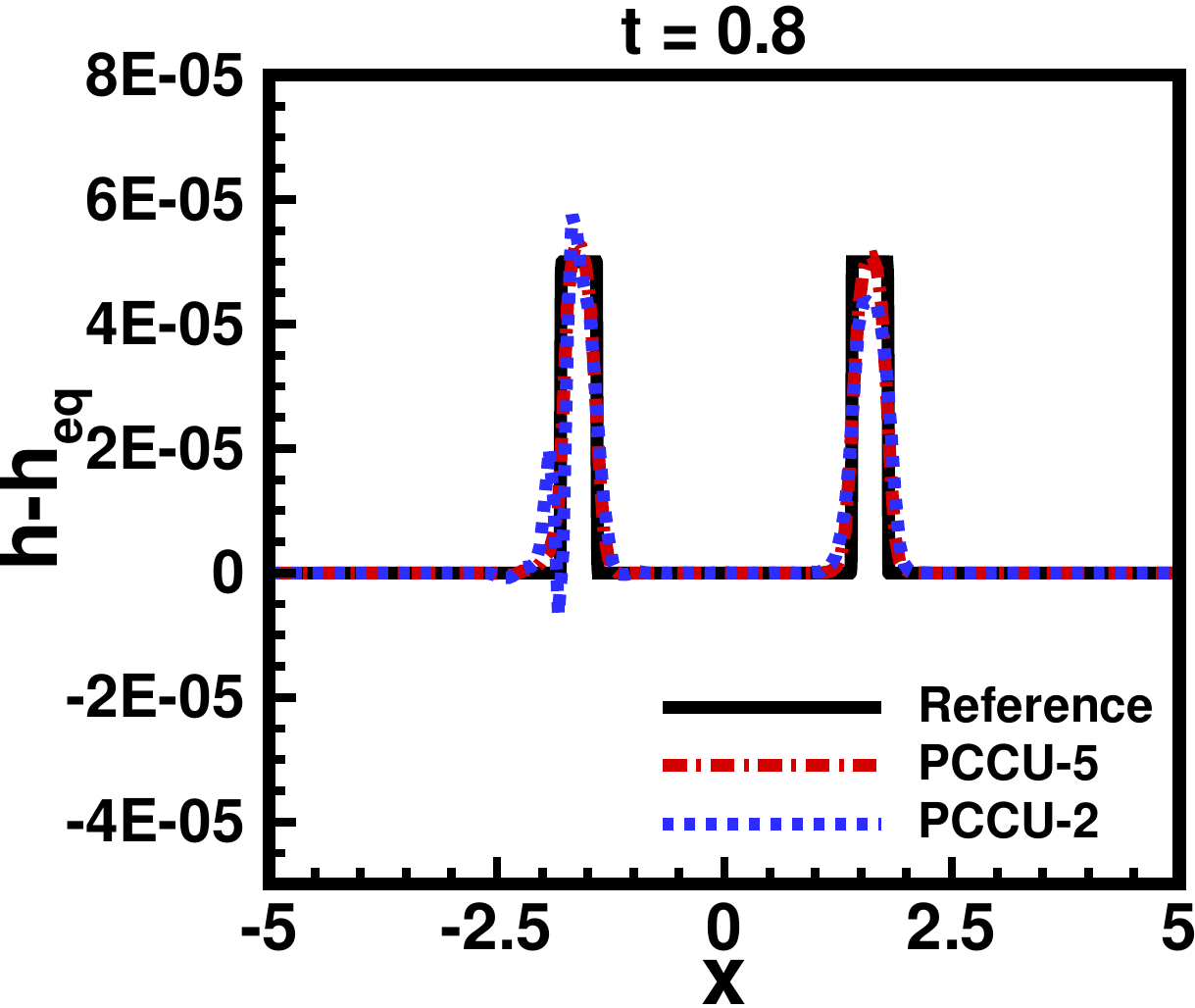}
\includegraphics[width=0.24\textwidth]{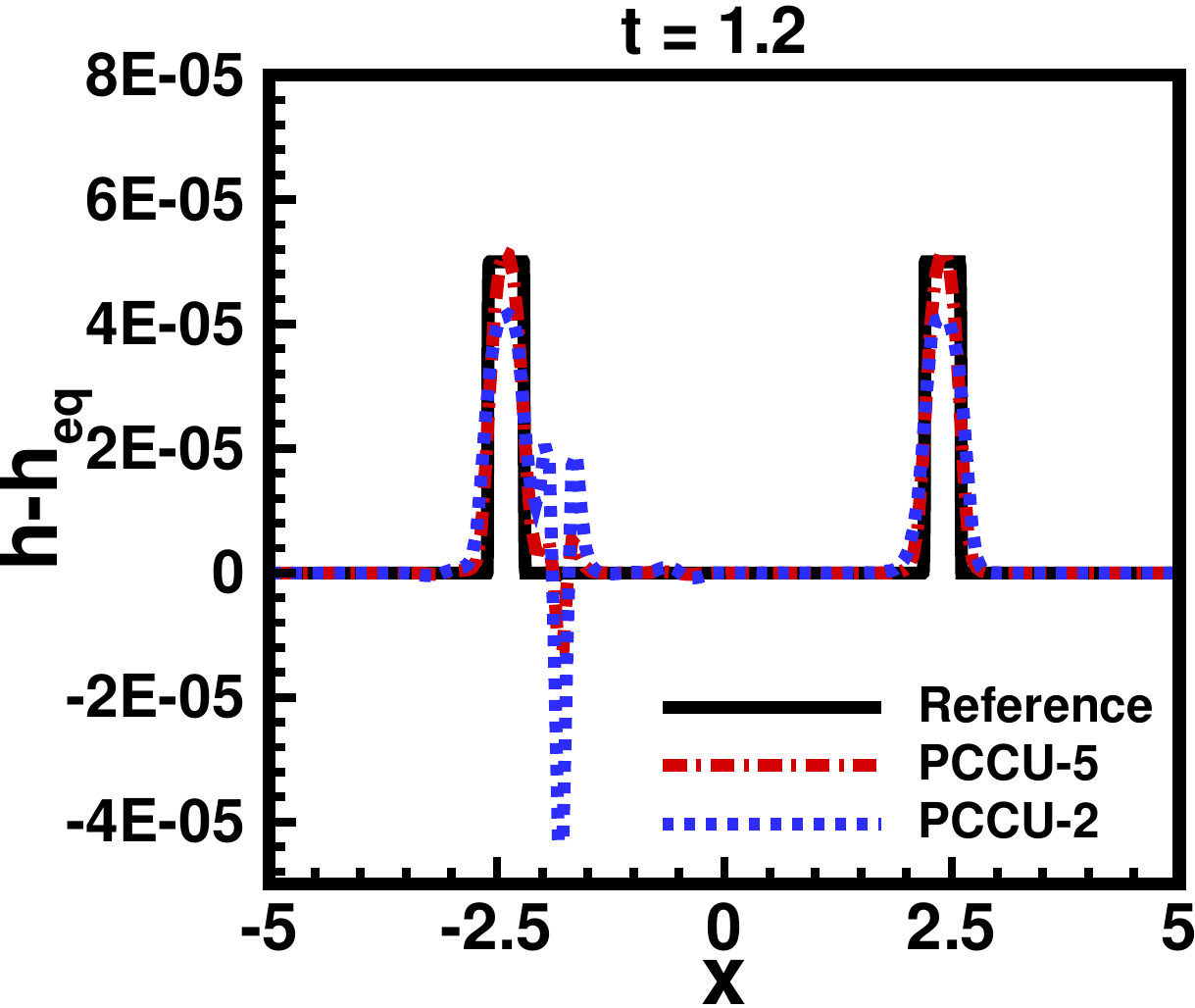}
\includegraphics[width=0.24\textwidth]{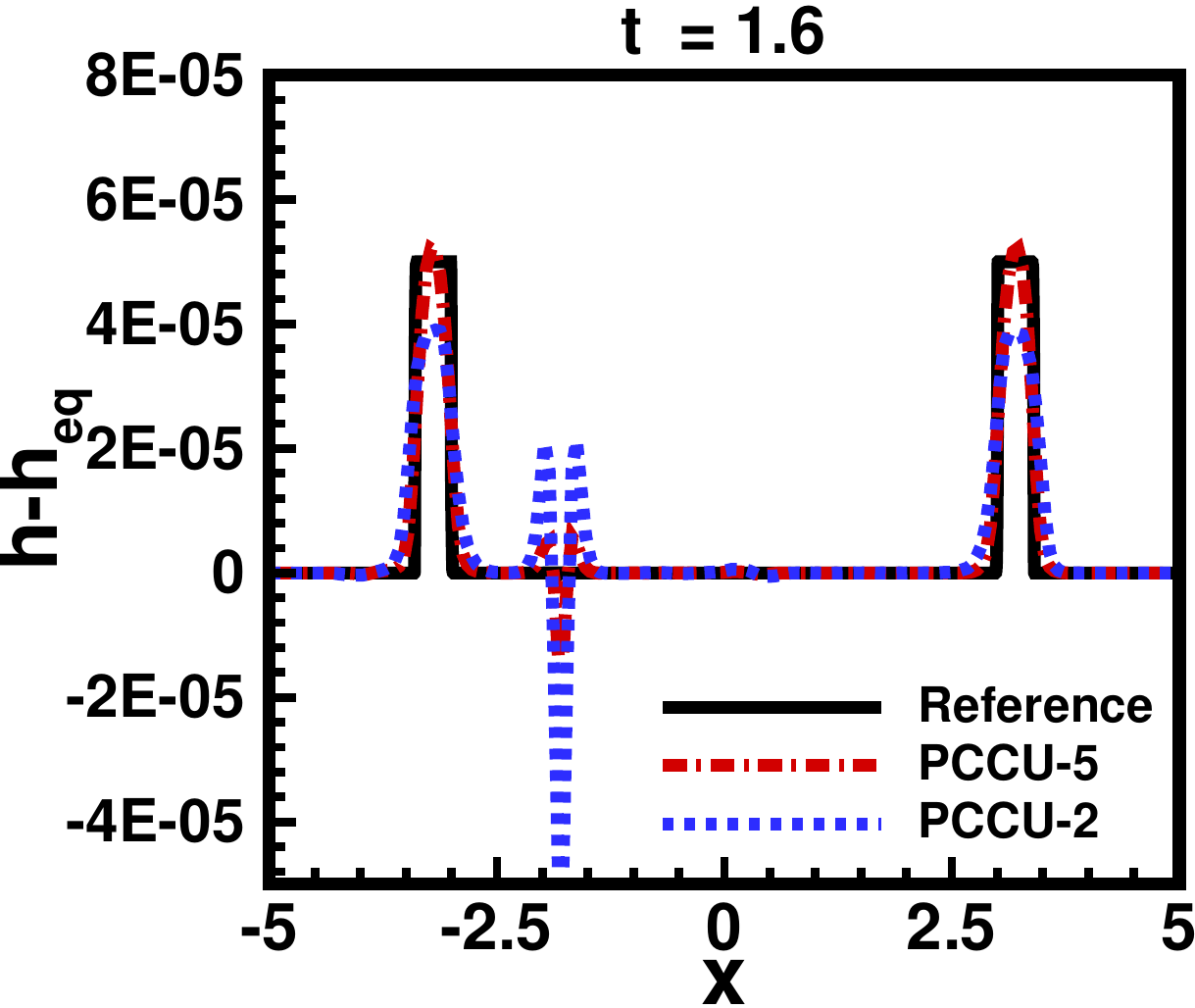}
\vskip5pt
\includegraphics[width=0.24\textwidth]{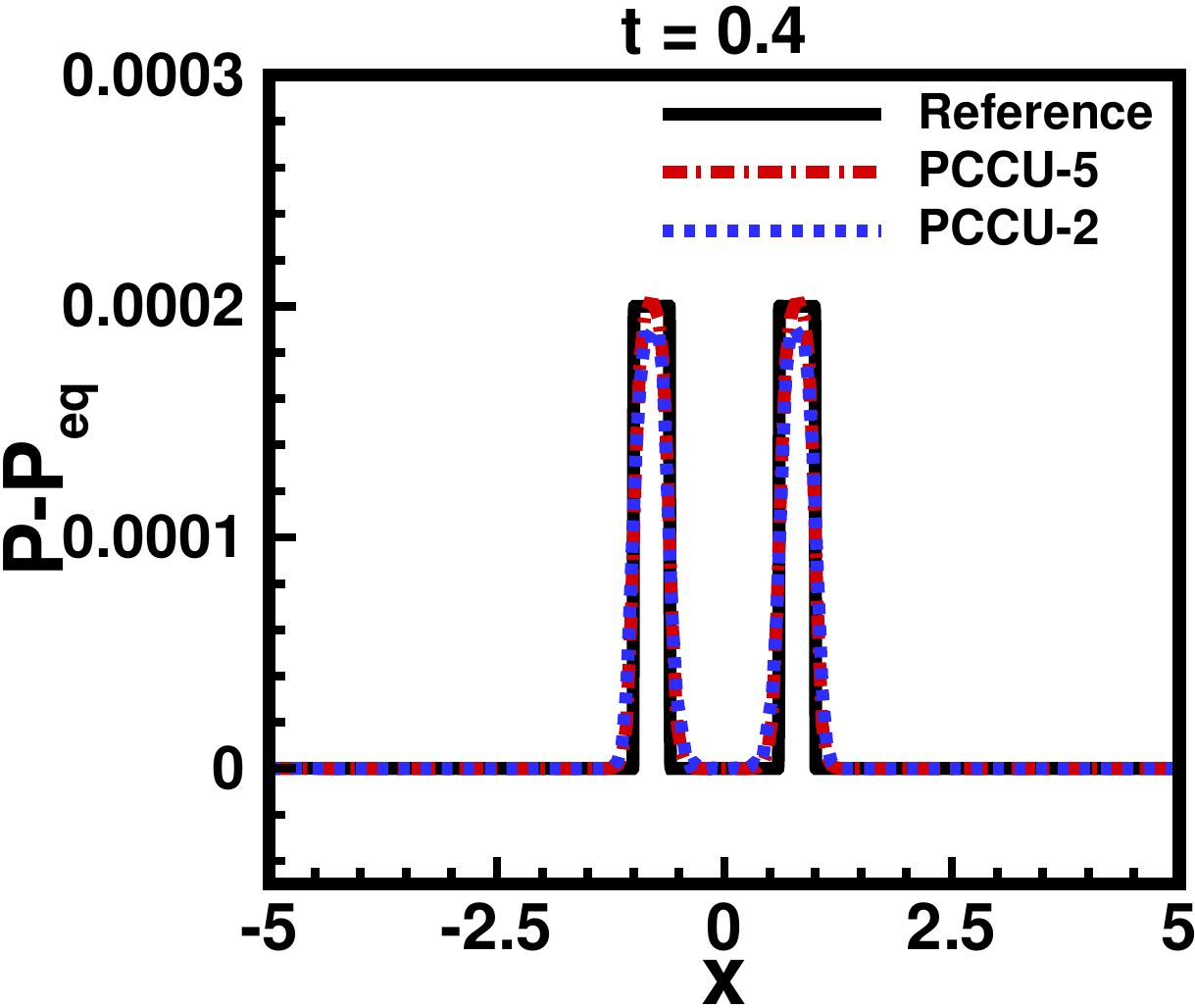}
\includegraphics[width=0.24\textwidth]{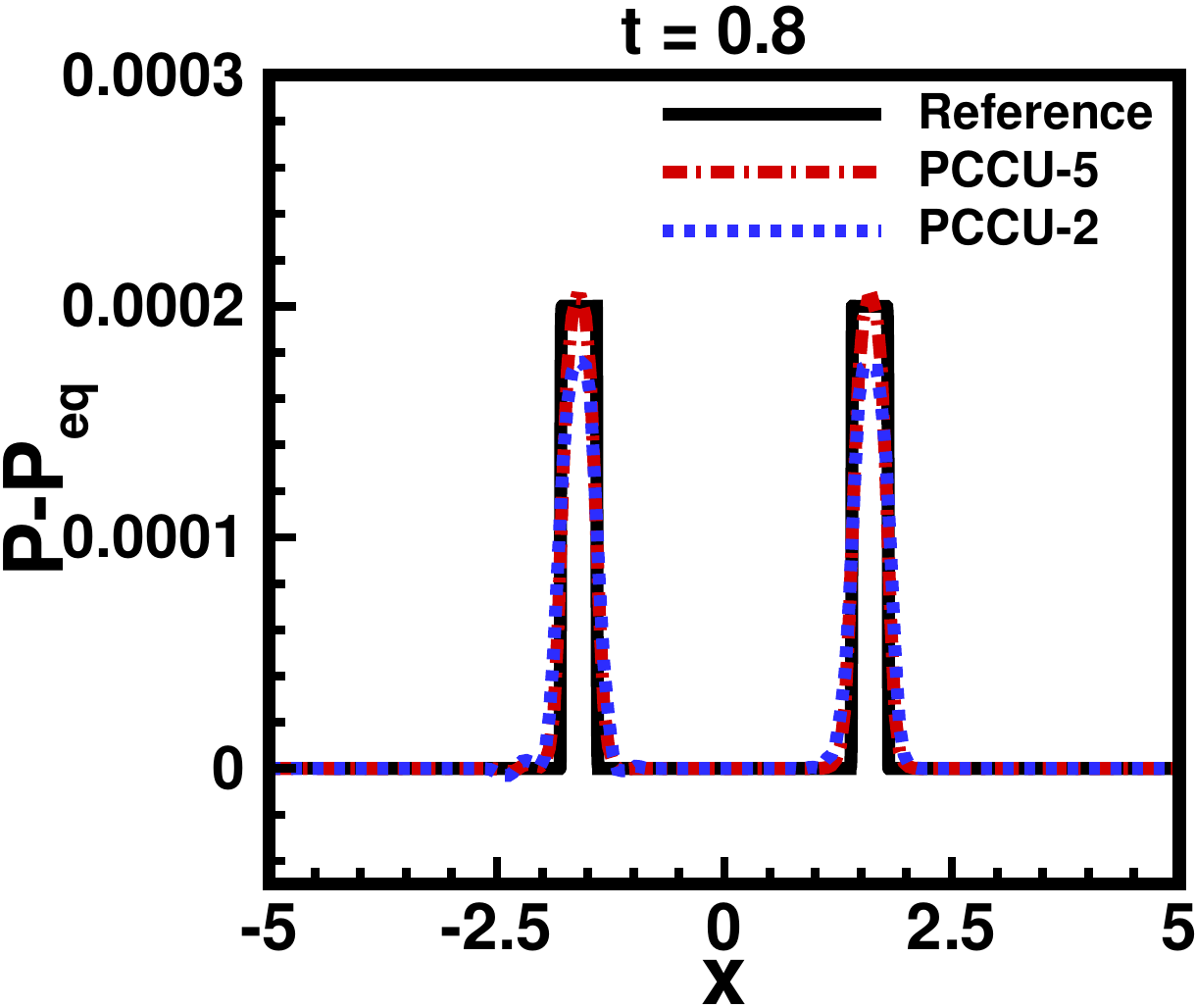}
\includegraphics[width=0.24\textwidth]{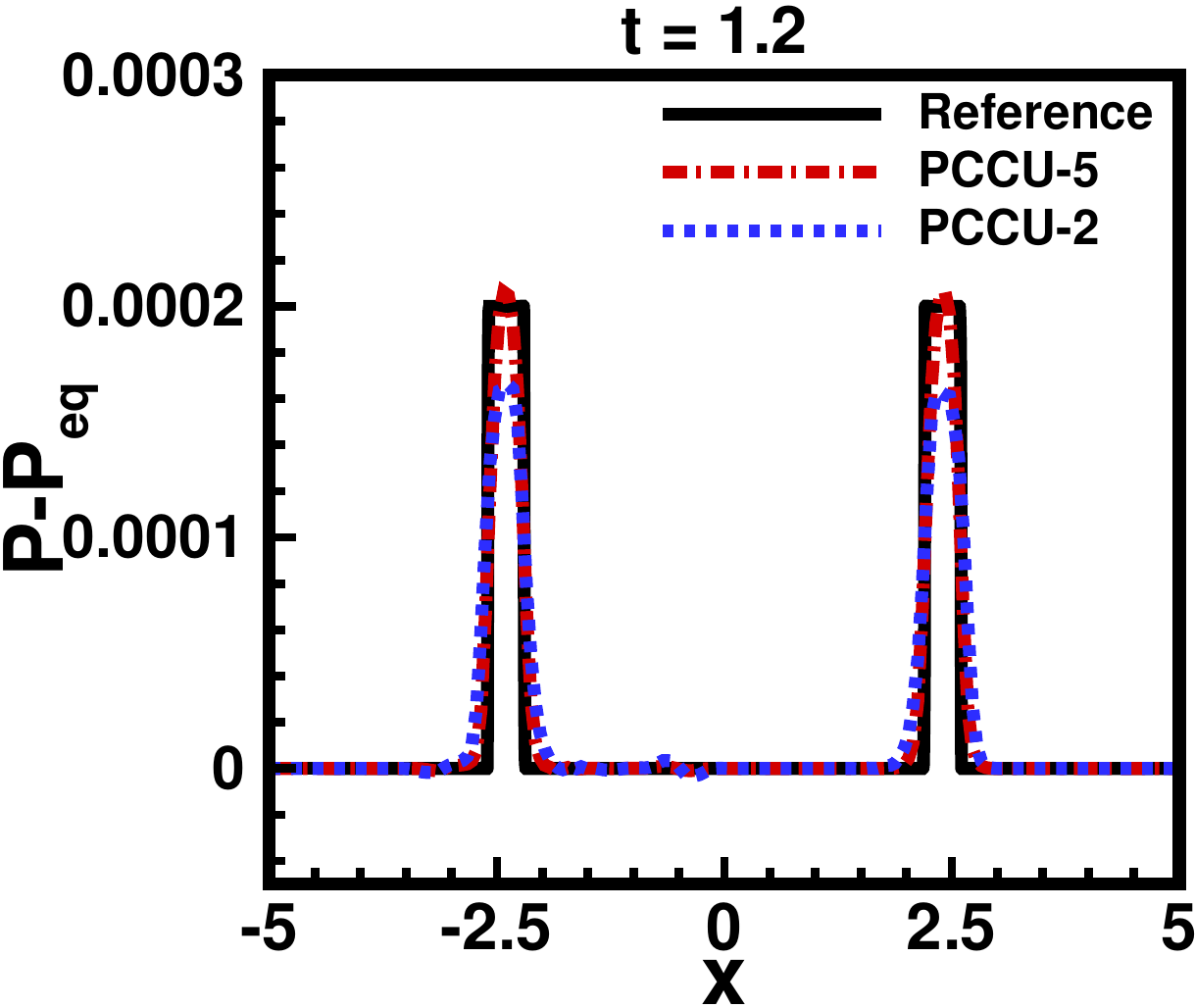}
\includegraphics[width=0.24\textwidth]{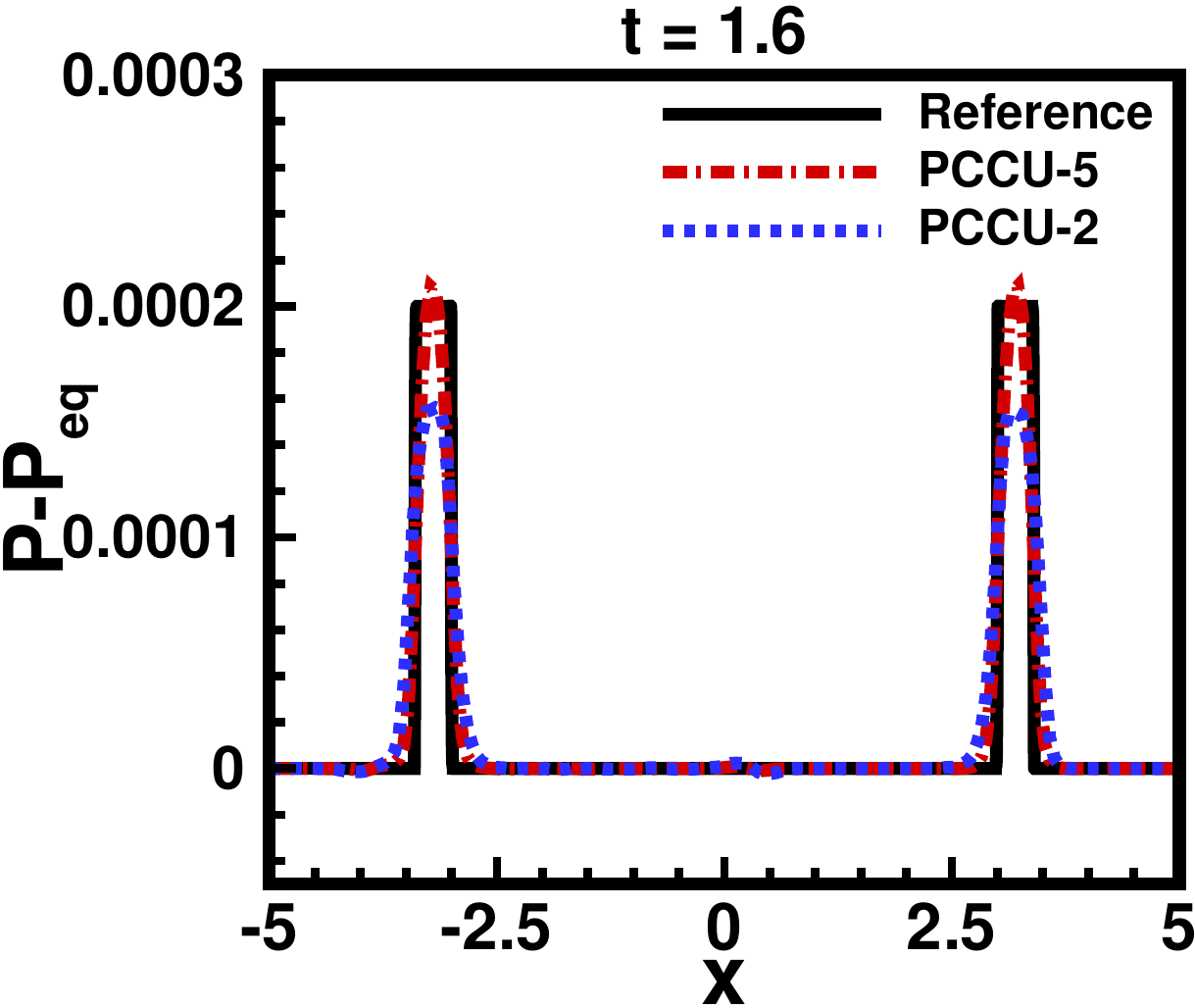}
\vskip-1mm
\caption{Example \ref{ex65}: $h-h_{\rm eq}$ (top row) and $P-P_{\rm eq}$ (bottom row) computed at times $t=0.4$, $0.8$, $1.2$, and $1.6$ by
the PCCU-5 and PCCU-2 schemes.}\label{fig67}
\end{figure}
\end{example}

\subsection{Two-Dimensional Examples}
\begin{example}[2-D numerical accuracy test]\label{ex66}
\rm In the first 2-D example taken from \cite{ZXX}, we test the experimental order of accuracy of the proposed 2-D PCCU-5 scheme. The
following smooth initial conditions,
\begin{equation*}
\begin{aligned}
&h(x,y,0)=10+e^{\sin(2\pi x)}\cos(2\pi y),&&q^x(x,y,0)=\sin(\cos(2\pi x))\sin(2\pi y),\\
&q^y(x,y,0)=\cos(2\pi x)\cos(\sin(2\pi y)),&&\theta(x,y,0)=9.812(2+\sin(2\pi x)\cos(2\pi y)),
\end{aligned}
\end{equation*}
are prescribed in the computational domain $[0,1]\times[0,1]$ subject to the periodic boundary conditions. The bottom topography is also
smooth: $Z(x,y)\sin(2\pi x)+\cos(2\pi y)$, and the time step is selected as $\dt=0.45\Delta^{5/3}$, where $\Delta=\dx=\dy$, to ensure that
the spatial and temporal contributions of the error are of the same order.

We compute the numerical solution until the final time $t=0.01$ on a sequence of uniform meshes with $N_x=N_y=25$, $50$, $100$, $200$,
$400$, and $800$. We then compute the $L^\infty$-errors and estimate the experimental convergence rates using the following Runge formulae,
which are based on the solutions computed on the three consecutive uniform grids with the mesh sizes $\Delta$, $2\Delta$, and $4\Delta$ and
denoted by $(\cdot)^\Delta$, $(\cdot)^{2\Delta}$, and $(\cdot)^{4\Delta}$, respectively:
$$
{\rm Error}(\Delta)\approx\frac{\delta_{12}^2}{|\delta_{12}-\delta_{24}|},\quad
{\rm Rate}(\Delta)\approx\log_2\left(\frac{\delta_{24}}{\delta_{12}}\right).
$$
Here, $\delta_{12}:=\|(\cdot)^\Delta-(\cdot)^{2\Delta}\|_{\infty}$ and $\delta_{24}:=\|(\cdot)^{2\Delta}-(\cdot)^{4\Delta}\|_{\infty}$.
The obtained results for $h$, $q^x$, $q^y$, and $h\theta$ are reported in Table \ref{tab64}, where one can clearly see that the fifth order
of accuracy is achieved.
\begin{table}[!htb]
\centering\small
\caption{$L^\infty$-errors and experimental orders of convergence for $h$, $q^x$, $q^y$, and $h\theta$.}\label{tab64}.
\begin{tabular}{ccccccccccccc}
\hline
&&\multicolumn{2}{c}{$h$} && \multicolumn{2}{c}{$q^x$}&&\multicolumn{2}{c}{$q^y$}&&\multicolumn{2}{c}{$h\theta$}\\
\cline{3-4}\cline{6-7}\cline{9-10}\cline{12-13}
$\Delta$&&Error   &Rate&&Error   &Rate&&Error   &Rate&&Error   &Rate\\\hline
$1/160$ &&3.37E-05&4.72&&4.40E-04&4.69&&5.39E-04&4.65&&4.32E-04&5.08\\
$1/320$ &&1.01E-06&4.89&&1.20E-06&4.93&&1.47E-05&4.92&&1.48E-05&4.98\\
$1/640$ &&3.02E-08&4.97&&3.61E-07&5.00&&4.38E-07&4.99&&5.36E-07&4.89\\
\hline
\end{tabular}
\end{table}
\end{example}

\begin{example}[2-D still-water equilibrium]\label{ex67}
\rm In this example, we demonstrate the ability of the proposed PCCU-5 scheme to exactly preserve 2-D still-water equilibria.

We consider the initial conditions,
\begin{equation*}
h(x,y,0)=3-Z(x,y),\quad u(x,y,0)=v(x,y,0)\equiv0,\quad\theta(x,y,0)\equiv\frac{39.248}{3},
\end{equation*}
which are prescribed in the computational domain $[-1,1]\times[-1,1]$ subject to the periodic boundary conditions. The nonsmooth bottom
topography in this example is
\begin{equation*}
Z(x,y)=\left\{\begin{aligned}
&0.5\,e^{-100\left((x+0.5)^2+(y+0.5)^2\right)},&&x<0,\\
&0.6\,e^{-100\left((x-0.5)^2+(y-0.5)^2\right)},&&\mbox{otherwise}.
\end{aligned}\right.
\end{equation*}

We compute the solution by the PCCU-5 scheme until the final time $t=1$ on a uniform mesh with $N_x\times N_y=20\times20$, $50\times50$ and
$100\times100$ cells and present the differences $\|h(\cdot,\cdot,1)-h(\cdot,\cdot,0)\|_\infty$,
$\|q^x(\cdot,\cdot,1)-q^x(\cdot,\cdot,0)\|_\infty$, $\|q^y(\cdot,\cdot,1)-q^y(\cdot,\cdot,0)\|_\infty$, and
$\|\theta(\cdot,\cdot,1)-\theta(\cdot,\cdot,0)\|_\infty$ in Table \ref{tab65}. As one can see, all entries in the table are close to the machine errors; hence, the proposed PCCU-5 scheme can preserve the still-water equilibrium.
\begin{table}[!htb]
\centering\small
\caption{Example \ref{ex67}: Errors ($\|h(\cdot,\cdot,1)-h(\cdot,\cdot,0)\|_\infty$, $\|q^x(\cdot,\cdot,1)-q^x(\cdot,\cdot,0)\|_\infty$,
$\|q^y(\cdot,\cdot,1)-q^y(\cdot,\cdot,0)\|_\infty$, and $\|\theta(\cdot,\cdot,1)-\theta(\cdot,\cdot,0)\|_\infty$) in $h$, $q^x$, $q^y$, and
$\theta$.}\label{tab65}
\begin{tabular}{cccccccccccccc}
\hline
$N_x\times N_y$&&Error in $h$&&Error in $q^x$&&Error in $q^y$&&Error in $\theta$\\\hline
$20\times20$   &&2.53E-14&&5.69E-14&&7.06E-14&&1.63E-13\\
$50\times50$   &&4.22E-14&&1.39E-13&&8.83E-14&&3.30E-13\\
$100\times100$ &&6.71E-14&&1.78E-13&&2.12E-13&&5.26E-13\\
\hline
\end{tabular}
\end{table}
\end{example}

\begin{example}[Small perturbation of the 2-D still-water equilibrium]\label{ex68}
\rm In this example taken from \cite{ZXX}, we test the ability of the PCCU-5 scheme to capture the propagation of a small perturbation of a
still-water steady state.

The initial data,
\begin{equation*}
(h,u,v,\theta)\Big|_{(x,y,0)}=\left\{\begin{aligned}
&\left(6-Z(x,y)+\eta,0,0,\frac{235.488}{6+\eta}\right),&&0.05\le x\le0.15,\\
&\hspace*{0.15cm}\left(6-Z(x,y),0,0,39.248\right),&&\quad~\mbox{otherwise},
\end{aligned}\right.
\end{equation*}
are prescribed in the computational domain $[-2,2]\times[0,1]$ subject to the free boundary conditions. The bottom topography contains an
isolated elliptical-shaped hump:
\begin{equation*}
Z(x,y)=3\,e^{-5(x-0.9)^2-50(y-0.5)^2}.
\end{equation*}

We compute the solution by the PCCU-5 scheme until the final time $t=0.08$ on a uniform mesh with $200\times100$ cells. In Figure
\ref{fig68}, we present time snapshots at $t=0.02$, $t=0.05$, and $t=0.08$ of the obtained $h+Z-6$ (deviation of the water surface from the
equilibrium level), $q^x$, $q^y$, and $h\theta$. As one can see, the right-going disturbance propagates past the hump and develops
small-scale flow structures, which are resolved by the PCCU-5 scheme in a sharp and non-oscillatory manner. The results are consistent with
those reported in \cite{ZXX}.
\begin{figure}[!htb]
\centering
\mbox{
\includegraphics[height=0.105\textheight]{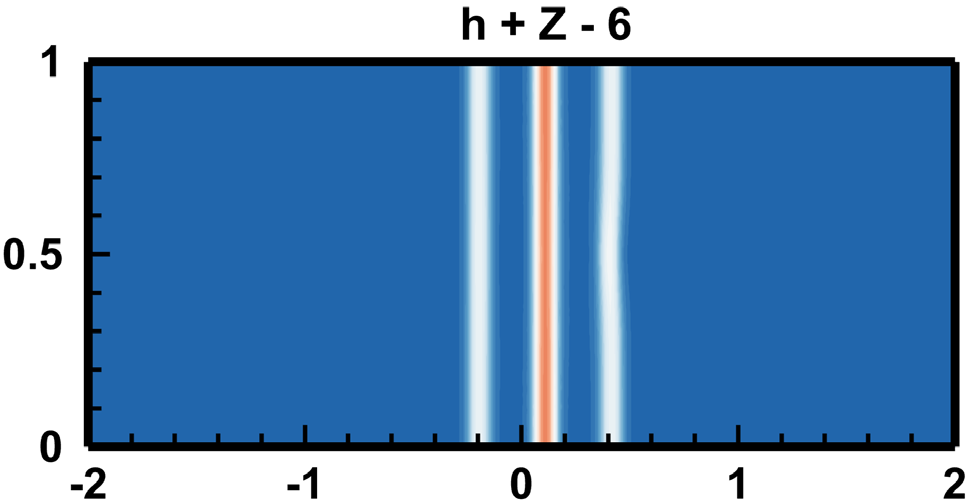}
\includegraphics[height=0.105\textheight]{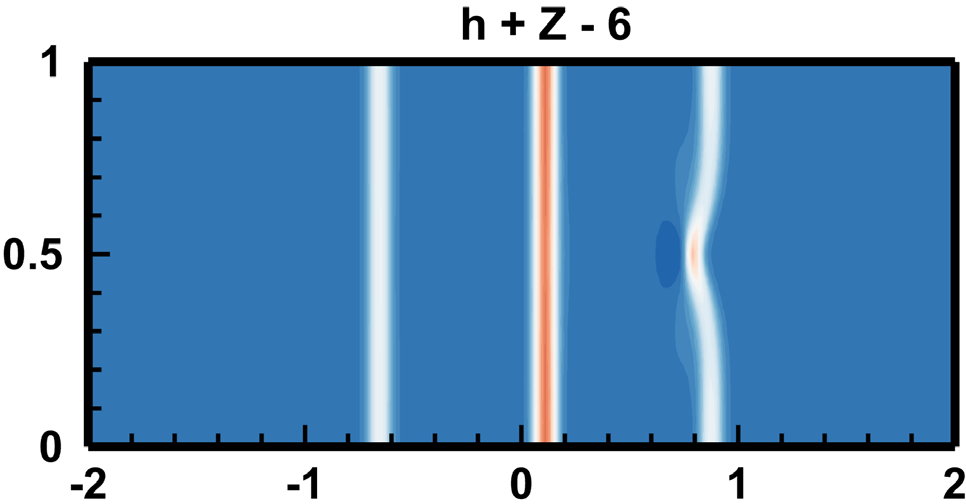}
\includegraphics[height=0.105\textheight]{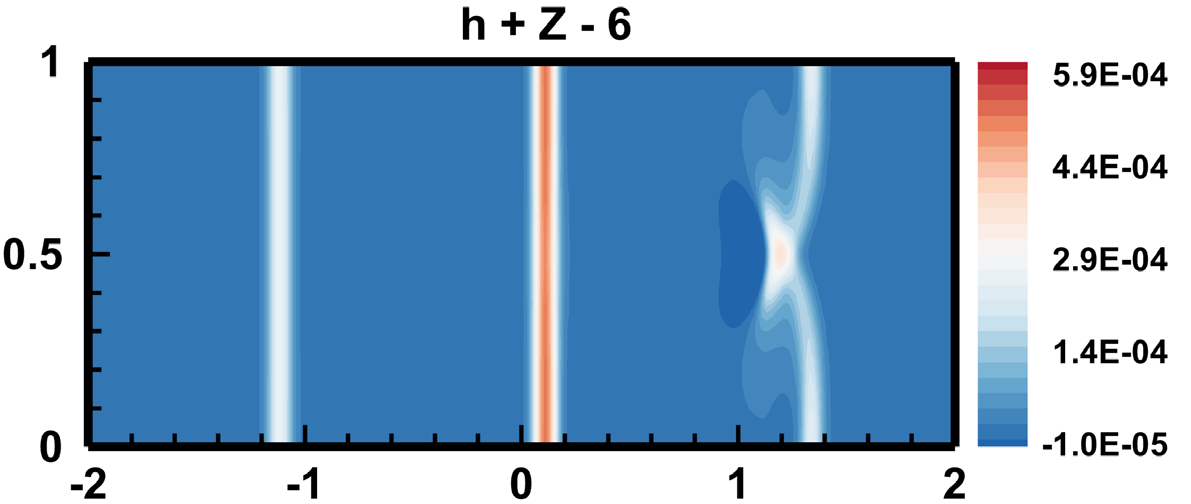}}
\vskip2pt
\mbox{
\includegraphics[height=0.105\textheight]{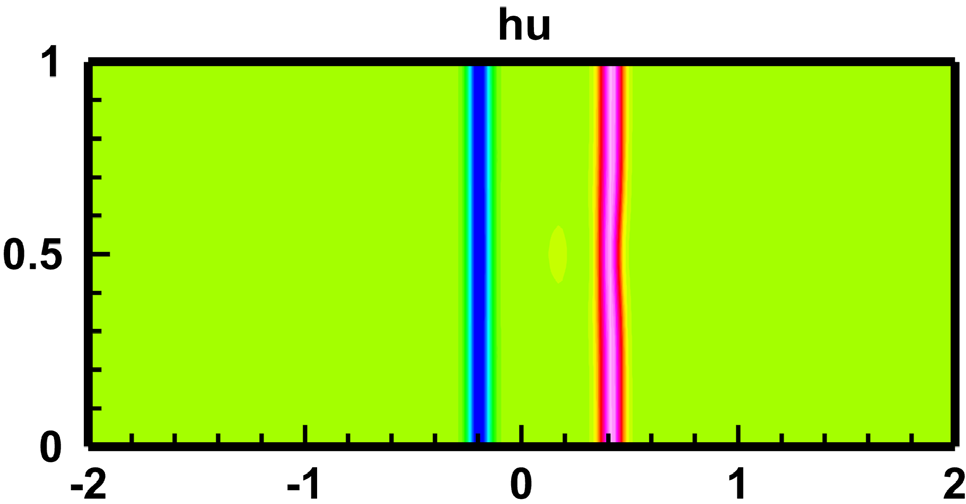}
\includegraphics[height=0.105\textheight]{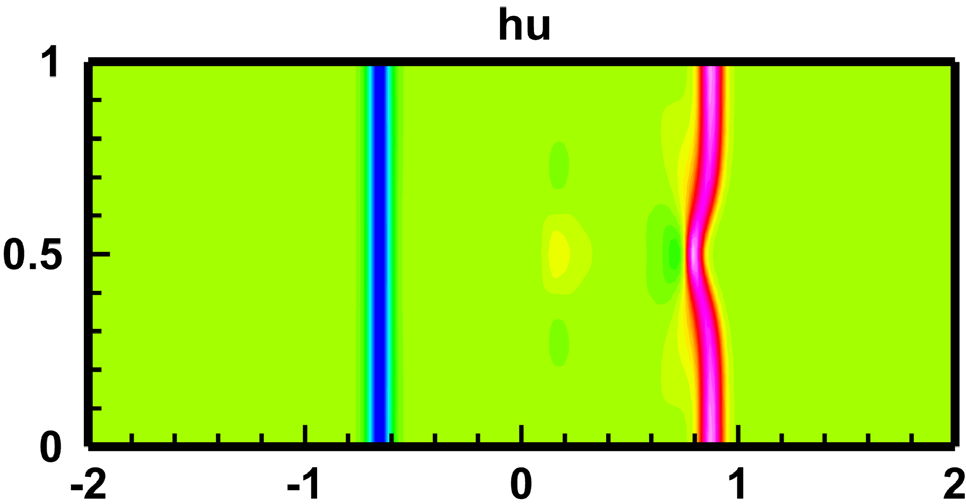}
\includegraphics[height=0.105\textheight]{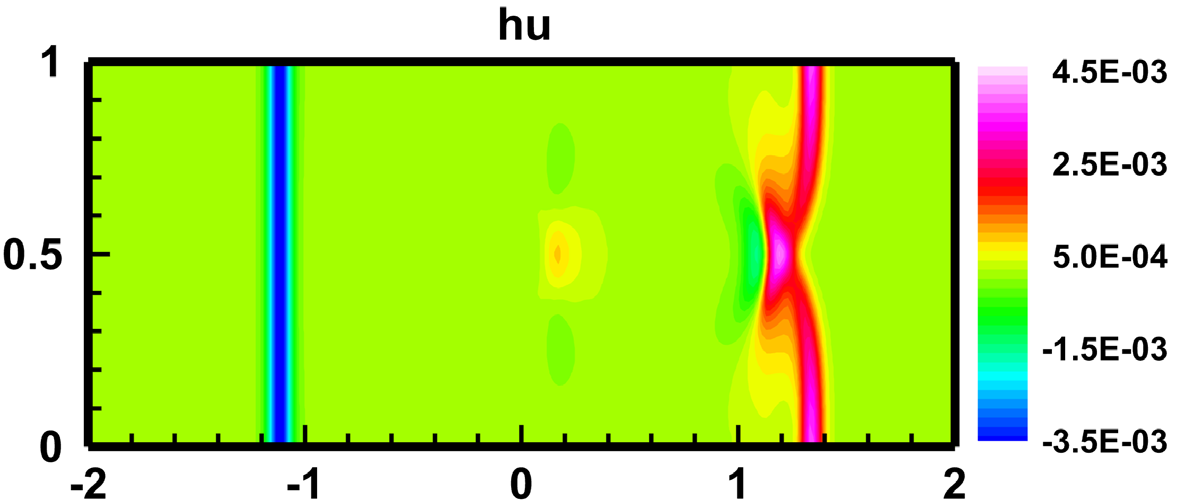}}
\vskip2pt
\mbox{
\includegraphics[height=0.105\textheight]{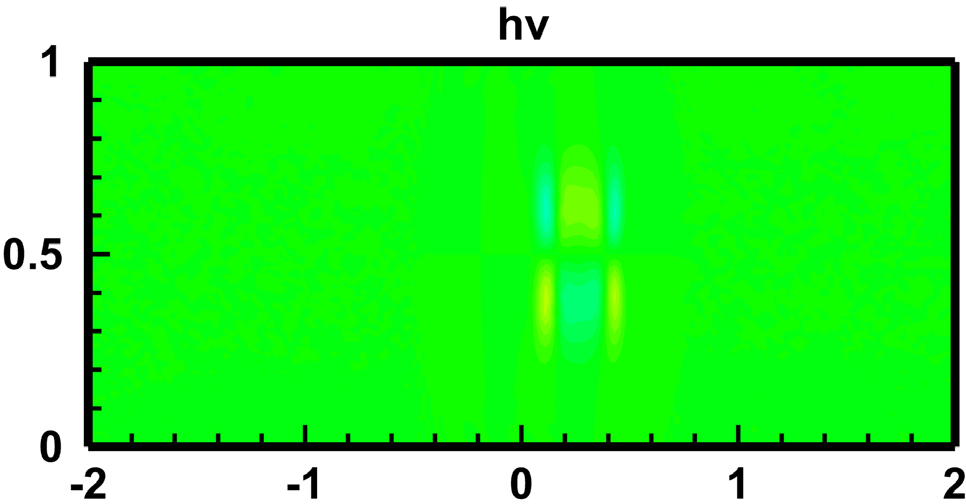}
\includegraphics[height=0.105\textheight]{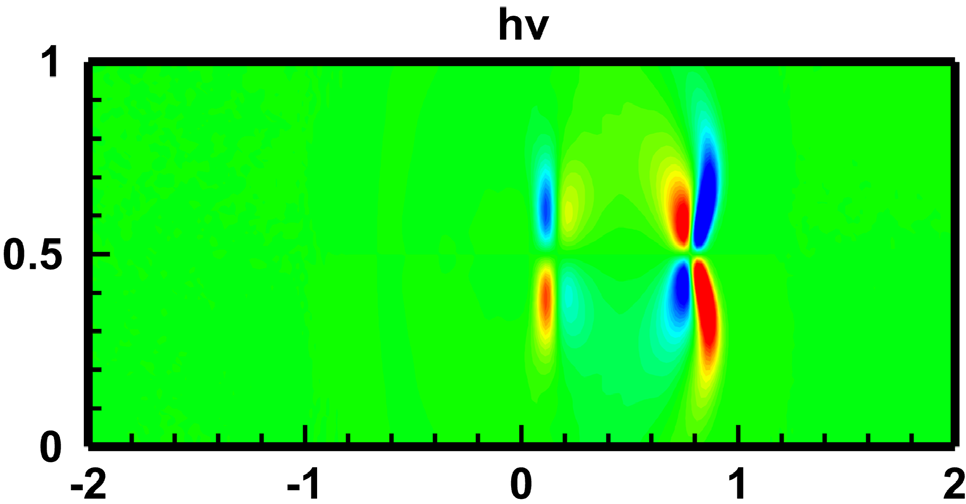}
\includegraphics[height=0.105\textheight]{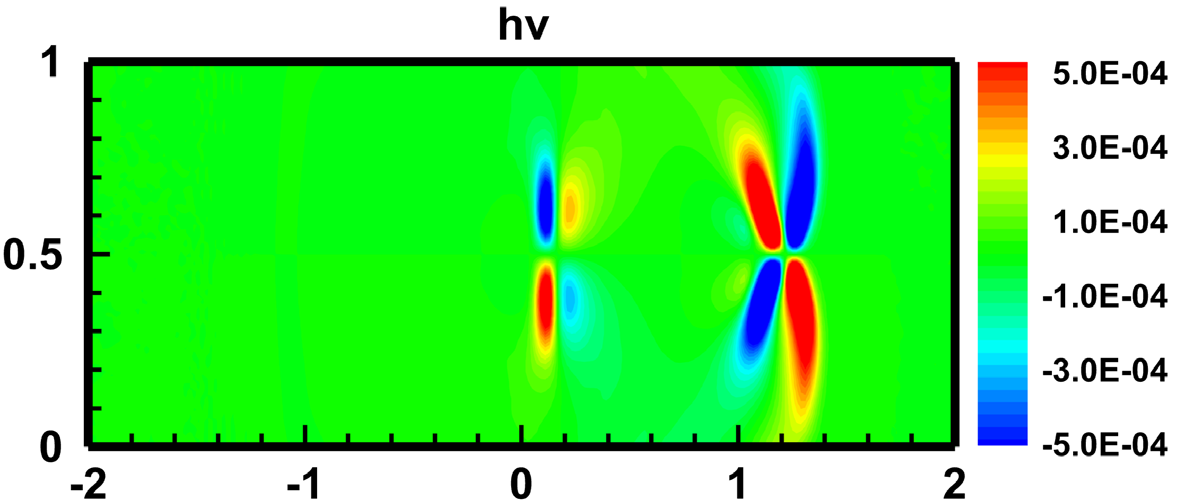}}
\vskip2pt
\mbox{
\includegraphics[height=0.105\textheight]{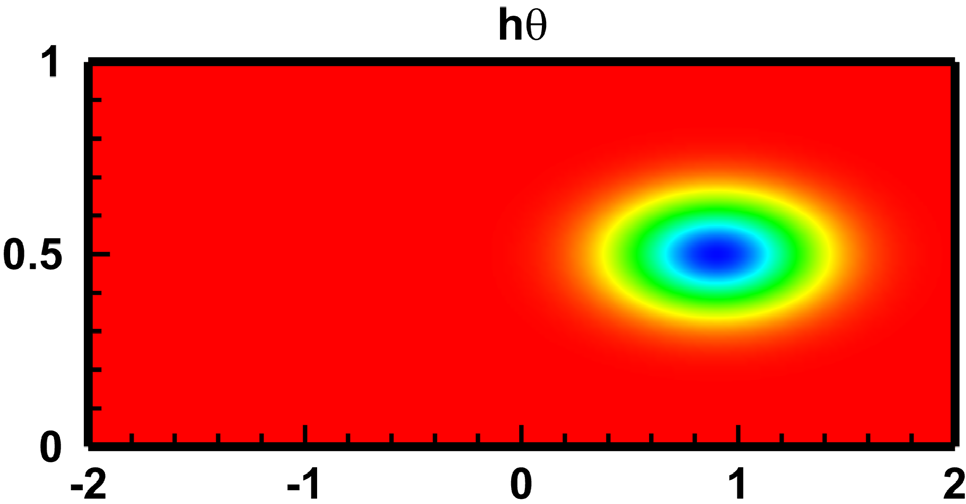}
\includegraphics[height=0.105\textheight]{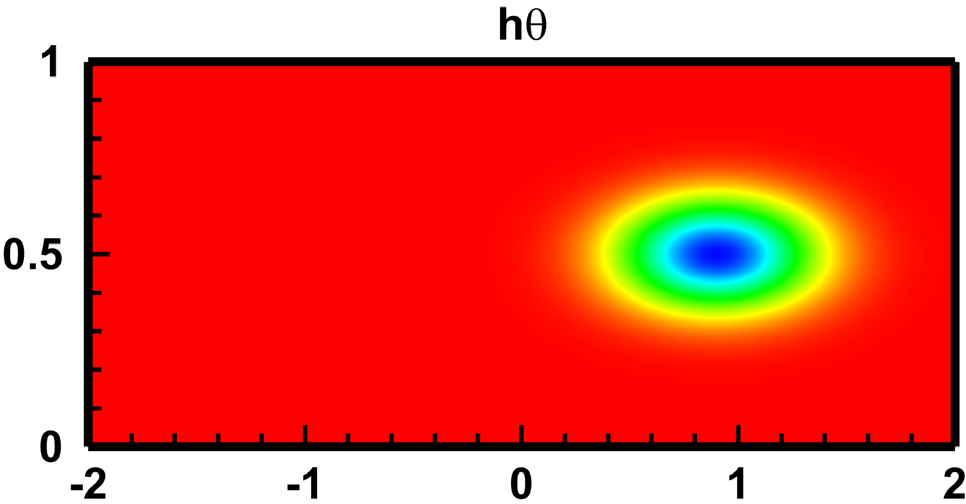}
\includegraphics[height=0.105\textheight]{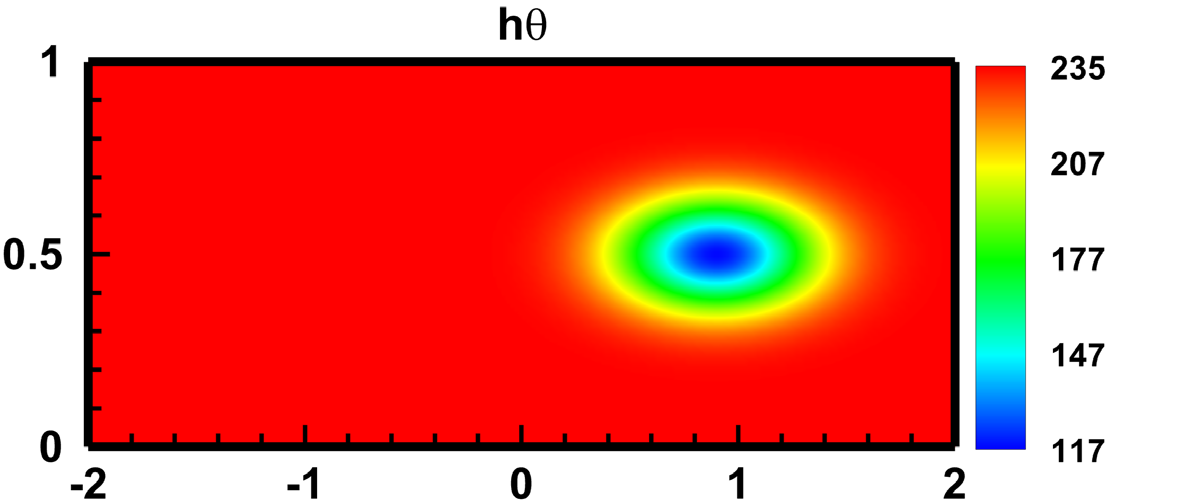}}
\vskip-1mm
\caption{Example \ref{ex68}: $h+Z-6$, $q^x$, $q^y$, and $h\theta$ computed by the PCCU-5 scheme at times $t=0.02$ (left column), $0.05$
(middle column), and $0.08$ (right column).}
\label{fig68}
\end{figure}
\end{example}

\begin{example}[Radial dam-break over a flat bottom topography]\label{ex69}
\rm In the next example taken from \cite{Liu22}, we consider a dam-break problem over a flat bottom topography $Z(x,y)\equiv0$. The initial
conditions,
\begin{equation*}
(h,u,v,\theta)\Big|_{(x,y,0)}=\left\{\begin{aligned}
&(2,0,0,9.812),&& x^2+y^2\le0.25, \\
&(1,0,0,14.718),&&\quad\mbox{otherwise},
\end{aligned}\right.
\end{equation*}
are prescribed in the computational domain $[-1,1]\times[-1,1]$ subject to the periodic boundary conditions.

We compute the solution by the PCCU-5 scheme until the final time $t=0.05$ on a uniform mesh with $200\times200$ cells. As shown in Figure
\ref{fig69}, the obtained solution is both sharp and oscillation-free, and it agrees well with the results reported in \cite{Liu22}.
\begin{figure}[!htb]
\centering
\includegraphics[width=0.32\textwidth]{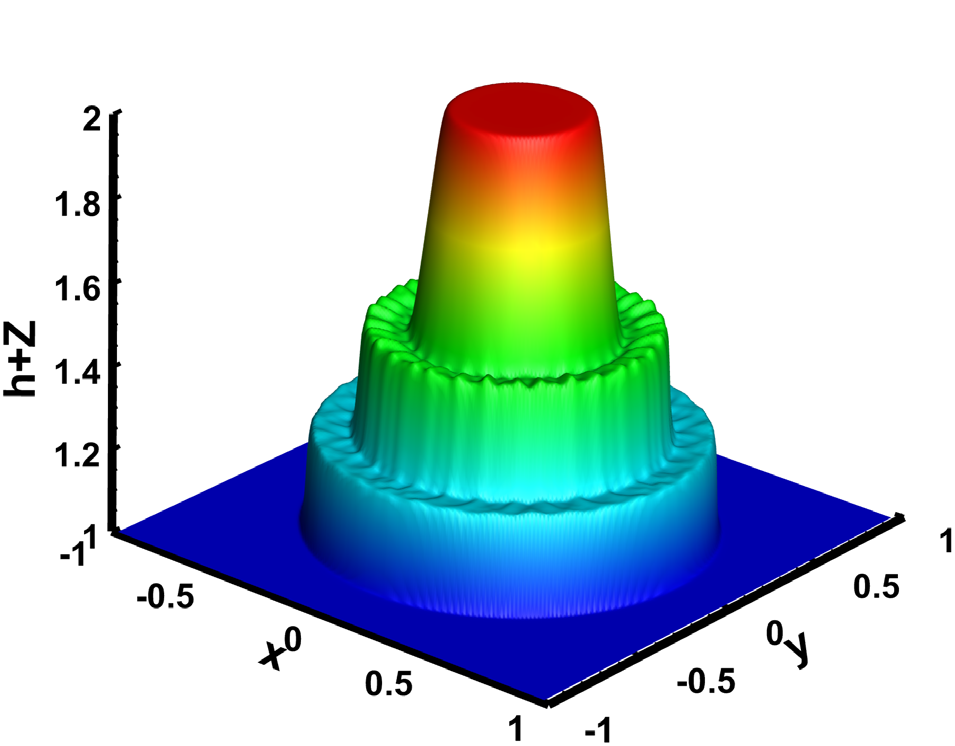}
\includegraphics[width=0.32\textwidth]{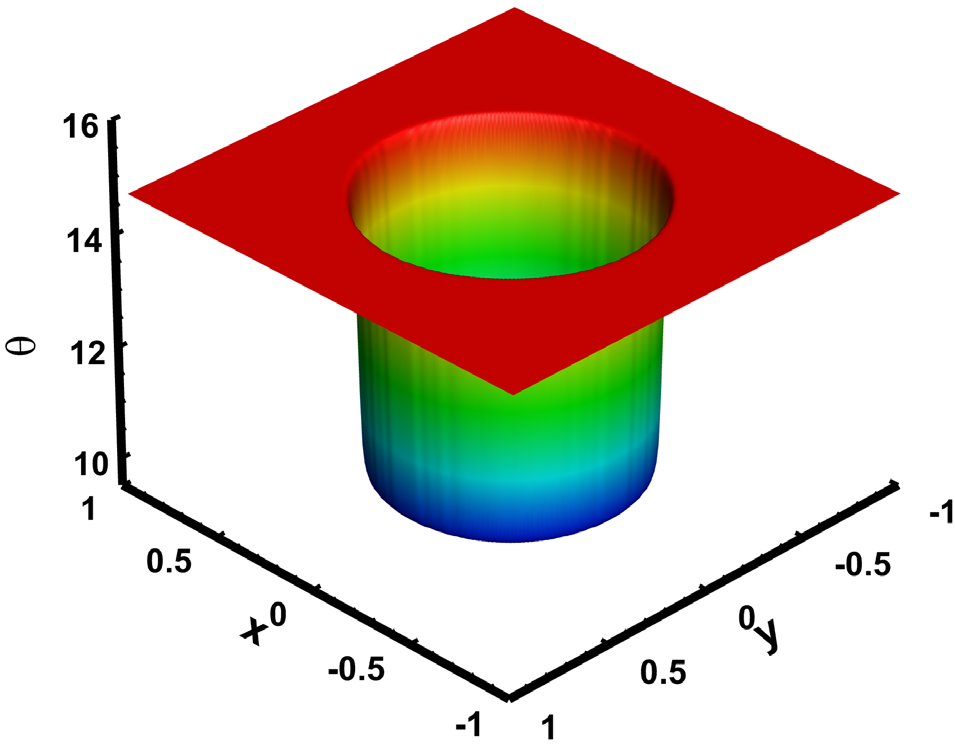}
\includegraphics[width=0.32\textwidth]{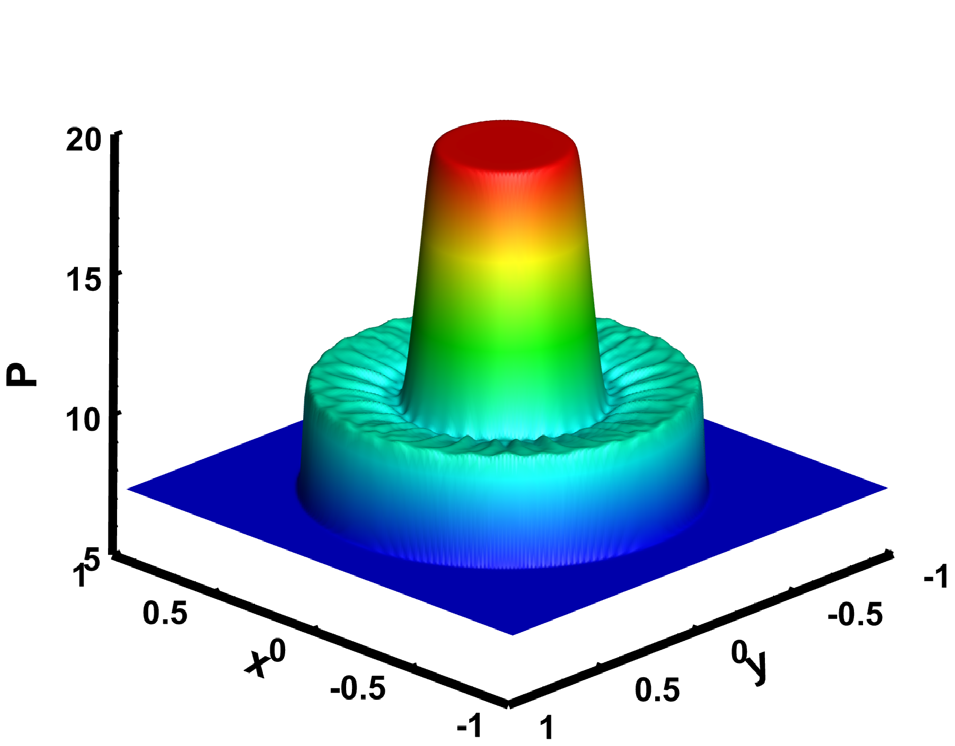}
\vskip-1mm
\caption{Example \ref{ex69}: $h+Z$, $\theta$, and $P$ computed by the PCCU-5 scheme.}
\label{fig69}
\end{figure}
\end{example}

\begin{example}[2-D steady-state solution]\label{ex610}
\rm In this example also taken from \cite{ZXX}, the computational domain is $[-1,1]\times[-1,1]$, the bottom topography consists of two
Gaussian humps:
$$
Z(x,y)=\left\{\begin{aligned}
0.5\,e^{-100\left((x+0.5)^2+(y+0.5)^2\right)},&&x<0,\\
0.6\,e^{-100\left((x-0.5)^2+(y-0.5)^2\right)},&&x>0,
\end{aligned}\right.
$$
the initial conditions are
\begin{equation}
(h,u,v,\theta)\Big|_{(x,y,0)}=\left\{\begin{aligned}
&\left(3-Z(x,y),0,0,\frac{4}{3}\right),&&x^2+y^2\le0.25,\\
&~\left(2-Z(x,y),0,0,3\right),&&\quad\mbox{otherwise},
\end{aligned}\right.
\label{6.7}
\end{equation}
and the periodic boundary conditions are imposed. This setting corresponds to a still-water equilibrium containing a temperature jump.

We compute the solution by the PCCU-5 scheme until the final time $t=0.12$ on a uniform mesh with $100\times100$ cells and present the
obtained solutions ($h+Z$, $\theta$, and $P$) in Figure \ref{fig610}. As one can see, the steady state is preserved by the proposed WB
scheme.
\begin{figure}[!htb]
\centering
\includegraphics[width=0.32\textwidth]{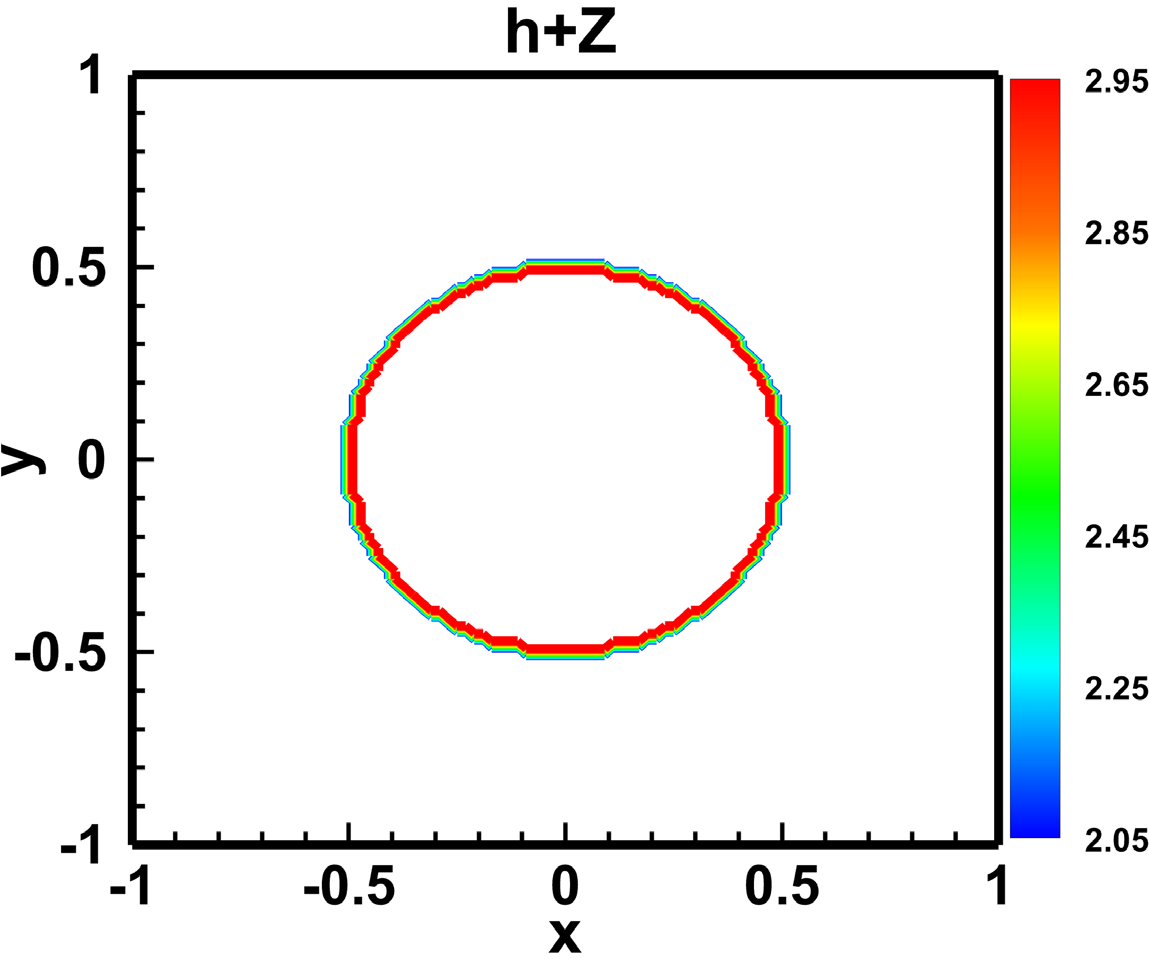}
\includegraphics[width=0.32\textwidth]{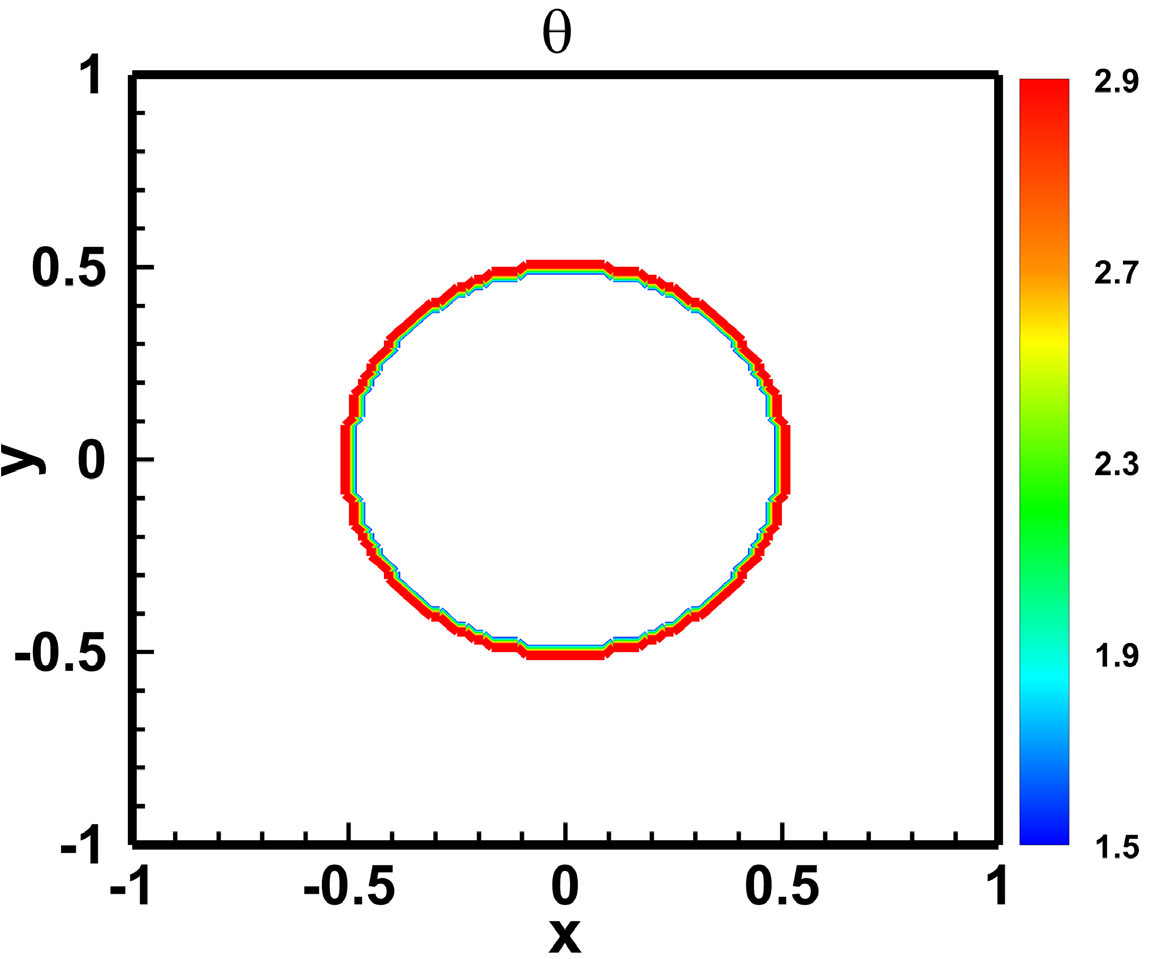}
\includegraphics[width=0.32\textwidth]{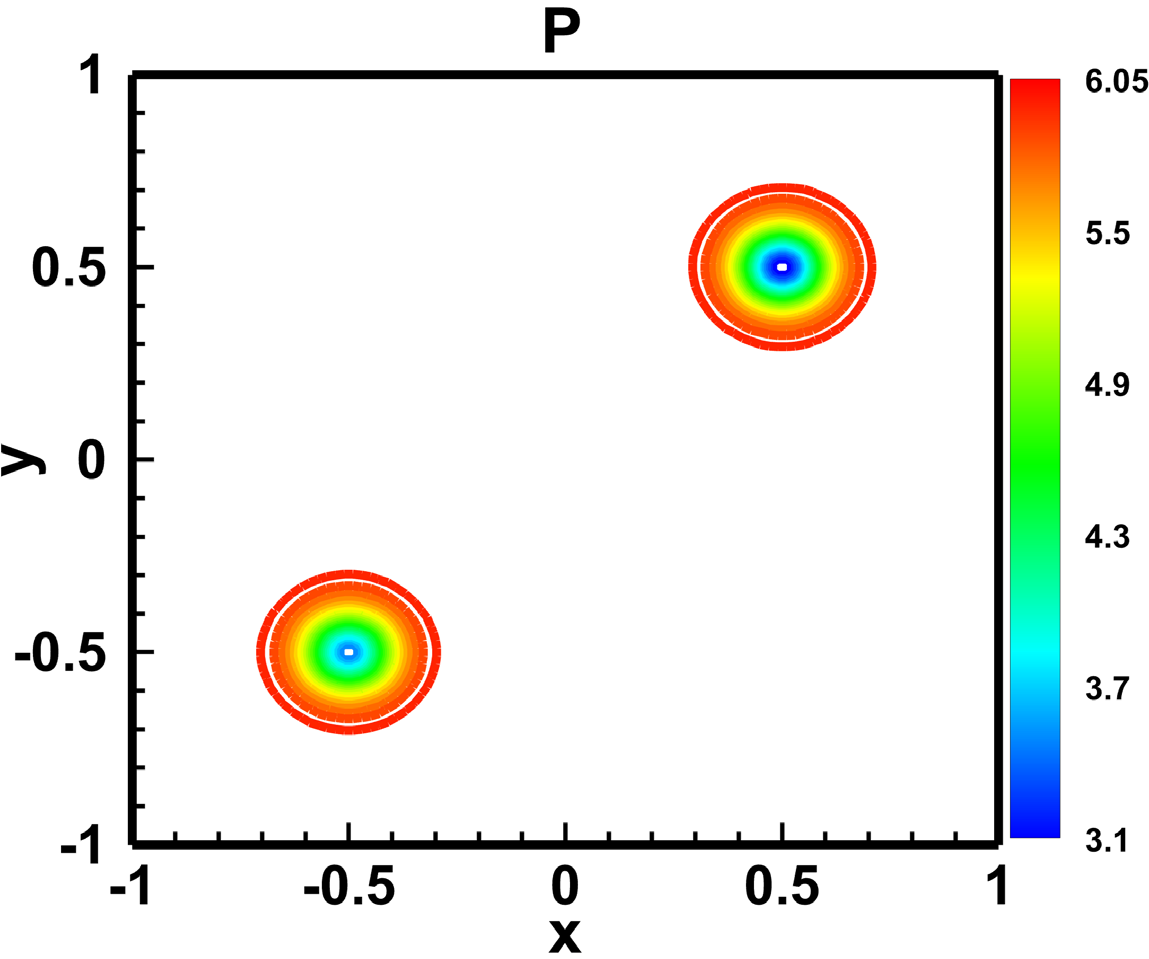}
\vskip-1mm
\caption{Example \ref{ex610}: Contours plots of $h+Z$, $\theta$, and $P$: $23$ uniformly spaced contours in the ranges $[2.05,2.95]$ (for
$h+Z$), $[1.5,2.9]$ (for $\theta$), and $[3.1,6.05]$ (for $P$).}
\label{fig610}
\end{figure}
\end{example}

\begin{example}[Circular perturbation of 2-D steady-state solution]\label{ex611}
\rm In this example, we follow \cite{ZXX} and add a perturbation of size $0.1$ inside the annulus $0.01\le x^2+y^2\le0.09$ to the initial
$h$ given in \eref{6.7}. The rest of the problem settings are the same as in Example \ref{ex610}. The solution, computed by the PCCU-5
scheme at time $t=0.1$, is presented in Figure \ref{fig611}, where one can observe that the perturbation is accurately captured without any
spurious oscillations.
\begin{figure}[!htb]
\centering
\includegraphics[width=0.32\textwidth]{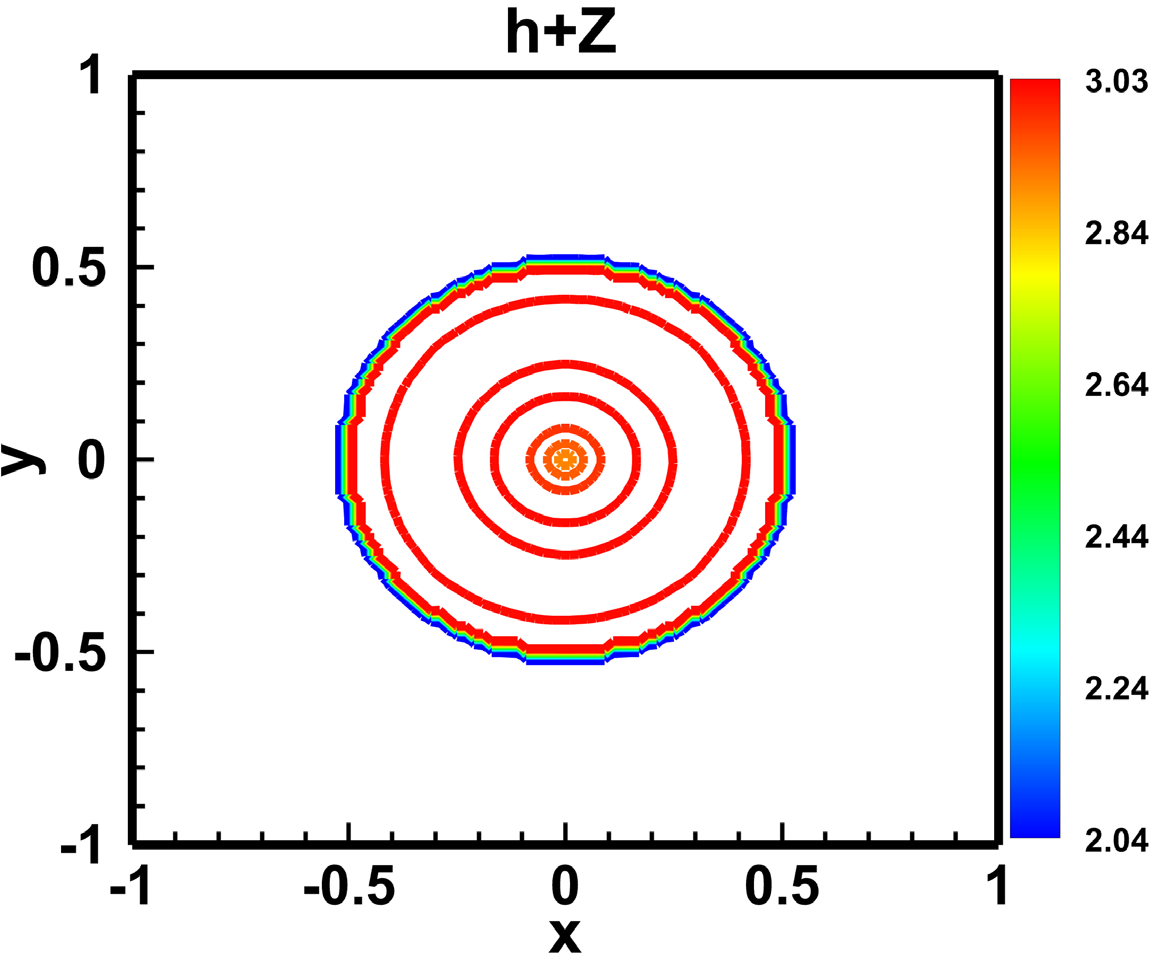}
\includegraphics[width=0.32\textwidth]{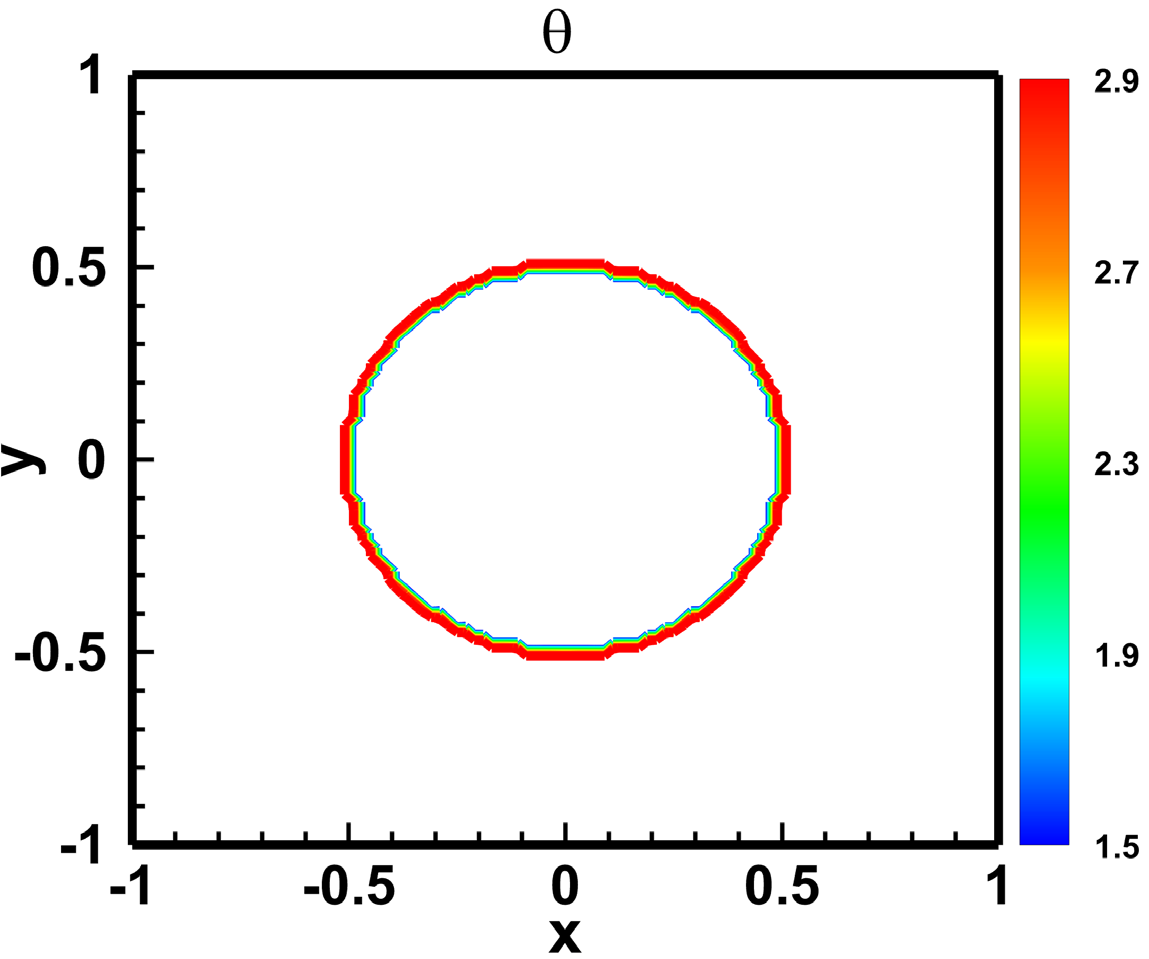}
\includegraphics[width=0.32\textwidth]{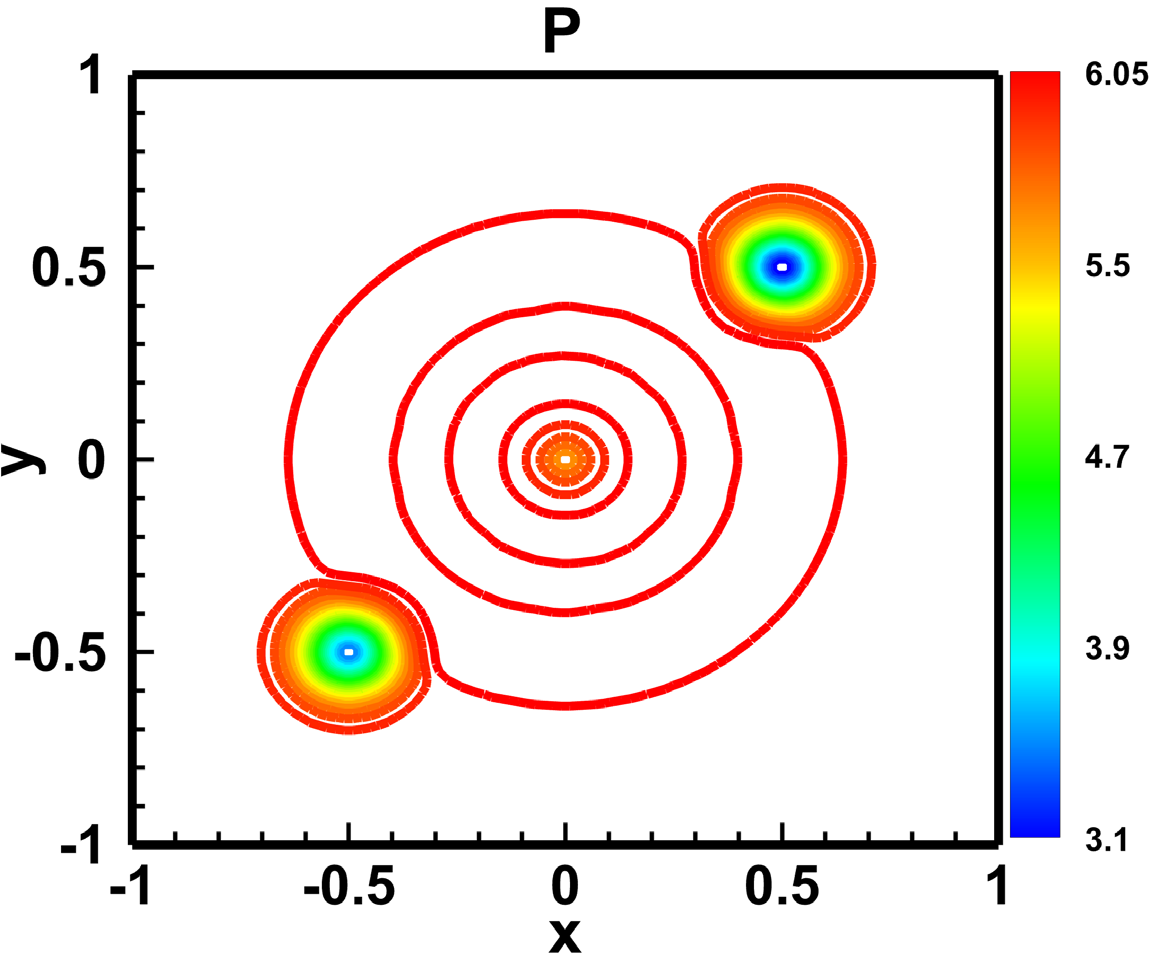}
\vskip-1mm
\caption{Example \ref{ex611}: Contours plots of $h+Z$, $\theta$, and $P$: $23$ uniformly spaced contours in the ranges $[2.05,2.95]$ (for
$h+Z$), $[1.5,2.9]$ (for $\theta$), and $[3.1,6.05]$ (for $P$).}
\label{fig611}
\end{figure}
\end{example}

\begin{example}[Quasi 1-D moving-water equilibria their small perturbations]\label{ex612}
\rm In this example, we test the ability of the proposed PCCU-5 scheme to exactly preserve quasi 1-D moving-water equilibria and accurately
capture their small generically 2-D perturbations.

The computational domain in this example is $[0,25]\times[0,10]$, and we take a continuous bottom topography, which varies in the
$x$-direction only: $Z(x,y)=Z(x)$ with the 1-D function $Z(x)$ given by \eref{6.1}. We consider the following two sets of initial and
boundary conditions, which correspond to the subcritical and supercritical equilibria:

\medskip
\noindent
$\bullet$ {\bf Subcritical equilibrium}:
\begin{equation}
\begin{aligned}
&(\Ep^x,q^x,q^y,\theta)\Big|_{(x,y,0)}=
\big(\Ep^x_{\rm eq},q^x_{\rm eq},q^y_{\rm eq},\theta_{\rm eq}\big)\Big|_{(x,y)}\equiv(110.33025,4.42\sqrt{5},0,49.06),\\
&h(25,y,t)=2,~~q^x(0,y,t)=4.42\sqrt{5};
\end{aligned}
\label{6.8}
\end{equation}

\medskip
\noindent
$\bullet$ {\bf Supercritical equilibrium}:
\begin{equation}
\begin{aligned}
&(\Ep^x,q^x,q^y,\theta)\Big|_{(x,y,0)}=
\big(\Ep^x_{\rm eq},q^x_{\rm eq},q^y_{\rm eq},\theta_{\rm eq}\big)\Big|_{(x,y)}\equiv(458.12,24\sqrt{5},49.06),\\
&h(0,y,t)=2,~~q^x(0,t)=24\sqrt{5}.
\end{aligned}
\label{6.9}
\end{equation}

In both cases, the rest of the boundary conditions in the $x$-direction are free (zero-order extrapolation), and the solid wall boundary
conditions are imposed in the $y$-direction.

As in the 1-D case (Example \ref{ex62}), the initial data in \eref{6.8}--\eref{6.9} are prescribed in terms of the equilibrium variables,
while $h_{\rm eq}$ is obtained by solving the cubic equations
\begin{equation*}
\hf\left(\frac{(q^x_{\rm eq})_{j,k}}{(h_{\rm eq})_{j,k}}\right)^2+{(\theta_{\rm eq})_{j,k}}\big[(h_{\rm eq})_{j,k}+Z_{j,k}\big]-
(\Ep^x_{\rm eq})_{j,k}=0
\end{equation*}
for $(h_{\rm eq})_{j,k}$ so that we set the initial data $h_{j,k}(0)=(h_{\rm eq})_{j,k}$.

We compute the solution by the PCCU-5 scheme until the final time $t=20$ on a uniform mesh with $100\times40$ cells and present the
differences $\|h(\cdot,\cdot,20)-h_{\rm eq}\|_\infty$, $\|q^x(\cdot,\cdot,20)-q^x_{\rm eq}\|_\infty$,
$\|q^y(\cdot,\cdot,20)-q^y_{\rm eq}\|_\infty$, $\|(h\theta)(\cdot,\cdot,20)-h_{\rm eq}\theta_{\rm eq}\|_\infty$, and
$\|\Ep^x(\cdot,\cdot,20)-\Ep^x_{\rm eq}\|_\infty$ in Table \ref{tab67}. As one can see, all of the entries in the table are close to the
machine errors, and hence the proposed PCCU-5 scheme can preserve moving-water equilibria in both of the studied cases.
\begin{table}[!htb]
\centering\small
\caption{Example \ref{ex612}: Errors ($\|h(\cdot,\cdot,20)-h_{\rm eq}\|_\infty$, $\|q^x(\cdot,\cdot,20)-q^x_{\rm eq}\|_\infty$,
$\|q^y(\cdot,\cdot,20)-q^y_{\rm eq}\|_\infty$, $\|(h\theta)(\cdot,\cdot,20)-h_{\rm eq}\theta_{\rm eq}\|_\infty$, and
$\|\Ep^x(\cdot,\cdot,20)-\Ep^x_{\rm eq}\|_\infty$) in $h$, $q^x$, $q^y$, $h\theta$, and $\Ep^x$.}\label{tab67}
\begin{tabular}{lcccccccccccccccc}
\hline
Case         &&Error in $h$&&Error in $q^x$&&Error in $q^y$&&Error in $h\theta$&&Error in $\Ep^x$\\\hline
Subcritical  &&5.77E-15&&4.08E-14&&0.00E+00&&3.84E-13&&2.56E-13\\
Supercritical&&1.06E-14&&1.92E-13&&0.00E+00&&1.23E-12&&1.48E-12\\
\hline
\end{tabular}
\end{table}

We then add small genuinely 2-D circular perturbations to the water depth in the above equilibria without changing the initial $q^x$, $q^y$,
or $h\theta$. The initial $h$ are:

\medskip
\noindent
$\bullet$ {\bf Subcritical case}:
\begin{equation*}
h(x,y,0)=h_{\rm eq}(x,y)+\left\{\begin{aligned}
&0.01,&&(x-6)^2+(y-5)^2<0.25,\\
&0,&&\mbox{otherwise};
\end{aligned}\right.
\end{equation*}

\medskip
\noindent
$\bullet$ {\bf Supercritical case}:
\begin{equation*}
h(x,y,0)=h_{\rm eq}(x,y)+\left\{\begin{aligned}
&0.05,&&(x-3)^2+(y-5)^2<0.25,\\
&0,&&\mbox{otherwise}.
\end{aligned}\right.
\end{equation*}

We compute the solution by the PCCU-5 scheme until the final time $t=0.4$ on three different uniform meshes with $100\times40$,
$200\times80$, and $400\times160$ cells. We plot the obtained differences $h(x,y,0.4)-h_{\rm eq}(x,y)$ in Figures \ref{fig612} and
\ref{fig613} for the subcritical and supercritical cases, respectively. As one can see, the small perturbation of the quasi 1-D moving
steady states are captured in a non-oscillatory manner even on the coarsest mesh, and the mesh convergence can be observed in both cases.
\begin{figure}[!htb]
\centering
\includegraphics[width=0.31\textwidth]{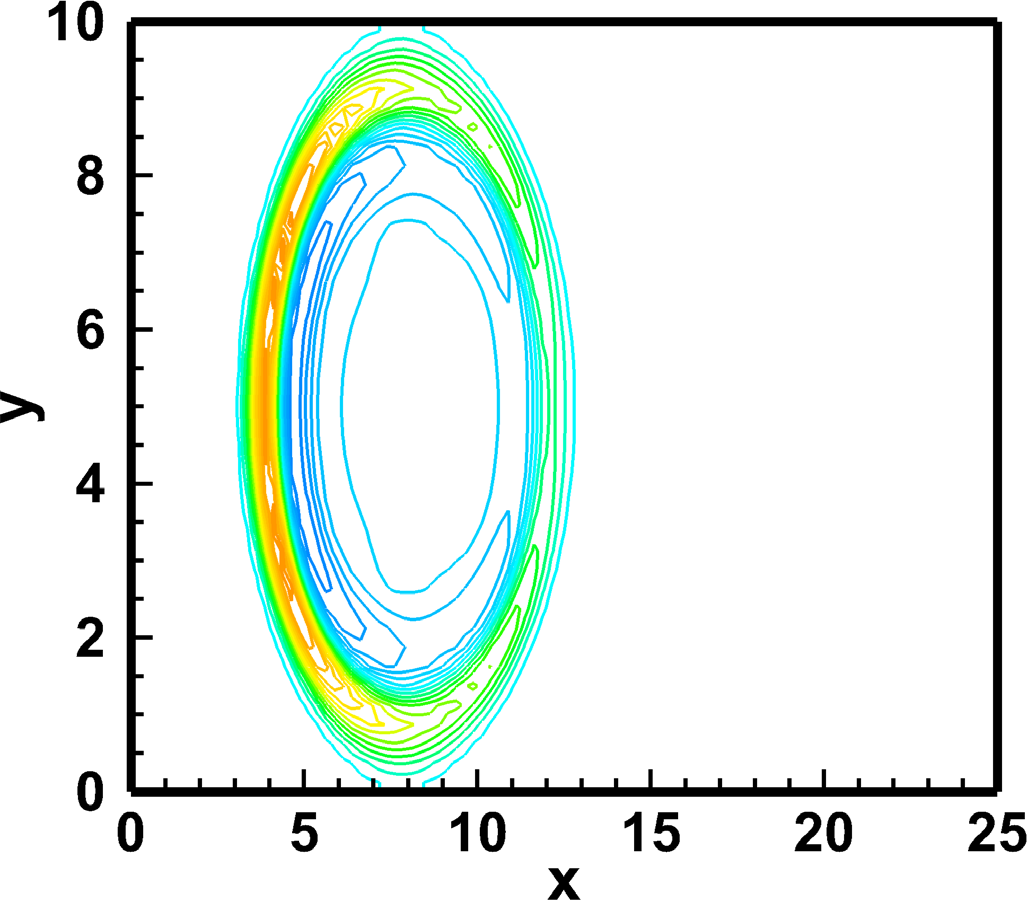}
\includegraphics[width=0.31\textwidth]{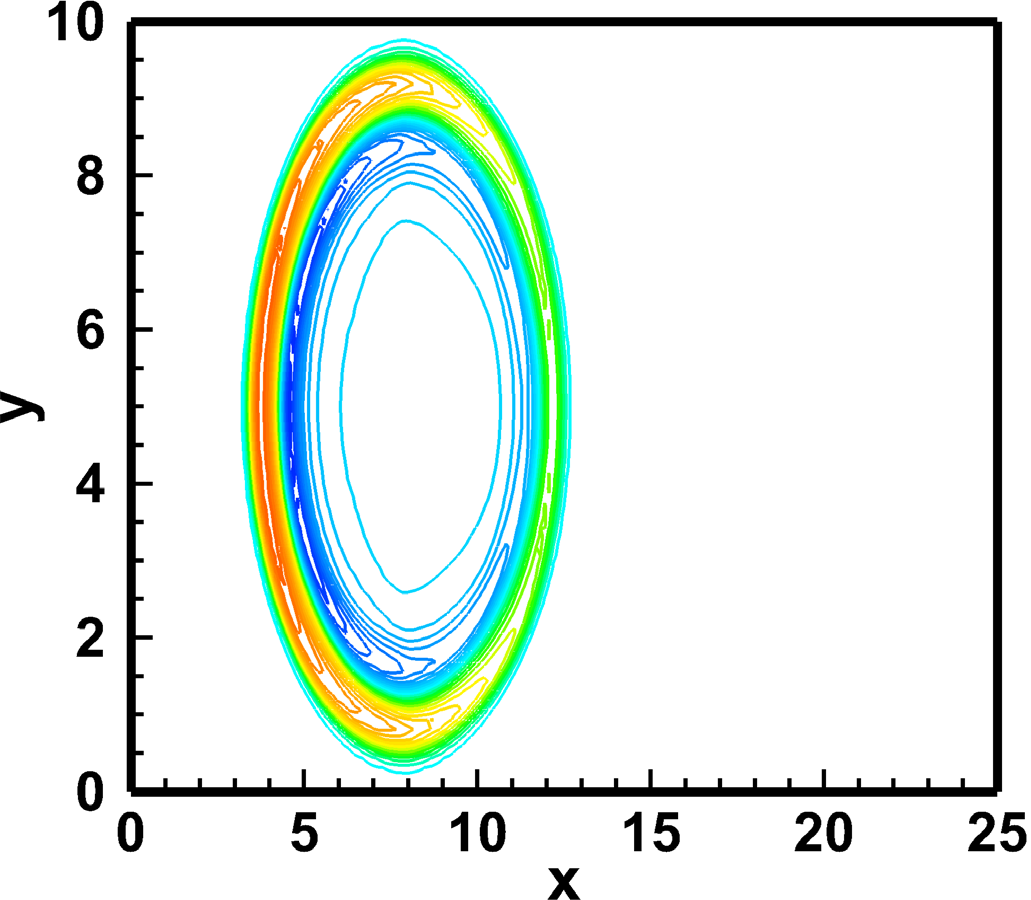}
\includegraphics[width=0.36\textwidth]{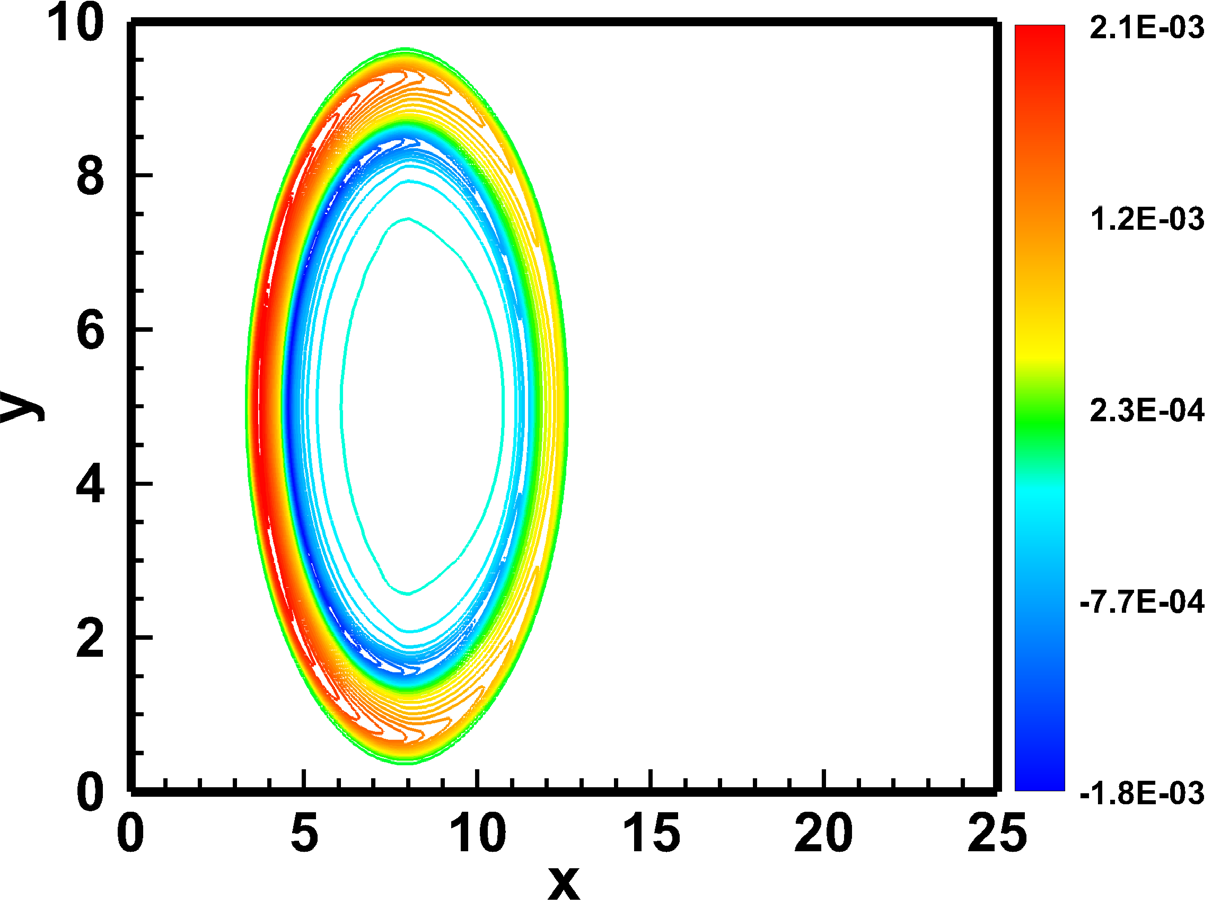}
\vskip-1mm
\caption{Example \ref{ex612}, subcritical case: Contours plots of $h(x,y,0.4)-h_{\rm eq}(x,y)$ computed by the PCCU-5 scheme on the
$100\times40$ (left), $200\times80$ (middle), $400\times160$ (right) meshes: $33$ uniformly spaced contours in the range $[-0.0017,0.0021]$.
}
\label{fig612}
\end{figure}
\begin{figure}[!htb]
\centering
\includegraphics[width=0.31\textwidth]{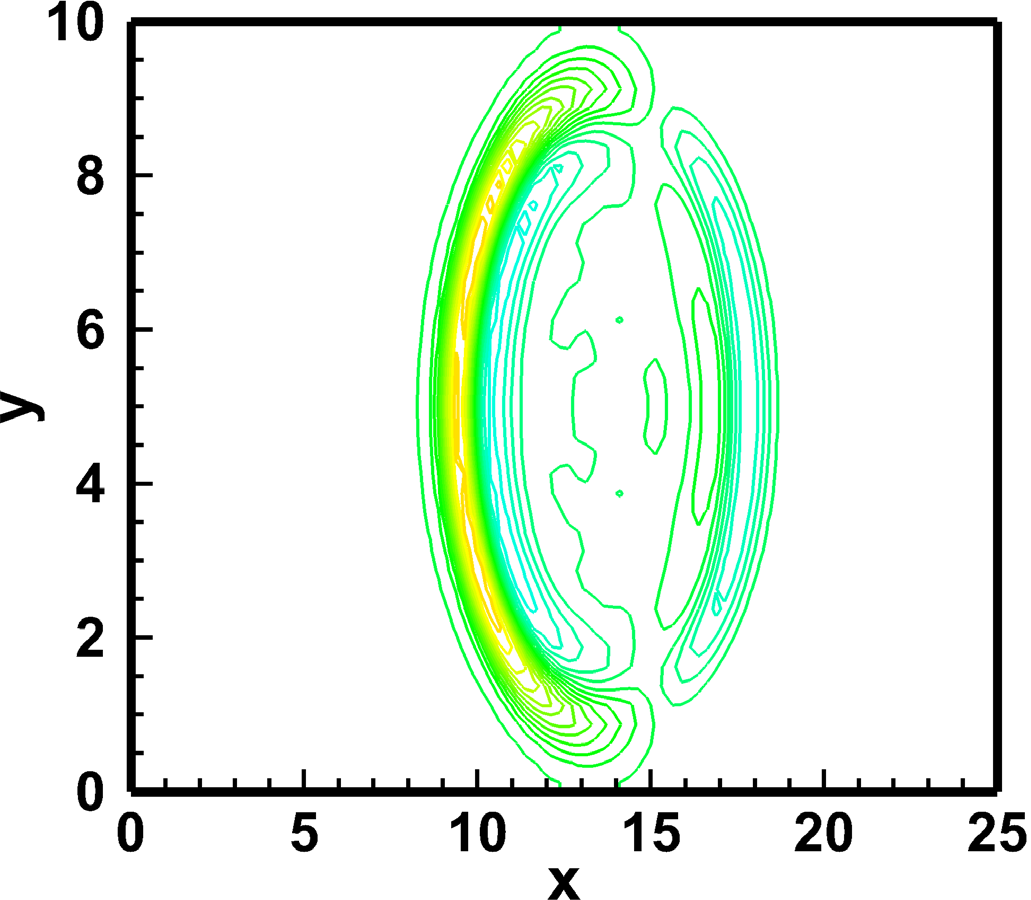}
\includegraphics[width=0.31\textwidth]{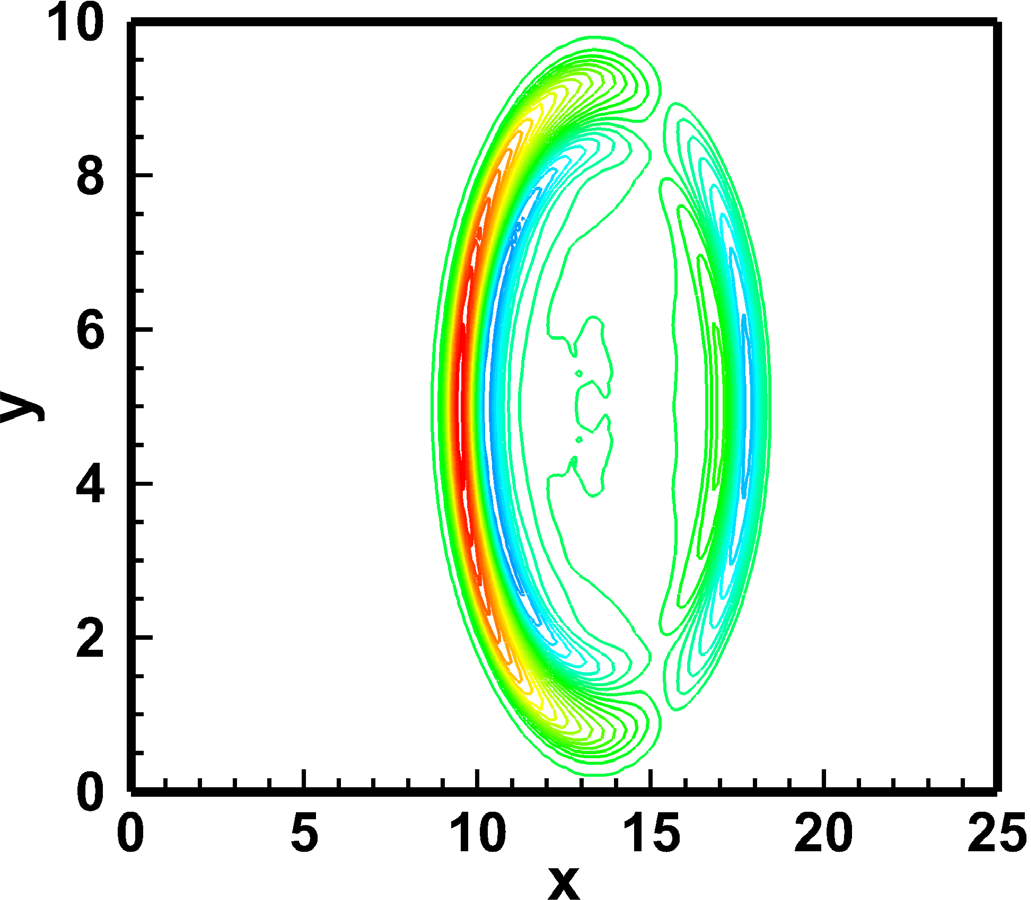}
\includegraphics[width=0.36\textwidth]{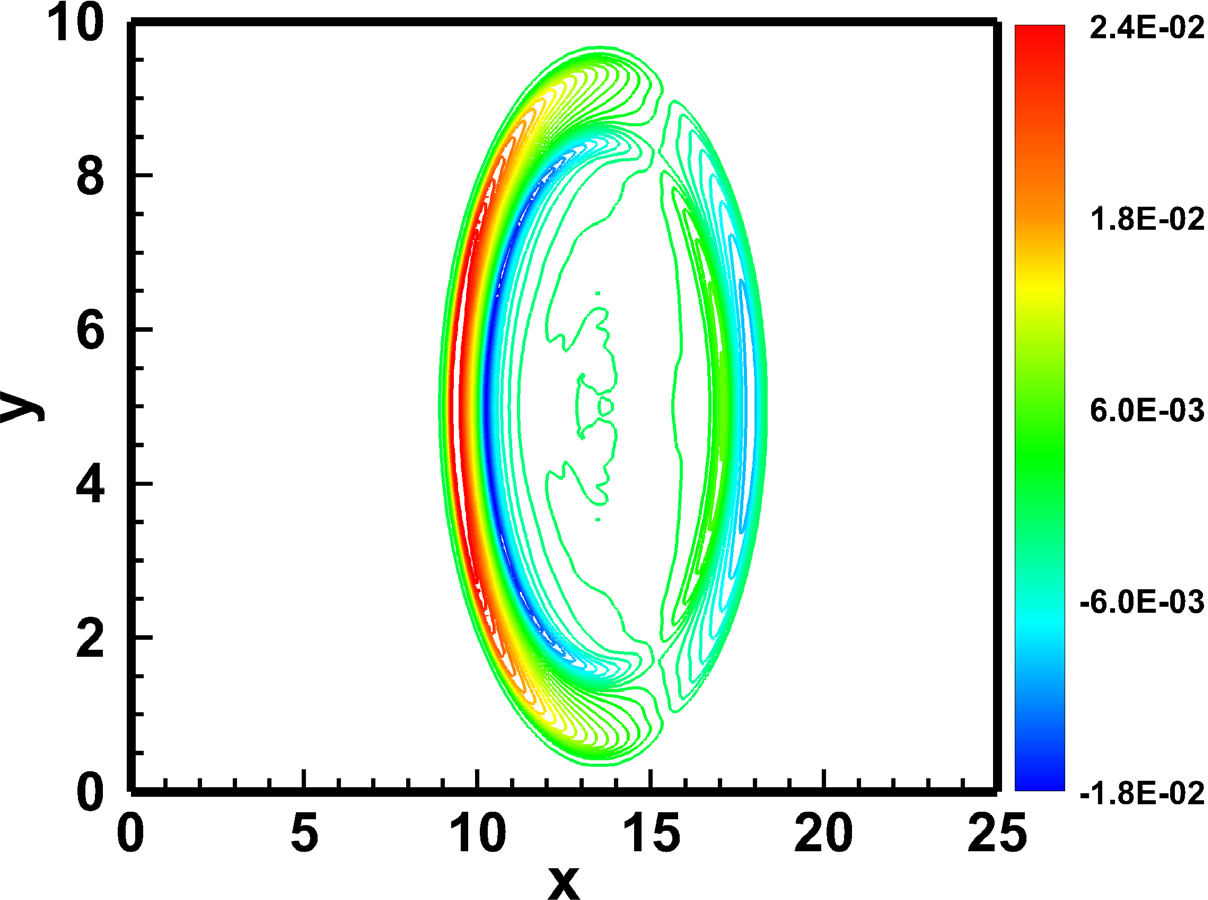}
\vskip-1mm
\caption{Example \ref{ex612}, supercritical case: Contours plots of $h(x,y,0.4)-h_{\rm eq}(x,y)$ computed by the PCCU-5 scheme on the
$100\times40$ (left), $200\times80$ (middle), $400\times160$ (right) meshes: $33$ uniformly spaced contours in the range $[-0.018,0.024]$.}
\label{fig613}
\end{figure}
\end{example}

\begin{example}[2-D discrete steady states]\label{ex613}
\rm The example, which is a modification of the test studied in \cite{CheKur}, aims to check the capability of the proposed PCCU-5 scheme to
converge to genuine 2-D equilibria and to capture their small perturbations.

We consider the bottom topography containing a Gaussian hump,
\begin{equation*}
Z(x,y)=0.1\,e^{-15(x-1.25)^2-90(y-0.5)^2},
\end{equation*}
and the initial conditions,
\begin{equation*}
(h,q^x,q^y,\theta)\Big|_{(x,y,0)}=(1-Z(x,y),0,0,1+0.1\,e^{-100(x-1)^2-90(y-0.5)^2}),
\end{equation*}
prescribed in the computational domain $[0,2]\times[0,1]$ subject to the following two different sets of boundary conditions, which
correspond to different flow regimes:

\medskip
\noindent
$\bullet$ {\bf Subcritical flow}:
\begin{equation*}
h(2,y,t)=1,\quad q^x(0,y,t)=0.3,\quad q^y(0,y,t)=0,\quad\theta(2,y,t)=1;
\end{equation*}

\medskip
\noindent
$\bullet$ {\bf Supercritical flow}:
\begin{equation*}
h(0,y,t)=1,\quad q^x(0,y,t)=3,\quad q^y(0,y,t)=0,\quad\theta(0,y,t)=1.
\end{equation*}

\noindent
In both cases, the rest of the boundary conditions in the $x$-direction are free (zero-order extrapolation), and the solid wall boundary
conditions are imposed in the $y$-direction.

We compute the numerical solutions by the PCCU-5 scheme until the final times $t=750$ (for the subcritical flow) and $t=20$ (for the
supercritical flow) on a uniform mesh with $50\times25$ cells. In both cases, the numerical solutions approach the corresponding steady
states, denoted by $h_{\rm eq},q^x_{\rm eq},q^y_{\rm eq},(h\theta)_{\rm eq}$. In Figure \ref{fig614}, we show how the difference
$\|h(\cdot,\cdot,t)-h(\cdot,\cdot,t-\dt)\|_{L^1}$ decays in time: this clearly confirms convergence to the discrete steady states, which are
plotted in Figure \ref{fig615}.
\begin{figure}[!htb]
\centering
\includegraphics[width=0.38\textwidth]{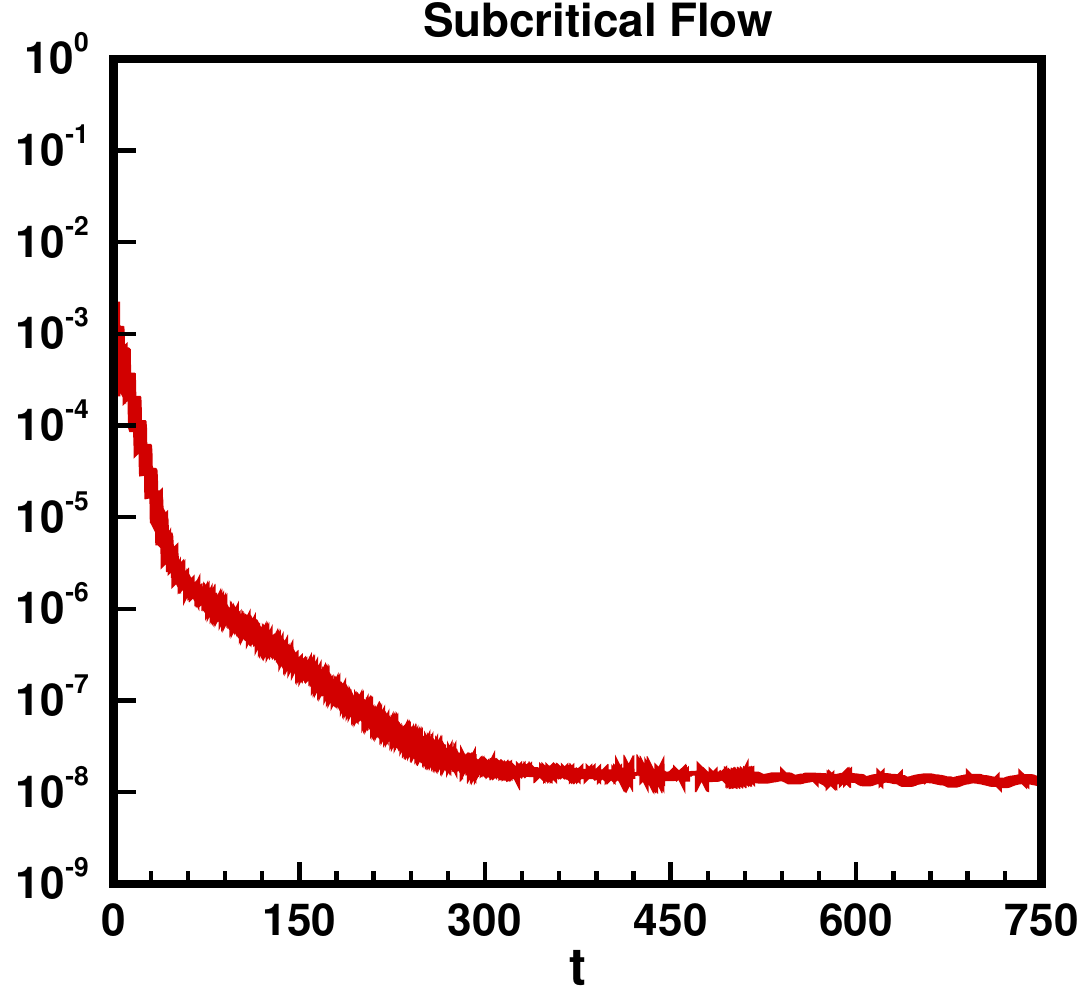}\hspace*{1cm}
\includegraphics[width=0.38\textwidth]{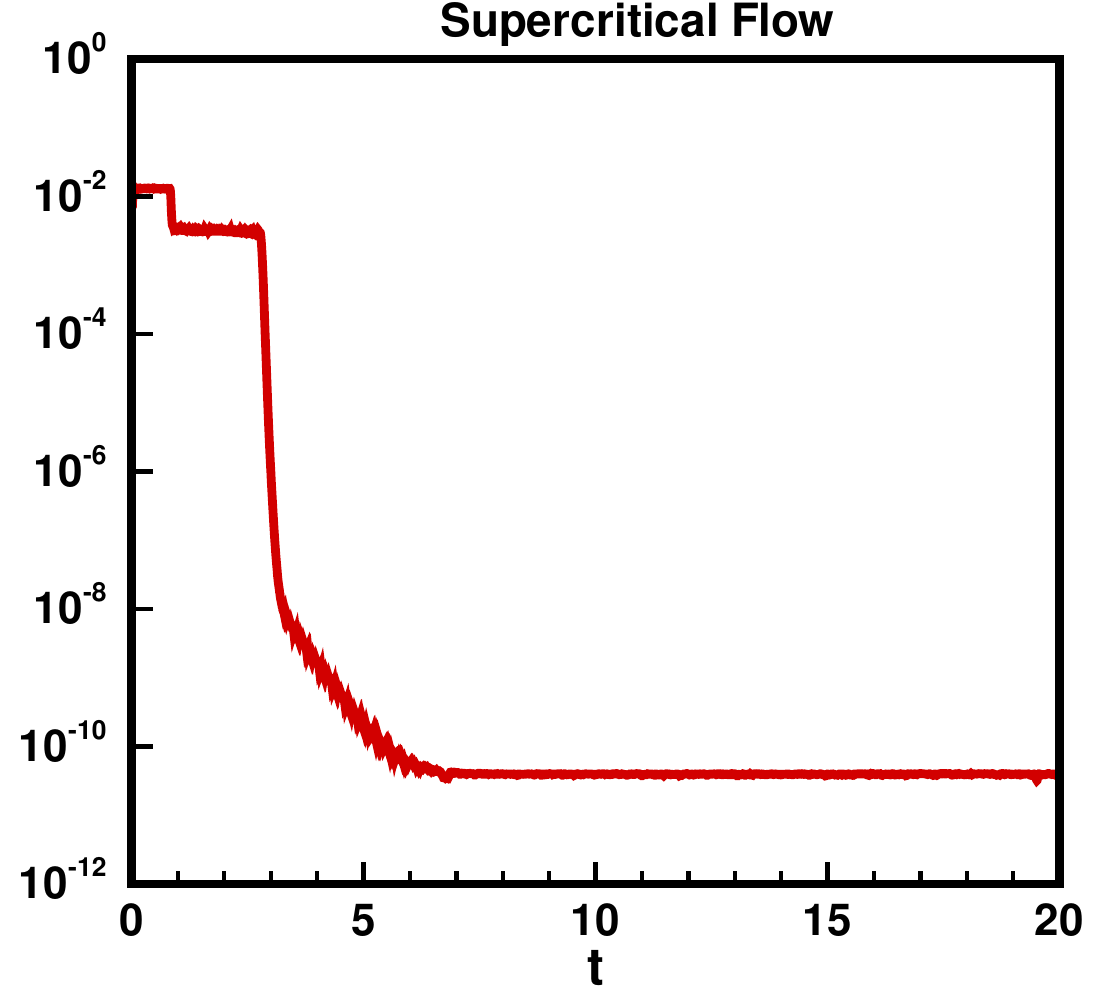}
\vskip-1mm
\caption{Example \ref{ex613}: $\|h(\cdot,\cdot,t)-h(\cdot,\cdot,t-\dt)\|_{L^1}$ (in a logarithmic scale) as a function of $t$ for the
subcritical (left) and supercritical (right) flows.}
\label{fig614}
\end{figure}
\begin{figure}[!htb]
\centering
\includegraphics[width=0.43\textwidth]{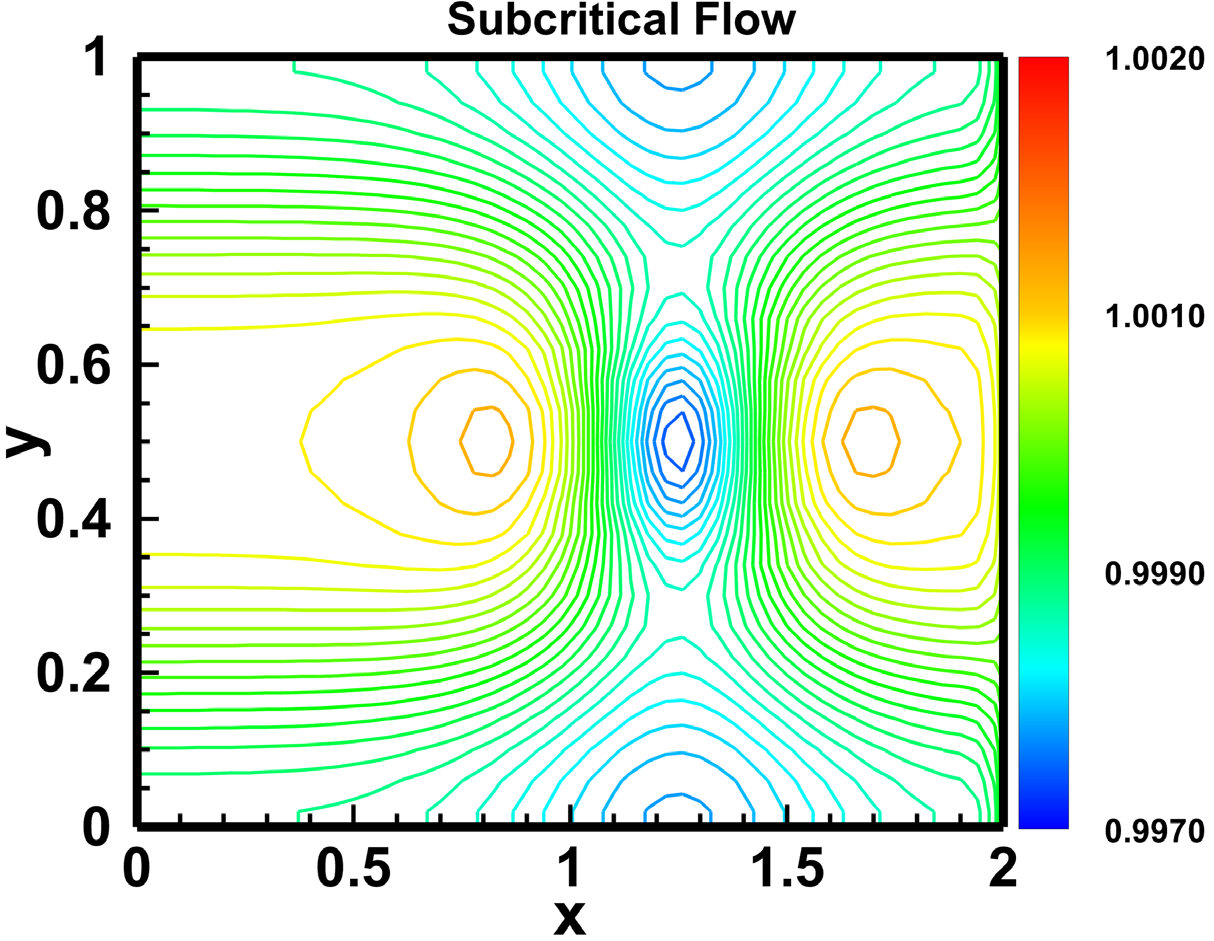}\hspace*{0.8cm}
\includegraphics[width=0.43\textwidth]{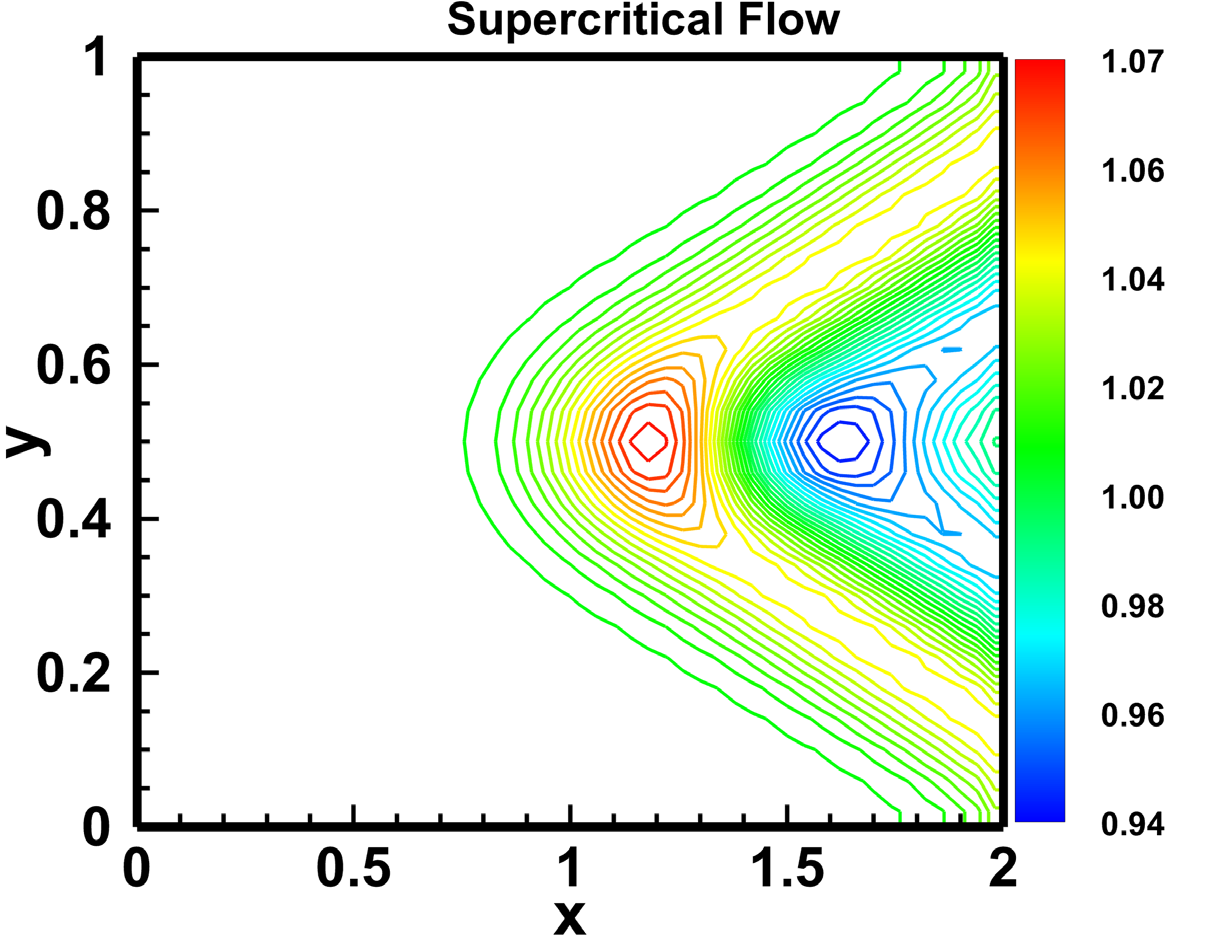}
\vskip-1mm
\caption{Example \ref{ex613}: Discrete subcritical (left) and supercritical (right) steady water surfaces $h_{\rm eq}+Z$.}
\label{fig615}
\end{figure}

We then add small circular perturbations to the water depth and take the following initial data:
\begin{equation*}
\begin{aligned}
&h(x,y,0)=h_{\rm eq}(x,y)+\left\{\begin{aligned}
&0.01,&&(x-0.25)^2+(y-0.55)^2<0.75^2,\\
&0,&&\mbox{otherwise};
\end{aligned}\right.\\
&q^x(x,y,0)=q^x_{\rm eq}(x,y),\quad q^y(x,y,0)=q^y_{\rm eq}(x,y),\quad(h\theta)(x,y,0)=(h\theta)_{\rm eq}(x,y).
\end{aligned}
\end{equation*}

We compute the solutions by the PCCU-5 scheme until the final time $t=0.28$ on a uniform mesh with $50\times25$ cells. In Figure
\ref{fig616}, we present the captured perturbations $h(x,y,0.28)-h_{\rm eq}(x,y)$, which are clearly oscillation-free.
\begin{figure}[!htb]
\centering
\includegraphics[width=0.43\textwidth]{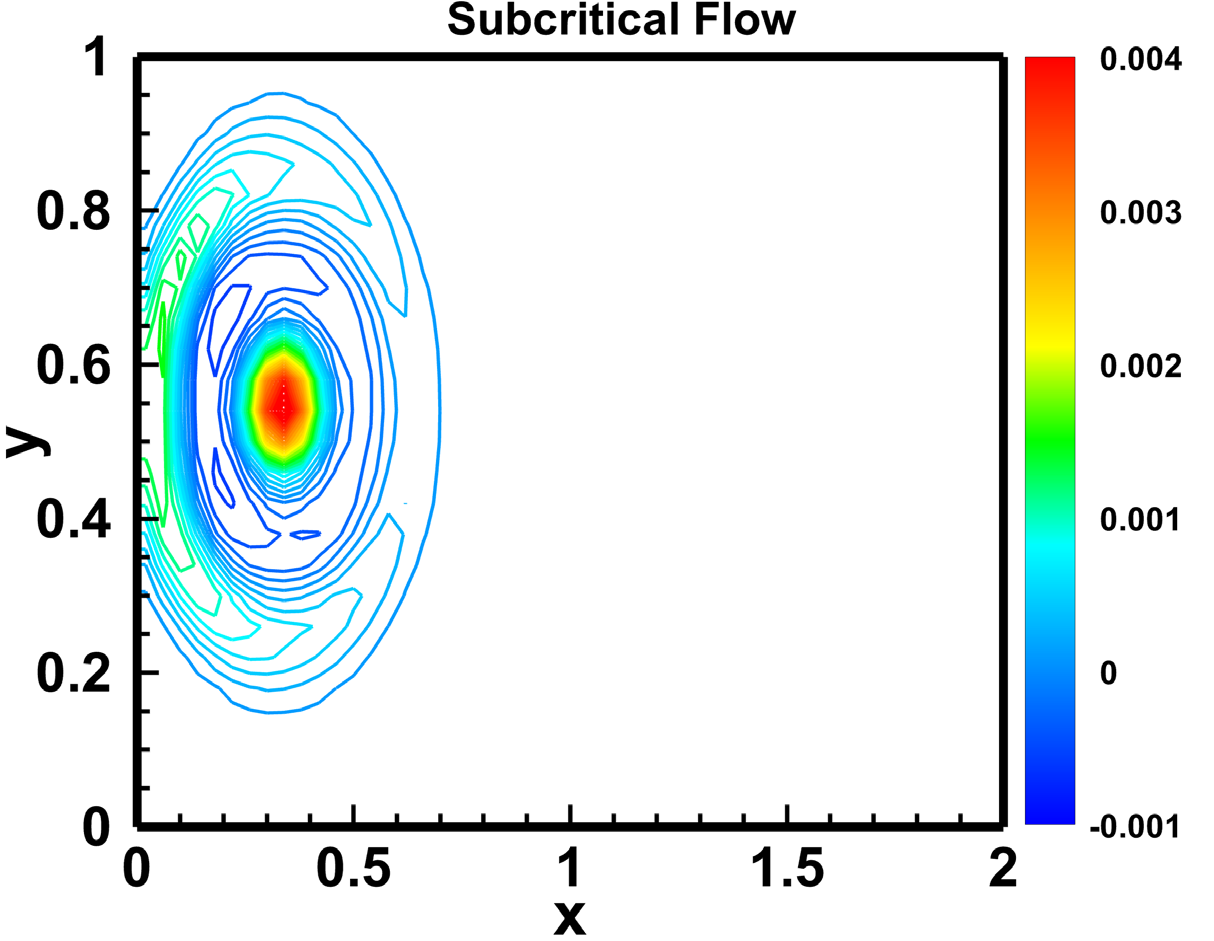}\hspace*{0.8cm}
\includegraphics[width=0.43\textwidth]{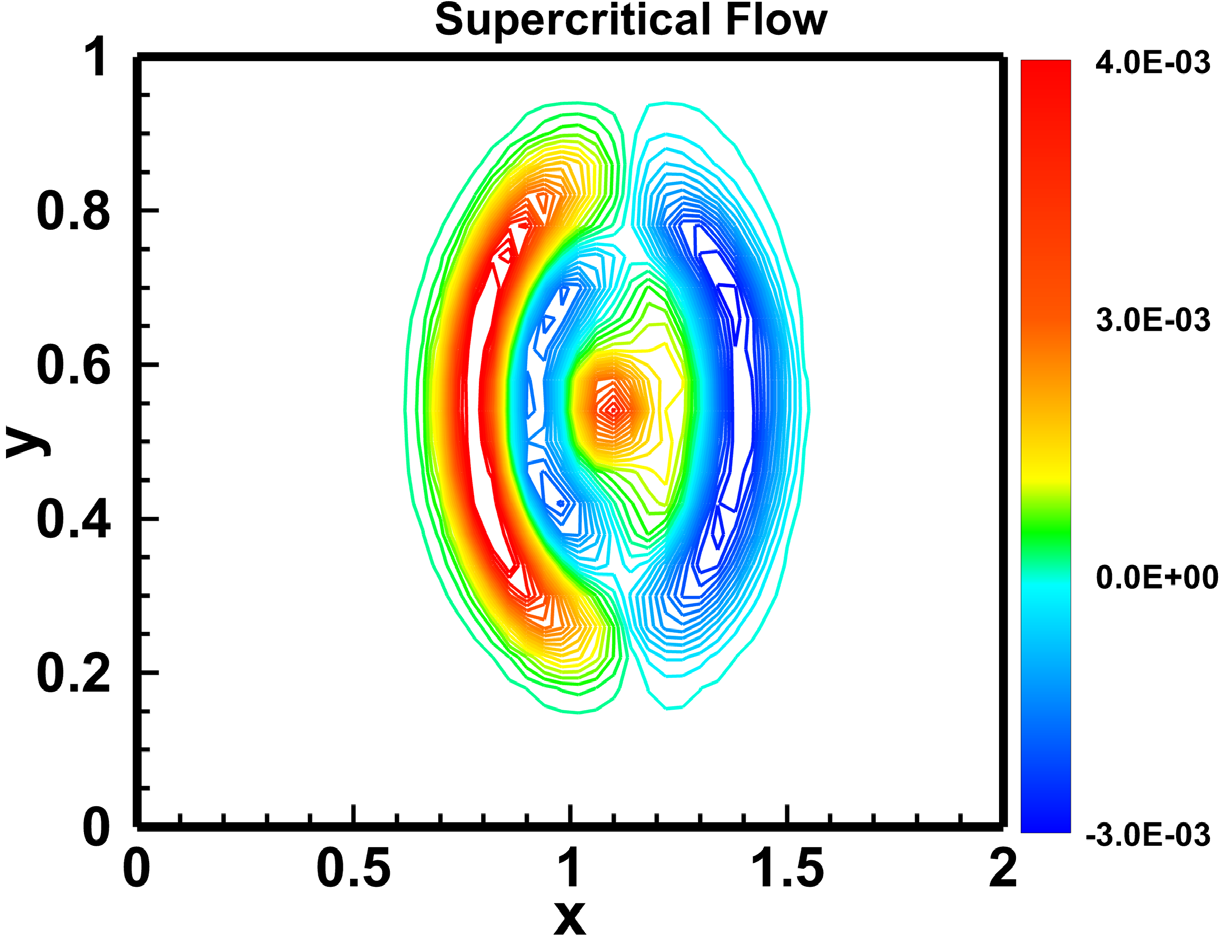}
\vskip-1mm
\caption{Example \ref{ex613}: Small perturbation ($h(x,y,0.28)-h_{\rm eq}(x,y)$) of the subcritical (left) and supercritical (right)
equilibria.}
\label{fig616}
\end{figure}
\end{example}

\section{Conclusions}\label{sec7}
In this work, we have developed a fifth-order well-balanced (WB) path-conservative A-WENO scheme with the central-upwind numerical fluxes
(PCCU-5) for the one- and two-dimensional Ripa models. A WB property has been enforced within the flux globalization framework. A
non-oscillatory nature of the resulting WB scheme is guaranteed with the help of a WENO interpolation applied to the local characteristic
equilibrium variables. The proposed PCCU-5 scheme has been tested on a large number of numerical examples. We have demonstrated that the
expected fifth order of accuracy is achieved, that the PCCU-5 scheme is capable of preserving a wide variety of steady states as well as
accurately capturing their small perturbations, and that it can also handle large magnitude shock waves.

\bigskip
\noindent {\bf Acknowledgments.} The work of Y.-P. Qiu and Z. Gao was partially supported by NSFC grant 12371435, Taishan Scholars Program
(No. tsqn202211059) and Shandong Provincial Natural Science Foundation (No. ZR2023MA043). The work of A. Kurganov was supported in part by
NSFC grant W2431004. The work of B.-S. Wang was supported in part by NSFC grant 12301530, the Shandong Provincial Qingchuang Science and
Technology Project, and the startup funding provided by the Ocean University of China. The work of X. Wen was partially supported by
Shandong Provincial Natural Science Foundation (No. ZR2022QA002).

\appendix
\section{Solution of \eref{3.7} and \eref{3.8}}\label{appA}
We describe a solution algorithm for obtaining $h_\jph^+$ from the ``$+$'' equation in \eref{3.7} (the ``$-$'' equation in \eref{3.7} as
well as the equations in \eref{3.8} can be solved similarly).

We rewrite the studied equation as
\begin{equation}
f\big(h_\jph^+\big):=\big(h_\jph^+\big)^3+a_0\big(h_\jph^+\big)^2+a_2=0,
\label{A.1}
\end{equation}
where
\begin{equation}
a_0=Z_\jph^++\frac{{\cal Q}_\jph^+-\Ep_\jph^+}{\theta_\jph^+},\quad a_2=\frac{\big(q_\jph^+\big)^2}{2\theta_\jph^+}.
\label{A.2}
\end{equation}
Equation \eref{A.1}--\eref{A.2} is a cubic equation, and we solve it exactly using the method outlined in \cite{Cheng19}.

We first note that for $h_\jph^+>0$,  $a_2\ge0$ and $a_0<0$, and we observe that positive solution of \eref{A.1}--\eref{A.2} exists only if
\begin{equation}
a_2 \leq -\frac{4}{27}a_0^3.
\label{A.3}
\end{equation}
If \eref{A.3} is not satisfied (which is a very unlikely case), we apply the Ai-WENO-Z interpolation to the water surface field
$w:=h+Z$ and then set
\begin{equation}
h_\jph^+=w_\jph^+-Z_\jph^+.
\label{A.4}
\end{equation}
If \eref{A.3} is satisfied, then there are two possibilities. If $q_\jph^+=0$, the unique positive solution is $h_\jph^+=-a_0$, while if
$q_\jph^+\ne0$, the three roots of the cubic equation \eref{A.1}--\eref{A.2} are given by
\begin{equation}
h_\jph^+=-\frac{1}{3}a_0\left[2\cos\left(\frac{\phi+2\pi k}{3}\right)+1\right],\quad k=1,2,3,
\label{A.5}
\end{equation}
where
\begin{equation}
\cos\phi=1+\frac{27a_2}{2a_0^3},\quad\phi\in(0,\pi).
\label{A.6}
\end{equation}
Out of these three roots, one is negative and two are positive (the latter ones correspond to subcritical and supercritical cases), and we
identify the physically relevant solution by choosing the positive root in \eref{A.5}--\eref{A.6}, which is closer to the value $h_\jph^+$
in \eref{A.4}.

\end{document}